\documentclass{amsart}
\usepackage{amsmath}
\usepackage{amsthm}
\usepackage{amssymb}
\usepackage{makecell} 
\usepackage{appendix}
\usepackage[german,english]{babel}
\usepackage[utf8]{inputenc}
\usepackage[new]{old-arrows}
\usepackage{epsfig,fancybox,color}
\usepackage{setspace}
\usepackage[linesnumbered,ruled]{algorithm2e}
\usepackage{tikz}
\usepackage{tikz-cd}
\usetikzlibrary{cd}
\usepackage[toc,page]{}
\usepackage{adjustbox}
\usepackage{todonotes}
\usepackage{yhmath}
\usepackage{mathtools}
\usepackage{stmaryrd}
\usepackage{fancyvrb} \VerbatimFootnotes
\usepackage{hyperref}
\usepackage{wasysym} 
\definecolor{darkred}{RGB}{160,0,0}
\definecolor{darkblue}{RGB}{0,0,160}
\hypersetup{
  colorlinks=true,
  citecolor=darkgreen,
  citecolor=darkblue,
  filecolor=black,
  linkcolor=darkblue,
  urlcolor=darkblue
}

\usepackage{xr}
\usepackage[margin=1.25in]{geometry} 

\tikzcdset{scale cd/.style={every label/.append style={scale=#1},
    cells={nodes={scale=#1}}}}

\DeclareMathAlphabet\mathbfcal{OMS}{cmsy}{b}{n} 

\DeclareFontFamily{U}{mathx}{\hyphenchar\font45}
\DeclareFontShape{U}{mathx}{m}{n}{
      <5> <6> <7> <8> <9> <10>
      <10.95> <12> <14.4> <17.28> <20.74> <24.88>
      mathx10
      }{}
\DeclareSymbolFont{mathx}{U}{mathx}{m}{n}
\DeclareMathSymbol{\bigtimes}{1}{mathx}{"91}

\newcommand{\excise}[1]{}

\theoremstyle{plain}
\newtheorem{thm}{Theorem}[section]
\newtheorem{lem}[thm]{Lemma}

\newtheorem{cor}[thm]{Corollary}
\newtheorem{prop}[thm]{Proposition}

\newtheorem{introthm}{Theorem}[section]

\theoremstyle{definition}
\newtheorem{defn}[thm]{Definition}

\newtheorem{prob}[thm]{Problem}

\newtheorem{construction}[thm]{Construction}
\newtheorem{condition}[thm]{Condition}

\theoremstyle{remark}

\newtheorem{remark}[thm]{Remark}

\newtheorem{notation}[thm]{Notation}

\numberwithin{equation}{section}

\renewcommand{\>}{\rangle}
\newcommand{\<}{\langle}
\newcommand{\CC}{\mathbb{C}}
\newcommand{\NN}{\mathbb{N}}

\newcommand{\A}{\mathbb{A}}

\newcommand{\maps}{\rightarrow}

\DeclareMathOperator{\Ban}{\mathbf{Ban}} 

\DeclareMathOperator{\E}{\mathbf{E}} 

\DeclareMathOperator{\LH}{\mathbf{LH}} 

\DeclareMathOperator{\TT}{\mathbb{T}} 

\DeclareMathOperator{\FarguesFontaine}{FF} 

\DeclareMathOperator{\act}{a} 

\DeclareMathOperator{\sh}{sh} 

\DeclareMathOperator{\D}{\mathbf{D}} 
\DeclareMathOperator{\Ho}{H} 

\DeclareMathOperator{\IndBan}{\mathbf{IndBan}} %

\DeclareMathOperator{\id}{id} 

\DeclareMathOperator{\Sp}{Sp} 
\DeclareMathOperator{\Hom}{Hom} 
\DeclareMathOperator{\intHom}{\underline{\Hom}} 
\DeclareMathOperator{\shHom}{\underline{\mathcal{H}\kern -.5pt\textit{om}}} 
\DeclareMathOperator{\RshHom}{R\underline{\mathcal{H}\kern -.5pt\textit{om}}} 
\DeclareMathOperator{\Ext}{Ext} 
\DeclareMathOperator{\shExt}{\underline{\mathcal{E}\kern -1,75pt\textit{xt}}} 

\DeclareMathOperator{\psh}{psh} 

\DeclareMathOperator{\Sh}{Sh} 
\DeclareMathOperator{\Loc}{Loc} 

\DeclareMathOperator{\I}{\mathbf{I}} 
\DeclareMathOperator{\C}{\mathbf{C}} 

\DeclareMathOperator{\op}{op} 

\DeclareMathOperator{\Dcap}{\mathcal{\wideparen{D}}} 
\DeclareMathOperator{\Dcapbd}{\mathcal{\wideparen{D}}\kern 1,00pt^{b}} 

\DeclareMathOperator{\Spa}{Spa} 
\DeclareMathOperator{\dR}{dR} 
\DeclareMathOperator{\pdR}{pdR} 
\DeclareMathOperator{\proet}{\text{pro\'{e}t}} %
\DeclareMathOperator{\et}{\text{\'{e}t}} 

\DeclareMathOperator{\ad}{ad} 

\DeclareMathOperator{\BB}{\mathbb{B}}
\DeclareMathOperator{\OB}{\cal{O}\kern -1,00pt\mathbb{B}} 
\DeclareMathOperator{\normalOB}{O\kern -1,00pt B} 
\DeclareMathOperator{\normalOA}{O\kern -1,00pt A} 
\DeclareMathOperator{\normalOtildeA}{O\kern -1,00pt {\tilde{A}}} 
\DeclareMathOperator{\OA}{\cal{O}\kern -1,00pt\mathbb{A}} 
\DeclareMathOperator{\OBcap}{\wideparen{\cal{O}\kern -1,00pt\mathbb{B}}} 
\DeclareMathOperator{\OBdR}{\cal{O}\kern -1,00pt\mathbb{B}_{\dR}} 
\DeclareMathOperator{\OXBdR}{\cal{O}_{X}\kern -1,00pt\mathbb{B}_{\dR}} 
\DeclareMathOperator{\OXxXBdR}{\cal{O}_{X\times X}\kern -1,00pt\mathbb{B}_{\dR}} 
\DeclareMathOperator{\OC}{\cal{O}\kern -0,90pt\mathbb{C}} %
\DeclareMathOperator{\CB}{\cal{C}\kern -1,00pt\mathbb{B}} 

\DeclareMathOperator{\Otheta}{\cal{O}\kern -1,00pt\theta} 
\DeclareMathOperator{\Ovartheta}{\cal{O}\kern -1,00pt\vartheta} 
\DeclareMathOperator{\Oiota}{\cal{O}\kern -1,00pt\iota} 

\DeclareMathOperator{\Der}{Der} 
\DeclareMathOperator{\TB}{\mathcal{T}\kern -1,00pt\mathbb{B}} 
\DeclareMathOperator{\OmegaB}{\Omega\kern -1,00pt\mathbb{B}} 

\DeclareMathOperator{\DcapB}{\Dcap\kern -1,00pt\BB}

\DeclareMathOperator{\RB}{R\kern -1,00pt\BB}
\DeclareMathOperator{\PB}{P\kern -1,00pt\BB}
\DeclareMathOperator{\rhoB}{\rho\kern -1,00pt\BB}

\DeclareMathOperator{\SB}{S\kern -1,50pt\BB} 

\DeclareMathOperator{\dRB}{\dR\kern -1,50pt\BB}
\DeclareMathOperator{\DDB}{\DD\kern -1,50pt\BB}

\DeclareMathOperator{\SolB}{\Sol\kern -0.5pt\BB} 
\DeclareMathOperator{\dRfunctorB}{\dRfunctor\kern -1.5pt\BB} 

\DeclareMathOperator{\Conn}{Conn} 
\DeclareMathOperator{\LocB}{\Loc\kern -1,50pt\BB} 
\DeclareMathOperator{\ConnOB}{\Conn\kern -1,50pt\OB} 

\DeclareMathOperator{\DB}{\mathcal{D}\kern -1,00pt\mathbb{B}} 
\DeclareMathOperator{\EB}{\mathcal{E}\kern -0.40pt\mathbb{B}} 

\DeclareMathOperator{\perf}{perf} 
\DeclareMathOperator{\bd}{b} 
\DeclareMathOperator{\Sol}{\mathcal{S}\kern -1.0pt\textit{ol}} 
\DeclareMathOperator{\nSol}{Sol} 
\DeclareMathOperator{\nSolB}{Sol\kern -0,5pt\mathbb{B}} 
\DeclareMathOperator{\Rec}{\mathcal{R}\kern -1.0pt\textit{ec}} 
\DeclareMathOperator{\nRec}{Rec} 

\DeclareMathOperator{\R}{R} 

\newcommand{\heart}{\ensuremath\heartsuit}

\DeclareMathOperator{\DD}{\mathbb{D}} 

\DeclareMathOperator{\cO}{\mathcal{O}} 

\DeclareMathOperator{\Ind}{\mathbf{Ind}} 

\newcommand{\cal}[1]{\mathcal{#1}}

\DeclareMathOperator{\ZZ}{\mathbb{Z}} 

\DeclareMathOperator{\holim}{holim} 
\DeclareMathOperator{\hocolim}{hocolim} 

\DeclareMathOperator{\rL}{L} 

\DeclareMathOperator{\qc}{qc} 

\DeclareMathOperator{\PD}{PD} 

\DeclareMathOperator{\isomap}{\stackrel{\cong}{\longrightarrow}} 

\DeclareMathOperator{\RH}{\mathcal{RH}} 
\DeclareMathOperator{\dRfunctor}{\textit{d}\kern +0.5pt\mathcal{R}} 

\DeclareMathOperator{\lh}{^{\small\heart}\kern -3.00pt}

\begin{document}

\selectlanguage{english}

\title{A fully faithful $p$-adic Riemann-Hilbert functor for coadmissible $\Dcap$-modules}

\author{Finn Wiersig}
\address{National University of Singapore}
\email{fwiersig@nus.edu.sg}
\date{\today}

\begin{abstract}
  This article establishes a Riemann--Hilbert correspondence in rigid-analytic geometry.
  We construct an explicit solution functor and prove that it
  is fully faithful on Ardakov--Wadsley's coadmissible $\Dcap$-modules.
  For vector bundles with flat connection, our functor is canonically
  identified with Scholze's horizontal sections functor.
\end{abstract}

\maketitle


\tableofcontents



\section{Introduction}
\label{ch:intro}


\subsection{Background}
\label{subsec:Background-recII}

The Riemann--Hilbert correspondence, originating from Hilbert's 21st problem~\cite{Equationsdifferentiellesapointssinguliersreguliers,KashiwaraKawaiholonomicIII,Kashiwaraconstructibilityholonomic,KashiwaraRHforholonomicsystems,Mebkhoutuneequivalence}, has had profound applications---from Saito's mixed Hodge modules~\cite{MR1000123,MR1047415} to the proofs of the Kazhdan--Lusztig conjecture~\cite{MR632980,MR610137}.
In its most standard version, independently proven by
Kashiwara~\cite{Kashiwaraconstructibilityholonomic,KashiwaraRHforholonomicsystems}
and Mebkhout~\cite{Mebkhoutuneequivalence}, the Riemann-Hilbert correspondence
asserts that the solution functor
is an equivalence of two triangulated categories associated to a complex manifold $Y$: the opposite of the
bounded derived category of regular holonomic modules over the sheaf of finite-order differential operators
$\mathcal{D}$, and the bounded derived category of constructible sheaves of complex vector spaces.
Its proof relies on the six-functor formalism for holonomic $\cal{D}$-modules.

Fix a prime $p$. The goal of this article is to construct a
$p$-adic analogue of the Riemann--Hilbert correspondence in
rigid-analytic geometry.

Let $k$ be a complete discretely valued field of mixed characteristic
$(0,p)$ with perfect residue field, and let $X$ be a smooth
rigid-analytic $k$-variety. The most direct analogue of the classical
theory would be a correspondence for holonomic modules over the sheaf
$\mathcal D$ of finite-order differential operators on $X$, obtained
by adapting the methods of Kashiwara and Mebkhout.

This approach encounters an immediate obstruction in the $p$-adic
setting. The classical proof relies on the six-functor formalism for
holonomic $\mathcal D$-modules, but no comparable formalism is
available for holonomic $\mathcal D$-modules on rigid-analytic spaces.
The problem is already visible on the disc. For example,
if $j\colon \DD^{\times}\hookrightarrow \DD=\Sp k\<z\>$ is the inclusion of the punctured
disc $\DD^{\times}=\DD\setminus\{0\}$ into the disc, then the $\cal{D}$-module pushforward of the structure sheaf $j_{+}\cO$
is not coherent as a $\cal{D}$-module on $\DD$. In particular, it cannot be holonomic.

One obtains better functoriality properties by working instead with a larger sheaf
of infinite-order differential operators. This idea belongs to the framework
of arithmetic $\mathcal{D}$-modules developed by Berthelot~\cite{BerthelotDmodulesArithmeticI} and further studied by
Huyghe with collaborators
\cite{HuygheDdaggerAffiniteProjectif,NootHuygheBB,HuygheSchmidtStrauchArithmeticStructures}.
In the rigid-analytic setting used here, we work with the sheaf
$\cal{D}\subset\Dcap$ constructed by Ardakov--Wadsley~\cite{AW19,MR3846550}.

  Concretely, on the unit disc one has
  \begin{equation*}
    \Dcap\left(\DD\right)=\varprojlim_{r>0}\cal{D}_{r}\left(\DD\right),
    \text{ where }
    \cal{D}_{r}\left(\DD\right)=\left\{
      \sum_{\alpha\geq0}f_{\alpha}\partial_{z}^{\alpha}\in\prod_{\alpha\geq0}k\<z\>\partial_{z}^{\alpha}\colon
      \text{$f_{\alpha}/p^{r\alpha}\to0$ as $\alpha\to\infty$}
    \right\}.
  \end{equation*}
  Equivalently, $\Dcap\left(\DD\right)$ consists of infinite-order diferential operators
  $\sum_{\alpha\geq0}f_{\alpha}\partial_{z}^{\alpha}\in\prod_{\alpha\geq0}k\<z\>\partial^{\alpha}$
  whose coefficients satisfy a rapid decay condition:
  $f_{\alpha}/p^{r\alpha}\to0$ for every $r>0$.

Already in this basic example, $\Dcap\left(\DD\right)$
is not noetherian, and hence the category of coherent
$\Dcap$-modules is not the correct substitute for coherent $\cal{D}$-modules.
Ardakov--Wadsley showed in~\cite{AW19}, and Bode showed in greater
generality in~\cite{MR3991380}, that the sections of $\Dcap$ over affinoids are 
Fréchet--Stein algebras in the sense of Schneider--Teitelbaum~\cite{ST03}.
That is, they are projective limits of noetherian Banach algebras with flat transition maps.
In the case of the disc, the Fréchet--Stein presentation is precisely the projective
system $\Dcap\left(\DD\right)=\varprojlim_{r>0}\cal{D}_{r}\left(\DD\right)$ displayed above.

This allows one to use Schneider--Teitelbaum's notion of coadmissible
modules. In the model case of the disc, a left $\Dcap(\DD)$-module $M$
is coadmissible if every $M_{r}=\cal{D}_{r}\left(\DD\right)\otimes_{\Dcap\left(\DD\right)}M$
is finitely generated as a $\cal{D}_{r}\left(\DD\right)$-module,
and the canonical map $M\to\varprojlim_{r>0}M_{r}$ is an isomorphism.
The resulting category of coadmissible $\Dcap(\DD)$-modules is abelian.

Ardakov--Wadsley globalise this construction and define the category of
coadmissible $\Dcap$-modules on a smooth rigid-analytic variety. This
category has the expected local-to-global behaviour. For example,
taking global sections gives an equivalence between coadmissible
$\Dcap$-modules on the disc $\DD$ and coadmissible
$\Dcap(\DD)$-modules. It also remedies the basic functorial failure
discussed above: for the inclusion
$j\colon \DD^\times\hookrightarrow \DD$, the pushforward
$j_+\mathcal O$ is a coadmissible $\Dcap$-module on
$\DD$. Thus coadmissibility plays for $\Dcap$-modules the role that
coherence plays for ordinary $\mathcal D$-modules.

The goal of this article is therefore to prove a Riemann-Hilbert correspondence for
$\Dcap$-modules. The classical strategy of Kashiwara and Mebkhout,
however, is not presently available in this setting. Indeed, despite substantial efforts,
there is currently no six-functor formalism for holonomic $\cal{D}$-modules on rigid-analytic spaces~\cite{TowardsRHArdakov2015,BitounBode+2021+97+118,MR4332493,Bo21,BodeOberwolfach2022,hallopeau2024microlocalisationdoperateursdifferentielsarithmetiques,hallopeau2024widehatmathcald0mathfrakxkmathbbqmodulesholonomes}.

Other techniques are available. Our approach follows Morita-theoretic ideas
of Prosmans-Schneiders in the Archimedean setting~\cite{PS00}. These apply in much larger generality
than Kashiwara's and Mebkhout's classical argument. Indeed, we prove
a Riemann-Hilbert correspondence not just for the holonomic $\Dcap$-modules,
but for \emph{all} coadmissible $\Dcap$-modules.

Coadmissible $\Dcap$-modules are also closely connected with locally
analytic representation theory. In the sense of
Schneider--Teitelbaum~\cite{ST01,STla02,ST03,ST05}, admissibility of a
locally analytic representation is expressed by the coadmissibility of
its strong dual. Ardakov's Beilinson--Bernstein-type localisation
theorem~\cite{EquivariantD2021} relates this representation-theoretic
notion of coadmissibility to coadmissible $\Dcap$-modules.

Ardakov--Wadsley~\cite{AW2024_GlobalsectionsoneqlinebundlesonOmega}
used this connection in
their study of locally analytic representations arising from
Drinfeld's first covering of the $p$-adic upper half-plane. They prove
admissibility results for these representations and describe their
duals in terms of explicit equivariant differential equations. This
leads naturally to the problem of understanding the local solutions of
such equations. Theorem~\ref{introthm:easymain} provides a foundation
for this problem by constructing a period-valued solution functor for
coadmissible $\Dcap$-modules and proving that it is fully faithful.


\subsection{Main result}

Our Theorem~\ref{introthm:easymain} is most naturally stated in terms of derived categories. While one might consider the category $\D_{\perf}^{\bd}\left(\Dcap\right)$, it is insufficient for applications to $p$-adic representation theory, as it does not contain all coadmissible $\Dcap$-modules. Instead, we work with Bode's triangulated category $\D_{\cal{C}}\left(\Dcap\right)$ of $\mathcal{C}$-complexes~\cite{Bo21}, which includes $\cal{D}_{\perf}^{\bd}\left(\Dcap\right)$ as a full subcategory and carries a t-structure whose heart is the abelian category of coadmissible $\Dcap$-modules.

A $\Dcap$-module may be viewed as a system of $p$-adic differential equations.
To obtain a sufficiently rich solution theory for such systems, we need to incorporate periods into our formalism.
Therefore, we introduce $\BB_{\pdR}^{\dag}$, the \emph{overconvergent almost de Rham\footnote{The acronym $\pdR$ refers to the french \emph{presque de Rham}.} period sheaf}, on the pro-étale site
$X_{\proet}$ of $X$. The absolute period ring is $B_{\pdR}^{\dag}:=B_{\dR}^{\dag,+}[1/t,\log t]$,
where $B_{\dR}^{\dag,+}$ is the stalk $\cal{O}_{\FarguesFontaine,\infty}$
of the structure sheaf of the analytic Fargues-Fontaine curve at the point $\infty\in\FarguesFontaine$
and $t$ is the $p$-adic $2\pi i$.
In particular, $B_{\pdR}^{\dag}$ is a subring of
the almost de Rham period ring $B_{\pdR}=B_{\dR}^{+}[1/t,\log t]$
defined by Fontaine~\cite{Fontaine2004Arithmetic}.

To bridge $\Dcap$-modules and $\BB_{\pdR}^{\dag}$-modules, we define a $\nu^{-1}\Dcap$-$\BB_{\pdR}^{\dag}$-bimodule $\OB_{\pdR}^{\dag}$ on $X_{\proet}$, where $\nu: X_{\proet} \to X$ is the canonical projection.

\begin{introthm}[Theorem~\ref{thm:pdRdagsolfunctor-fullyfaithful-reconstructionpaper} and Corollary~\ref{cor:pdRdagdRfunctor-fullyfaithful-reconstructionpaper}]
\label{introthm:easymain}
  The solution functor
  \footnote{See \S\ref{subsec:solutionderHamfunctorsDcapmod-solutionpaper} for the precise definition.}
  \begin{equation*}
    \SolB_{\pdR}^{\dag}\colon\D_{\cal{C}}\left(\Dcap\right)^{\op}\hookrightarrow\D\left(\BB_{\pdR}^{\dag}\right),\quad
    \cal{M}^{\bullet}\mapsto\R\shHom_{\nu^{-1}\Dcap}\left(\nu^{-1}\cal{M}^{\bullet},\OB_{\pdR}^{\dag}\right)
  \end{equation*}
  is a fully faithful embedding of triangulated categories. Equivalently, the de Rham functor
  \begin{equation*}
    \dRfunctorB_{\pdR}^{\dag}\colon\D_{\cal{C}}\left(\Dcap\right)\hookrightarrow\D\left(\BB_{\pdR}^{\dag}\right),\quad
    \cal{M}^{\bullet}\mapsto\SolB_{\pdR}^{\dag}\left(\DD\left(\cal{M}^{\bullet}\right)\right)[\dim X]
  \end{equation*}
  is a fully faithful embedding of triangulated categories.
\end{introthm}

Here, $\DD$ denotes duality functor introduced in~\cite{Bo21}.

\begin{remark}
  If $\cal{M}^{\bullet}$ is a vector bundle with flat connection, then
  $\dRfunctorB_{\pdR}^{\dag}\left(\cal{M}^{\bullet}\right)$
  is canonically identified with Scholze's horizontal sections functor~\cite[Theorem 7.6]{Sch13pAdicHodge}.
  We prove this claim in~\cite{WiersigCauchy},
  and explain that it is a generalisation of the classical $p$-adic Cauchy theorem.
\end{remark}


\subsection{Proof of Theorem~\ref{introthm:easymain}}

As explained in~\S\ref{subsec:Background-recII}, Theorem~\ref{introthm:easymain}
is analogous to the reconstruction theorem of Prosmans--Schneiders~\cite{PS00}
in the Archimedean setting. If $Y$ is a complex analytic manifold, they prove that
the restriction of the solution functor to bounded perfect complexes
\begin{equation*}
  \D_{\perf}^{\bd}\left(\mathcal{D}^{\infty}\right)^{\op}
  \hookrightarrow
  \D\left(\CC_{Y}\right),
  \qquad
  \mathcal{M}^{\bullet}
  \mapsto
  \R\shHom_{\mathcal{D}^{\infty}}\left(\mathcal{M}^{\bullet},\mathcal{O}\right),
\end{equation*}
is a fully faithful embedding of triangulated categories. Here
$\mathcal{D}^{\infty}$ denotes the sheaf of infinite-order differential operators
on $Y$, which plays the role of $\Dcap$ in the $p$-adic setting; moreover,
$\CC_{Y}$ is the constant sheaf and $\mathcal{O}$ is the sheaf of holomorphic
functions on $Y$.

The proof of Prosmans--Schneiders does not rely on a six-functor formalism.
Instead, it uses Morita-theoretic techniques to reduce the reconstruction theorem
to the following computation: the
$\mathcal{D}^{\infty}$-$\CC_{Y}$-bimodule structure on $\mathcal{O}$ induces an
isomorphism
\begin{equation*}
  \rho\colon
  \mathcal{D}^{\infty}
  \isomap
  \R\shHom_{\CC_{Y}}\left(\mathcal{O},\mathcal{O}\right).
\end{equation*}
In cohomological degree zero this had already been proved by Ishimura~\cite{Ishimura78},
namely one has an isomorphism
\begin{equation*}
  \varrho\colon
  \mathcal{D}^{\infty}
  \isomap
  \shHom_{\CC_{Y}}\left(\mathcal{O},\mathcal{O}\right).
\end{equation*}

Our proof follows the same strategy. We reduce Theorem~\ref{introthm:easymain}
to the following explicit computation.

\begin{introthm}[Corollary~\ref{cor:padic-ishimuraprosmansschneiders}]
\label{introthm:Ext-vanishing}
  The $\nu^{-1}\Dcap$-$\BB_{\pdR}^{\dag}$-bimodule structure on
  $\OB_{\pdR}^{\dag}$ induces an isomorphism
  \begin{equation*}
    \rho\colon
    \Dcap
    \isomap
    \R\nu_{*}
    \R\shHom_{\BB_{\pdR}^{\dag}}
    \left(\OB_{\pdR}^{\dag},\OB_{\pdR}^{\dag}\right).
  \end{equation*}
\end{introthm}

Theorem~\ref{introthm:Ext-vanishing} implies
Theorem~\ref{introthm:easymain} for bounded perfect complexes of
$\Dcap$-modules. For coadmissible $\Dcap$-modules and $\mathcal{C}$-complexes,
one needs a slight strengthening of Theorem~\ref{introthm:Ext-vanishing};
we refer to the beginning of
\S\ref{subsec:proofsofmainthms-reconstructionpaper} for the precise formulation.
For ease of exposition, we only explain here the proof of
Theorem~\ref{introthm:Ext-vanishing}.

The proof has two parts. First, we prove that $\rho$ is an isomorphism in
cohomological degree zero. Second, we prove that the right-hand side has no
higher cohomology.

For the degree-zero statement, we have to show that the
$\nu^{-1}\Dcap$-$\BB_{\pdR}^{\dag}$-bimodule structure on $\OB_{\pdR}^{\dag}$
induces an isomorphism
\begin{equation*}
  \varrho\colon
  \Dcap
  \isomap
  \nu_{*}
  \shHom_{\BB_{\pdR}^{\dag}}
  \left(\OB_{\pdR}^{\dag},\OB_{\pdR}^{\dag}\right).
\end{equation*}
This assertion is local on $X$, so we may assume that $X$ admits étale coordinates
$z_{1},\dots,z_{d}$. By
\cite[Theorem~\ref*{cor:localdescription--OBpdRdag-sheafy-reconstructionpaper}]{WiersigPeriods},
after passing to a suitable pro-étale cover $\widetilde{X}\to X$, the period
structure sheaf admits the explicit description
\begin{equation}\label{eq:OBpdRlocal-intrp-recpaper}
  \varinjlim_{q\geq0}
  \BB_{\pdR}^{\dag}|_{\widetilde{X}}
  \left\<\frac{Z_{1},\dots,Z_{d}}{p^{q}}\right\>
  \cong
  \OB_{\pdR}^{\dag}|_{\widetilde{X}},
\end{equation}
for suitable coordinates $Z_{1},\dots,Z_{d}$. Under this description, $\rho$
sends $\partial_{z_{i}}\in\Dcap(X)$ to $d/dZ_{i}$. The injectivity of $\rho$
is then immediate.

For surjectivity, let
$\psi\colon\OB_{\pdR}^{\dag}\to\OB_{\pdR}^{\dag}$ be a
$\BB_{\pdR}^{\dag}$-linear endomorphism. To $\psi$ one associates a formal differential operator
$\eta(\psi)$, in the same spirit as Ishimura's construction
in~\cite{Ishimura78}. The convergence condition in
\eqref{eq:OBpdRlocal-intrp-recpaper} implies that this Taylor series defines an
element $\eta(\psi)\in\Dcap(X)$. By construction, one then has
$\rho(\eta(\psi))=\psi$, which proves surjectivity.

It remains to prove that
\begin{equation*}
  \R\nu_{*}
  \R\shHom_{\BB_{\pdR}^{\dag}}
  \left(\OB_{\pdR}^{\dag},\OB_{\pdR}^{\dag}\right)
\end{equation*}
is concentrated in degree zero. The local description
\eqref{eq:OBpdRlocal-intrp-recpaper} gives effective control over the derived
homomorphism sheaves. The main point is therefore to control the derived
pushforward along $\nu$.

In Lemma~\ref{lem:prop:derivedRq-iso-positivedegree-reconstructionpaper}, we
reduce this to the computation of $\R\nu_{*}\OB_{\pdR}^{\dag}$. This computation
was carried out in
\cite[Theorem~\ref{thm:pushforwardOBpdRdag-introduction}]{WiersigPeriods}:
there is a canonical isomorphism
\begin{equation*}
  \mathcal{O}
  \isomap
  \R\nu_{*}\OB_{\pdR}^{\dag}.
\end{equation*}
This proves the vanishing of the higher cohomology sheaves on the right-hand
side of Theorem~\ref{introthm:Ext-vanishing}, and hence completes the proof of
Theorem~\ref{introthm:Ext-vanishing}.


\subsection{Structure of this article}

In \S\ref{sec:Recap-Periods-Rec}, we briefly recall the constructions
from~\cite{WiersigPeriods}: the period sheaves used in this article and the
main results from that work on which we rely. We do not review the definitions
in full, but only give the overview needed for the sequel.

The main body of the article begins in \S\ref{sec:sol-dR-functors-recpaper}. After
recalling some basic facts about $\Dcap$-modules in \S\ref{subsec:Dcap-modules}, we construct in
\S\ref{subsec:OBla-bimodule} the $\nu^{-1}\Dcap$-$\BB_{\pdR}^{\dag}$-bimodule
structure on $\OB_{\pdR}^{\dag}$. This leads to the solution and de Rham
functors, which are defined in
\S\ref{subsec:solutionderHamfunctorsDcapmod-solutionpaper}.

In \S\ref{ch:recthms-recpaper}, we study these functors. We first recall, in
\S\ref{subsec:C-recpaper}, the category of $\cal{C}$-complexes. We then
state the main results of the article in \S\ref{subsec:thms-recpaper}. Their
proofs occupy the following two sections: Theorem~\ref{introthm:Ext-vanishing}
is proved in \S\ref{subsec:proofsofmainthms-reconstructionpaper}, and
Theorem~\ref{introthm:easymain} is proved in
\S\ref{subsec:proofofsolfunctorfullyfaithfull-reconstructionpaper}.

Finally, in Appendix~\ref{appendix:reconstructionpaper}, we collect several
general results on monoidal categories, functional analysis, and sheaves which
are used in the main text.


\subsection{Acknowledgements}

Parts of this article are based on the author’s PhD thesis.
I am deeply grateful to my PhD supervisors, Konstantin Ardakov and Kobi Kremnitzer,
for introducing me to this subject and for their guidance and support:
thank you Konstantin and Kobi.
Furthermore, I would like to thank
Tomoyuki Abe,
Andreas Bode,
Hui Gao, and
David Hansen
for valuable conversations.
Part of this research was financially supported by a Mathematical Institute Award at Oxford University.


\subsection{Notations and conventions}\label{subsec:conventions-reconstructionpaper}
Throughout, $k$ denotes a complete discrete valuation field of mixed characteristic
$(0,p)$ with perfect residue field $\kappa$.
Fix a uniformiser $\pi\in k$ and write $k^{\circ}\subseteq k$
for the ring of power-bounded elements. Set $k_{0}:=W(\kappa)[1/p]$.
Let $C$ be the completion of a fixed algebraic closure $\overline{k}$ of $k$.

\begin{remark}
  We require $k$ to be discretely valued for two reasons. First, this
  hypothesis is used to obtain the local description
  \eqref{eq:OBpdRlocal-intrp-recpaper} of $\OB_{\pdR}^{\dag}$;
  see~\cite[Remark~\ref*{rem:whydiscvalued}]{WiersigPeriods}.
  Second, it allows us to apply
  \cite[Theorems~4.0.4 and~4.0.5]{BSSW2024_rationalizationoftheKnlocalsphere},
  which underlie the proof of
  Theorem~\ref{thm:BdRdagplusplus-Galoiscohomology-introduction}.
  We note that
  \cite{BSSW2024_rationalizationoftheKnlocalsphere} uses the term
  \emph{local field} in the sense specified at the beginning of its
  \S 4.
\end{remark}

$=$ denotes an equality, $\cong$ denotes an isomorphism, and $\simeq$ denotes an equivalence.
When an equality, isomorphism, or equivalence follows from a specific result, the reference
is written above the symbol denoting the equality, isomorphism, or equivalence. For example,
$X\stackrel{\text{\ref{lem:sheafinternalhom-over-monoid}}}{\cong}Y$ means:
Lemma~\ref{lem:sheafinternalhom-over-monoid} implies that $X$ and $Y$ are isomorphic.

The natural numbers are $\NN=\{0,1,2,3,\dots\}$.






\section{Overconvergent period sheaves}
\label{sec:Recap-Periods-Rec}

We recall the period sheaves from~\cite{Sch13pAdicHodge}
and~\cite{WiersigPeriods} which will be used throughout the article.

Let $X$ be a smooth locally Noetherian adic space over $\Spa\left(k,k^{\circ}\right)$.
The pro-étale site $X_{\proet}$ of $X$ consists of formal limits
$\text{``}\varprojlim\text{"}_{i\in I}U_{i}$ of $U_{i}\in X_{\et}$
such that $I$ is a cofiltered category and that the
transition maps $U_{j}\maps U_{i}$ are finite étale and surjective.
A collection $\{f_{i}\colon U_{i}\maps U\}_{i}$ in $X_{\proet}$ is
a covering if it is a pointwise covering, and a set-theoretic condition
is satisfied which we ignore for simplicity.
There is an obvious canonical projection of sites $\nu\colon X_{\proet}\to X$.

Let us spell out the basic example which will be used repeatedly. Let
\begin{equation*}
  X=\TT^{d}=\Spa\left(k\left\<T_{d}^{\pm},\dots,T_{d}^{\pm}\right\>,k^{\circ}\left\<T_{1}^{\pm},\dots,T_{d}^{\pm}\right\>\right)
\end{equation*}
be the $d$-dimensional torus. It has the family of étale coverings
\begin{equation*}
  \TT_{j}^{d}:=\Spa\left(k\left\<T_{1}^{\pm1/p^{j}},\dots,T_{1}^{\pm1/p^{j}}\right\>,
  k^{\circ}\left\<T_{1}^{\pm1/p^{j}},\dots,T_{1}^{\pm1/p^{j}}\right\>\right).
\end{equation*}
whose formal limit
$\widetilde{\TT}^{d}:=\text{``}\varprojlim_{j\in\NN}\text{"}\TT_{j}^{d}$
is an element of the site $\widetilde{\TT}_{\proet}^{d}$.
If $K$ is the completion of an algebraic extension of $k$ such that $K$ is perfectoid,
then the base-change $\widetilde{\TT}_{K}^{d}$ is an affinoid perfectoid object of
$\TT^{d}_{\proet}$. More generally, for arbitrary $X$, the
affinoid perfectoid objects form a basis for $X_{\proet}$.

In~\cite{WiersigPeriods}, we introduce several period sheaves on $X_{\proet}$.
The most important ones for the present article are
$\BB_{\dR}^{\dag,+}$, $\BB_{\dR}^{\dag}$, and $\BB_{\pdR}^{\dag}$.
They are filtered colimits of sheaves of Banach spaces. For instance,
\begin{equation*}
  \BB_{\dR}^{\dag,+}
  =
  \varinjlim_{q\geq0}\BB_{\dR}^{q,+}
  =
  \varinjlim_{q\geq0}\BB_{\dR}^{>q,+}.
\end{equation*}
In order to keep track of the topologies induced by these Banach structures,
we regard these period sheaves as sheaves with values in the category of
ind-Banach spaces. This allows us to use the formalism of~\cite{Sch99}
to obtain a robust formalism for the cohomology of ind-Banach valued sheaves.
We use this formalism freely throughout the article.

We also use the corresponding period structure sheaves
$\OB_{\dR}^{\dag,+}$,
$\OB_{\dR}^{\dag}$, and
$\OB_{\pdR}^{\dag}$.
By definition, there is a canonical morphism of sheaves of rings
\begin{equation*}
  \nu^{-1}\cal{O}\widehat{\otimes}_{k_{0}}\BB_{\dR}^{\dag,+}
  \to
  \OB_{\dR}^{\dag,+},
\end{equation*}
and $\OB_{\dR}^{\dag,+}$ is obtained as a completion of
$\nu^{-1}\cal{O}\widehat{\otimes}_{k_{0}}\BB_{\dR}^{\dag,+}$. Although this
completion is somewhat involved, it admits a simple local description.

\begin{thm}[{\cite[Theorem~\ref*{thm:localdescription-of-OBla}]{WiersigPeriods}}]
  Let $X$ be affinoid and equipped with an étale map $X\to\TT^{d}$. Let
  $\widetilde{X}:=X\times_{\TT^{d}}\widetilde{\TT}^{d}$
  be the induced pro-étale cover of $X$. Then the morphism
  \begin{equation*}
    \BB_{\dR}^{\dag,+}|_{\widetilde{X}}
    \left\<\frac{Z_{1},\dots,Z_{d}}{p^{\infty}}\right\>
    \stackrel{\cong}{\longrightarrow}
    \OB_{\dR}^{\dag,+}|_{\widetilde{X}},
    Z_{i}\mapsto T_{i}-\left[T_{i}^{\flat}\right]
  \end{equation*}
  is an isomorphism of sheaves of
  $\BB_{\dR}^{\dag,+}|_{\widetilde{X}}$-ind-Banach algebras.
\end{thm}

The analogous local descriptions hold for $\OB_{\dR}^{\dag}$ and
$\OB_{\pdR}^{\dag}$. The key result from~\cite{WiersigPeriods} used in this
article is the canonical isomorphism
\begin{equation*}
  \cal{O}\isomap\R\nu_{*}\OB_{\pdR}^{\dag}.
\end{equation*}
This isomorphism plays a crucial role in the proof of
Theorem~\ref{introthm:easymain}.


\section{Solution and de Rham functors for $\Dcap$-modules}
\label{sec:sol-dR-functors-recpaper}




\subsection{The sheaf $\Dcap$}
\label{subsec:Dcap-modules}

Let $X$ denote an arbitrary smooth rigid-analytic $k$-variety.
We recall the construction of the sheaf $\Dcap$ on $X$.
Our main reference is~\cite{AW19}. $X_{w}$ denotes the category
whose objects are the affinoid subdomains of $X$ and whose morphisms
are the inclusions, carrying the weak Grothendieck topology.
\begin{equation*}
  \TT^{d}=\Sp k\left\<T_{1}^{\pm},\dots,T_{d}^{\pm}\right\>
\end{equation*}
is the $d$-dimensional torus over $k$.
First, we assume that $X=\Sp A$ is affinoid and equipped with an étale
morphism $g\colon X\to\TT^{d}$.
Compute the $A$-module of $k$-linear differentials of $A$:
\begin{equation*}
  L:=\Der_{k}(A)=\bigoplus_{l=1}^{d}A\partial_{l},
\end{equation*}
where $\partial_{l}$ denotes the lift of
the canonical vector field $d/dT_{l}$
along the étale map
$g^{\#}\colon\mathcal{O}\left( \TT^{d} \right)=k\left\<T_{1}^{\pm},\dots,T_{d}^{\pm}\right\>\to\mathcal{O}\left( X \right)$.
These $\partial_{l}$ do not need to preserve
$\mathcal{A}:=\mathcal{O}(X)^{\circ}$ in general. But because they
are bounded, see~\cite[Lemma 3.1]{KBB18},
we can find $r_{g}\geq 1$ large enough such that
the $\pi^{r}\partial_{l}$ preserve $\mathcal{A}$
for every $r\geq r_{g}$. Define the $\mathcal{A}$-submodule
\begin{equation*}
  \mathcal{L}:=\bigoplus_{l=1}^{d}\mathcal{A}\partial_{l}\subseteq L.
\end{equation*}
$\mathcal{L}_{r}:=p^{r}\mathcal{L}$ is an \emph{$\mathcal{A}$-Lie lattice}
for every $r\geq r_{g}$, that is
$[\mathcal{L}_{r},\mathcal{L}_{r}]\subseteq\mathcal{L}_{r}$
and $\mathcal{L}_{r}(\mathcal{A})\subseteq\mathcal{A}$.

\begin{remark}
  A priori, it seems more natural to define $\mathcal{L}_{r}$ to be $\pi^{r}\mathcal{L}$.
  However, we choose $p$ over $\pi$ because it simplifies
  the proof of Theorem~\ref{thm:OBlaplus-bimodule}.
\end{remark}

We cite the following from~\cite[\S 3]{AW19}.
Fix $r\geq r_{g}$. Denote the pullback $\mathcal{O}(X)\to\mathcal{O}(U)$ by
$\omega$. An \emph{admissible} $k^{\circ}$-algebra is a
commutative $k^{\circ}$-algebra which is topologically
of finite type and flat over $k^{\circ}$.
An \emph{affine formal model} in $\mathcal{O}(U)$
is an admissible $k^{\circ}$-algebra $\mathcal{B}$ such that
$\mathcal{O}(U)\cong \mathcal{B}\otimes_{k^{\circ}}k$. $\mathcal{B}$
is \emph{$\mathcal{L}_{r}$-stable} if $\omega(\mathcal{A})\subseteq\mathcal{B}$ and the
action of $\mathcal{L}_{r}$ on $\mathcal{A}$ lifts to $\mathcal{B}$.
$U$ is \emph{$\mathcal{L}_{r}$-admissible} if it admits an $\mathcal{L}_{r}$-stable affine
formal model. $X_{r}=X_{g,r}:=X_{w}\left(\mathcal{L}_{r}\right)$ denotes
the full subcategory of $X_{w}$ consisting of the $\mathcal{L}_{r}$-admissible
affinoid subdomains. It is a site by~\cite[\S 3.2, Lemma]{AW19}.
The coverings are the finite admissible coverings
by objects in $X_{r}$.


\begin{defn}{\cite[\S 3.3, Definition]{AW19}}
  Let $X$ be affinoid and equipped with an étale morphism $g\colon X\to\TT^{d}$.
  Fix $r\geq r_{g}$. For any $\mathcal{L}_{r}$-admissible affinoid subdomain
  $U\subseteq X$ and any $\mathcal{L}_{r}$-stable affine formal model $\mathcal{B}$ in
  $\mathcal{O}(U)$, define
  \begin{equation*}
    \mathcal{D}_{r} (U)= \mathcal{D}_{g,r}(U) := \widehat{U\left( \mathcal{B}\otimes_{\mathcal{A}}\mathcal{L}_{r}\right)}\otimes_{k^{\circ}}k.
  \end{equation*}
  The symbol $U$ refers to the enveloping algebra
  of $\mathcal{B}\otimes_{\mathcal{A}}\mathcal{L}_{r}$ as a
  $\left( k^{\circ} , \mathcal{B}\right)$-algebra, see~\cite[\S 2.1]{AW19}.
  The completion is the $\pi$-adic one.
\end{defn}

Regard $\mathcal{D}_{r} (U)$ as a $k$-Banach algebra with unit ball
$\widehat{U\left( \mathcal{B}\otimes_{\mathcal{A}}\mathcal{L}_{t}\right)}$.

\begin{lem}\label{lem:Dt-sheaf}
  Let $X$ be affinoid and equipped with an étale morphism $g\colon X\to\TT^{d}$.
  Fix $r\geq r_{g}$. The assignment
  \begin{equation*}
    U\mapsto\mathcal{D}_{r}(U)
  \end{equation*}
  defines a sheaf of $k$-Banach algebras on $X_{r}$.
\end{lem}

\begin{proof}
  See~\cite[\S 3.5, Theorem]{AW19}
  and apply the open mapping theorem.
\end{proof}

\begin{defn}{\cite[\S 9.3, Definition]{AW19}}\label{defn:Dcapsubsections}
  $X$ is affinoid and it admits an étale morphism $X\to\TT^{d}$.
  For every $U\in X_{w}$, define the $k$-ind-Banach algebra
  \begin{equation*}
    \Dcap\left(U\right):=
      \varprojlim\mathcal{D}_{g,r}(U),
  \end{equation*}
  cf. Lemma~\ref{lem:monoidobject-lim}.
  The inverse limits runs over all étale maps $g\colon X\to\TT^{d}$
  and $r$ large enough such that $U\in X_{r}$.
\end{defn}

\begin{lem}\label{lem:Dcap-sheaf-on-a-basis}
  $X$ is affinoid and admits an étale morphism $X\to\TT^{d}$. Then
  \begin{equation*}
    0
    \to\Dcap\left(X\right)
    \to\prod_{j}\Dcap\left(U_{j}\right)
    \to\prod_{j_{1},j_{2}}\Dcap\left(U_{j_{1}}\times_{U} U_{j_{2}}\right)
  \end{equation*}
  is strictly exact for any finite covering $\left\{U_{j}\to X\right\}_{j}$ by affinoids.
\end{lem}

\begin{proof}
  This follows from Lemma~\ref{lem:Dt-sheaf}.
\end{proof}

\begin{defn}
  Let $X$ be an arbitrary smooth rigid-analytic $k$-variety.
  By abuse of notation, $\Dcap$ denotes the unique sheaf
  of $k$-ind-Banach algebras such that for every
  $U\in X_{w}$ admitting an étale map $U\to\TT^{d}$, $\Dcap\left(U\right)$
  is given as in Definition~\ref{defn:Dcapsubsections}.
\end{defn}

\begin{lem}\label{lem:Dcap-frechetalgebra}
  Suppose $X$ is affinoid and equipped with a fixed étale morphism $g\colon X\to\TT^{d}$.
  Then the canonical morphism
    $\Dcap\left(U\right)
    \stackrel{\cong}{\longrightarrow}
    \varprojlim_{r\geq r_{g}}\mathcal{D}_{g,r}\left(U\right)$
  is an isomorphism for any affinoid subdomain $U\subseteq X$.
\end{lem}

\begin{proof}
  This follows from~\cite[\S 6.1, Lemma, (b)]{AW19}.
\end{proof}


\subsection{The bimodule structure on $\OB_{\dR}^{\dag,+}$}
\label{subsec:OBla-bimodule}

We continue to fix a smooth rigid-analytic $k$-variety $X$.
Every affinoid $k$-algebra is strongly Noetherian,
cf.~\cite[Theorem 3.1.8.3]{ScholzeWeinstein2020} and~\cite[\S 6.1.1, Proposition 3]{BGR84},
thus the adic space $X^{\ad}$ associated to $X$ is locally Noetherian.
This allows to define $X_{\proet}:=\left(X^{\ad}\right)_{\proet}$,
the \emph{pro-étale site of $X$}.
The morphism of sites $\nu$ is
\begin{equation*}
  X_{\proet}=\left(X^{\ad}\right)_{\proet}
  \stackrel{\nu}{\longrightarrow}
  \left(X^{\ad}\right)_{\et}
  \longrightarrow X.
\end{equation*}
See~\cite[\S 2.1]{Huber96Etale} for the definition
of the morphism at the right-hand side.

The term \emph{module} always refers to a left module.
We view $\nu^{-1}\Dcap$
as a sheaf of $k$-ind-Banach algebras
and $\nu^{-1}\mathcal{O}$ as a sheaf
of $\nu^{-1}\Dcap$-ind-Banach modules.
This makes $\nu^{-1}\mathcal{O}\widehat{\otimes}_{k_{0}}\BB_{\dR}^{q,+}$
a $\nu^{-1}\Dcap\widehat{\otimes}_{k_{0}}\BB_{\dR}^{q,+}$-module object,
for any $q\in\NN$.
Since $\BB_{\dR}^{q,+}$ is commutative, it is in fact a
$\nu^{-1}\Dcap$-$\BB_{\dR}^{q,+}$-bimodule object,
cf. Definition~\ref{defn:bimoduleobject}.

\begin{thm}\label{thm:OBlaplus-bimodule}
  Let $X$ be affinoid. Fix an étale morphism $g\colon X\to\TT^{d}$
  and $q\geq q_{g}$ (see~\cite[Notation~\ref*{notation:qgbound}]{WiersigPeriods}).
  Then there exists a unique
  $\nu^{-1}\Dcap$-$\BB_{\dR}^{q,+}$-bimodule structure on $\OB_{\dR}^{q,+}$
  such that the canonical morphism
  \begin{equation}\label{eq:Otensorbla-to-OBla-epi}
    \nu^{-1}\mathcal{O}\widehat{\otimes}_{k_{0}}\BB_{\dR}^{q,+}
    \to\OB_{\dR}^{q,+}
  \end{equation}  
  is a morphism of $\nu^{-1}\Dcap$-$\BB_{\dR}^{q,+}$-bimodule objects.
\end{thm}

\begin{remark}\label{remark:OBdRqplus-bimodulestructureproof-overview}
  Fix the notation as in Theorem~\ref{thm:OBlaplus-bimodule}.
  We have the elements
  $\partial_{1},\dots,\partial_{d}\in\Dcap(X)$. By~\cite[Theorem~\ref*{thm:localdescription-of-OBqplus}]{WiersigPeriods},
  $\OB_{\dR}^{q,+}$ has the sections $Z_{1},\dots,Z_{d}$, locally on the proétale site. Then
  \begin{equation*}
    \partial_{j}\cdot Z_{l} := \frac{d}{dZ_{j}}\left( Z_{l} \right)=\delta_{jl}
  \end{equation*}
  defines the action of $\nu^{-1}\Dcap$, where $\delta_{jl}$ is the Kronecker delta.
  A large part of the proof concerns showing that this action is compatible with
  the canonical $\BB_{\dR}^{q,+}$-module structure.
\end{remark}

\begin{proof}[Proof of Theorem~\ref{thm:OBlaplus-bimodule}]
  Corollary~\ref{cor:epi-uniquebimoduleaction}
  and the Lemma~\ref{lem:Otensorbla-to-OBla-epi}
  imply uniqueness.

  \begin{lem}\label{lem:Otensorbla-to-OBla-epi}
    The morphism~(\ref{eq:Otensorbla-to-OBla-epi}) is an epimorphism.
  \end{lem}

  \begin{proof}
    For every $q\in\NN$, $U=\text{``}\varprojlim_{i\in I}\text{"}U_{i}\in X_{\proet}$, and $i\in I$,
    \begin{equation*}
      \mathcal{O}^{+}\left(U_{i}\right)\widehat{\otimes}_{W(\kappa)}\A_{\dR}^{q}(U)
      \to
      \left(\mathcal{O}^{+}\left(U_{i}\right)\widehat{\otimes}_{W(\kappa)}\A_{\inf}(U)\right)
      \left\<\frac{\ker\Otheta_{\inf}}{p^{q}}\right\>
    \end{equation*}
    has dense image. It is thus an epimorphism.
    Now use that $\widehat{\otimes}_{k^{\circ}}k$ preserves colimits, then
    sheafify and apply~\cite[Lemma 2.16]{Bo21} to find
    that~(\ref{eq:Otensorbla-to-OBla-epi}) is an epimorphism.
  \end{proof}

  It remains to construct the bimodule structure.  
  Let $\widetilde{X}\to X$ be the standard pro-étale covering.
  Let $U\in X_{\proet}/\widetilde{X}$.
  One can assume that $U$ is affinoid perfectoid
  by~\cite[Lemma~\ref*{lem:sheaves-on-Xproet-and-Xproetaffperfdfin}]{WiersigPeriods}.
  We aim to give
  $\OB_{\dR}^{q,+}\left(U\right)$ the structure of a
  $\nu^{-1}\Dcap(U)$-$\BB_{\dR}^{q,+}(U)$-bimodule object,
  functorially in $U$.
  This is equivalent to giving $\OB_{\dR}^{q,+}\left(U\right)$
  the structure of a $\Dcap(V)$-$\BB_{\dR}^{q,+}(U)$-bimodule object
  for every $V\in X_{w}$ with $U\to\nu^{-1}\left(V\right)$,
  functorial in $U$ and $V$. Note
  \begin{equation*}
    \OB_{\dR}^{q,+}\left(U\right)
    \cong
    \BB_{\dR}^{q,+}\left( U \right)\left\< \frac{Z_{1},\dots,Z_{d}}{p^{\infty}}\right\>
  \end{equation*}
  by~\cite[Corollary~\ref*{cor:localdescription-of-subsections-OBqplus}]{WiersigPeriods}.
  Apply Lemma~\ref{lem:Dcap-frechetalgebra} and
  ~\cite[Lemmata~\ref*{lem:localisation-normedotimes}, and~\ref*{lem:localise-torsionfree-modules}]{WiersigPeriods}
  to get
  \begin{align*}
    \Dcap(V)
    &\cong\varprojlim_{r\geq r_{g}}\mathcal{D}_{g,r}\left( V\right)
    =\varprojlim_{r\geq r_{g}}
    \widehat{U\left( \mathcal{B} \otimes_{\mathcal{A}}\mathcal{L}_{r}\right)}
    \otimes_{k^{\circ}}k \\
    &\cong
    \varprojlim_{r\geq r_{g}}
    \widehat{U\left( \mathcal{B} \otimes_{\mathcal{A}}\mathcal{L}_{r}\right)}\left[ 1/p \right]
    \cong\varprojlim_{r\geq r_{g}}
    \widehat{U\left( \mathcal{B} \otimes_{\mathcal{A}}\mathcal{L}_{r}\right)}
    \widehat{\otimes}_{W(\kappa)}k_{0}.
  \end{align*}
  $\mathcal{A}$ is an affine formal model in $\mathcal{O}(X)$,
  $\mathcal{L}:=\bigoplus_{l=1}^{d}\mathcal{A}\partial_{l}$
  where $\partial_{l}$ denotes the lift of
  the canonical vector field $d/dT_{l}$
  along the étale map
  $g^{\#}\colon\mathcal{O}\left( \TT^{d} \right)=k\left\<T_{1}^{\pm},\dots,T_{d}^{\pm}\right\>\to\mathcal{O}\left( X \right)$,
  $\mathcal{L}_{r}:=p^{r}\mathcal{L}$,
  and $\mathcal{B}$ is an $\mathcal{L}_{r}$-stable affine formal model
  in $\mathcal{O}(V)$. 
  Assume that $r$ is large, such that $V\subseteq X$
  is $\mathcal{L}_{r}$-admissible. This is possible by~\cite[\S 6.1, Lemma (b)]{AW19}.
  We construct
  $\widehat{U\left( \mathcal{B} \otimes_{\mathcal{A}}\mathcal{L}_{r}\right)}$-$\BB_{\dR}^{q,+}(U)$-bimodule
  structures on $\BB_{\dR}^{q,+}\left( U \right)\left\< \frac{Z_{1},\dots,Z_{d}}{p^{q}}\right\>$
  for $r\geq q$ large enough such that the following holds.
  
  \begin{condition}\label{condition:gluethemaps-along-limcolim}
  \hfill
  \begin{itemize}
    \item[(i)] For all $r^{\prime}\geq r$, the following diagrams commute:
    \begin{equation*}
      \begin{tikzcd}
        \left( \widehat{U\left( \mathcal{B} \otimes_{\mathcal{A}}\mathcal{L}_{r}\right)}
        \widehat{\otimes}_{W(\kappa)} \BB_{\dR}^{q,+}(U)\right)
        \widehat{\otimes}_{W(\kappa)}\BB_{\dR}^{q,+}\left( U \right)\left\< \frac{Z_{1},\dots,Z_{d}}{p^{q}}\right\>
        \arrow{r} &
        \BB_{\dR}^{q,+}\left( U \right)\left\< \frac{Z_{1},\dots,Z_{d}}{p^{q}}\right\> \\
        \left( \widehat{U\left( \mathcal{B} \otimes_{\mathcal{A}}\mathcal{L}_{r^{\prime}}\right)}
        \widehat{\otimes}_{W(\kappa)} \BB_{\dR}^{q,+}(U)\right)
        \widehat{\otimes}_{W(\kappa)}\BB_{\dR}^{q,+}\left( U \right)\left\< \frac{Z_{1},\dots,Z_{d}}{p^{q}}\right\>
        \arrow{r}\arrow{u} &
        \BB_{\dR}^{q,+}\left( U \right)\left\< \frac{Z_{1},\dots,Z_{d}}{p^{q}}\right\>.\arrow[equal]{u}  
      \end{tikzcd}
    \end{equation*}
    Here, the horizontal maps denote the bimodule actions.
    \item[(ii)] The units
      $W(\kappa)\to\BB_{\dR}^{q,+}\left( U \right)\left\< \frac{Z_{1},\dots,Z_{d}}{p^{q}}\right\>$
      of all the bimodule structures on
      $\BB_{\dR}^{q,+}\left( U \right)\left\< \frac{Z_{1},\dots,Z_{d}}{p^{q}}\right\>$ coincide.
  \end{itemize}
  \end{condition}
  
  Once these bimodule structures are constructed,
  invert $p$
  to get a $\mathcal{D}_{g,r}(V)$-$\BB_{\dR}^{q,+}(U)$-bimodule structure
  on $\BB_{\dR}^{q,+}(U)\left\<\frac{Z_{1},\dots,Z_{d}}{p^{q}}\right\>$.
  Lemma~\ref{lem:bimodulestructure-limcolim} applies because of
  Condition~\ref{condition:gluethemaps-along-limcolim}, giving
  the desired $\Dcap(V)$-$\BB_{\dR}^{q,+}(U)$-bimodule
  structure on $\OB_{\dR}^{q,+}\left(U\right)$.
  
  Simplify notation:
    $O^{\circ}:=\mathcal{O}(V)^{\circ}$,
    $L_{r}:=\mathcal{B} \otimes_{\mathcal{A}}\mathcal{L}_{r}$,
    $B^{q,+}:=\BB_{\dR}^{q,+}\left(U\right)$, and
    $O\kern -1,00pt B^{q}:=\BB_{\dR}^{q,+}\left(U\right)\left\<\frac{Z_{1},\dots,Z_{d}}{p^{q}}\right\>$.
  We explained above that we aim to construct
  a $\widehat{U\left(L_{r}\right)}$-$B^{q,+}$-bimodule structure on
  $O\kern -1,00pt B^{q,+}$. This is equivalent to giving a morphism
  \begin{equation}\label{eq:Ulr-toHomAqOAqOAq}
    \widehat{U\left( L_{r}\right)}
    \to\intHom_{B^{q,+}}\left( O\kern -1,00pt B^{q,+}, O\kern -1,00pt B^{q,+}\right)
  \end{equation}
  of $W(\kappa)$-Banach algebras. We also explained that
  the morphism has to be suitably functorial in $V$ and $U$.
  This will be obvious from the construction. Thus we choose
  to omit both $U$ and $V$ in the simplified notation.
  
  View $L_{r}$ together with the canonical anchor map
  $\sigma\colon L_{r}\to\Der_{k^{\circ}}\left(O^{\circ}\right)$
  as a $\left(k^{\circ},O^{\circ}\right)$-Lie-algebra, cf.~\cite[\S 2.1]{AW19}.
  We construct~(\ref{eq:Ulr-toHomAqOAqOAq}) via the universal
  property of the completed enveloping algebra. Write
  \begin{align*}
    j_{O^{\circ}}\colon O^{\circ} &\to \intHom_{B^{q,+}}\left( O\kern -1,00pt B^{q,+}, O\kern -1,00pt B^{q,+}\right), \\
    f&\mapsto\left( h\mapsto \epsilon\left(f\right)h\right),
  \end{align*}
  where $\epsilon$ denotes the map constructed in~\cite[Lemma~\ref*{lem:imagetildenu-Ainf}]{WiersigPeriods},
  but we omit the index $i$. Next, we define
  \begin{align*}
    j_{L_{r}}\colon L_{r} &\to \intHom_{B^{q,+}}\left( O\kern -1,00pt B^{q,+}, O\kern -1,00pt B^{q,+}\right), \\
    \sum_{l=1}^{d}f_{l}p^{r}\partial_{l}&\mapsto\sum_{l=1}^{d}j_{O^{\circ}}\left(f_{l}\right)p^{r}\frac{d}{dZ_{l}}.
  \end{align*}
  We remark that this map has bounded image because
  $\j_{O^{\circ}}$ has bounded image (because
  $\epsilon$ is bounded by~\cite[Lemma~\ref*{lem:imagetildenu-Ainf}]{WiersigPeriods}),
  and because
  \begin{align*}
    \left\| j_{L_{r}}\left(p^{r}\partial_{l}\right)\left( Z^{\alpha}\right) \right\|
    =\left\| p^{r}\alpha_{l} Z^{\alpha-e_{l}}\right\|
    \leq p^{-r} \cdot 1 \cdot p^{-q|\alpha-1|}
    \leq p^{-r+q} p^{-q|\alpha|}
    \leq p^{-q|\alpha|}
    \left\|Z^{\alpha}\right\|
  \end{align*}
  for every $l=1,\dots, d$, $\alpha=\left(\alpha_{1},\dots,\alpha_{d}\right)\in\NN^{d}$
  with $\alpha_{l}>0$, and $e_{l}=\left(0,\dots,1,\dots,0\right)$ the unit vector.

  To exploit the universal property of $U\left( L_{r}\right)$, we have to check:
  \begin{condition}\label{condition:UP-envelopingalgebra}
  \hfill
    \begin{itemize}
      \item[(a)] $j_{O^{\circ}}$ is a homomorphism of
      $k^{\circ}$-algebras.
      \item[(b)] $j_{L_{r}}$ is an $O^{\circ}$-Lie algebra
      homomorphism.
      \item[(c)] For all $f\in O^{\circ}$ and $P\in L_{r}$,
      $j_{L_{r}}(fP)=j_{O^{\circ}}\left(f\right)j_{L_{r}}(P)$.
      \item[(d)] For all $f\in O^{\circ}$ and $P\in L_{r}$,
      $\left[j_{L_{r}}(P),j_{O^{\circ}}(f)\right]=j_{O^{\circ}}\left(\sigma(P)(f)\right)$.
    \end{itemize}
  \end{condition}
  
  By~\cite[\S 2.1]{AW19}, this
  would indeed give a map
  \begin{equation*}
    U\left( L_{r}\right)
    \to\intHom_{B^{q,+}}\left( O\kern -1,00pt B^{q,+}, O\kern -1,00pt B^{q,+}\right)
  \end{equation*}
  of $W(\kappa)$-algebras. From the discussion above,
  it follows that it has bounded image. Therefore, we can extend
  this map by continuity to the desired
  morphism~(\ref{eq:Ulr-toHomAqOAqOAq}).
  
  Conditions~\ref{condition:UP-envelopingalgebra} (a) and (c) are obvious. It remains to check
  (b) and (c).
    
  \begin{lem}\label{lem:DtBlam-Dtalgebraobject-derivations-vanish-on-coordinates}
    A bounded $W(\kappa)$-linear derivation $D\colon O^{\circ}\to O\kern -1,00pt B^{q,+}$
    is a bounded $W(\kappa)$-linear map which satisfies the Leibniz rule.
    If $D|_{\mathcal{O}\left(\TT^{d}\right)^{\circ}}=0$, then $D=0$.
  \end{lem}
    
  \begin{proof}
    Equip both $O\kern -1,00pt B^{q,+}$ and
    $O^{\circ}=\mathcal{O}(V)^{\circ}$ with the $p$-adic norms.
    $D$ is still continuous because it is $W(\kappa)$-linear.
    Recall Lemma~\ref{lem:localise-torsionfree-modules}. The following
    is a bounded derivation between two $k$-Banach spaces:
    \begin{equation*}
      D[1/p]\colon
      \mathcal{O}(V)^{\circ}
      \to
      O\kern -1,00pt B^{q,+}.
    \end{equation*}
    By~\cite[\S 3.6, page 64]{FvdP04}, it is identified
    with a bounded linear map
    \begin{equation}\label{eq:DtBlam-Dtalgebraobject-derivations-vanish-on-coordinates-1}
      \Omega_{\mathcal{O}\left(V\right)/k}\to O\kern -1,00pt B^{q,+}.
    \end{equation}
    \emph{Loc. cit.} identifies
    the composition of~(\ref{eq:DtBlam-Dtalgebraobject-derivations-vanish-on-coordinates-1})
    and
    $\Omega_{\mathcal{O}\left(\TT^{d}\right)/k}\cong\Omega_{\mathcal{O}\left(V\right)/k}$
    with
    \begin{equation*}
      D[1/p]|_{\mathcal{O}\left(\TT^{d}\right)}
      =D|_{\mathcal{O}\left(\TT^{d}\right)^{\circ}}[1/p]
      =0.
    \end{equation*}
    Therefore (\ref{eq:DtBlam-Dtalgebraobject-derivations-vanish-on-coordinates-1})
    is the zero map and so is $D[1/p]$.
    $D=0$ follows.
  \end{proof}
  
  \begin{lem}\label{lem:DtBlam-Dtalgebraobject-3}
    The following diagram commutes:
     \begin{equation*}
       \begin{tikzcd}[column sep = huge]
         \empty &
         \intHom_{B^{q,+}}\left(O\kern -1,00pt B^{q,+},O\kern -1,00pt B^{q,+}\right) \arrow{d}{-\circ\epsilon} \\
         L_{r} \arrow{ru}{j_{L_{r}}} \arrow[swap]{r}{P \mapsto \epsilon \circ \sigma(P)} &
         \intHom_{W(\kappa)}\left(O^{\circ},O\kern -1,00pt B^{q,+}\right).
       \end{tikzcd}
     \end{equation*}
   \end{lem}
  
   \begin{proof}
     Fix $P\in L_{r}$. Both $j_{L_{r}}\left(P\right)\circ \epsilon$ and
     $\epsilon\circ\sigma(P)$ are bounded $W(\kappa)$-linear derivations.
     Because of Lemma~\ref{lem:DtBlam-Dtalgebraobject-derivations-vanish-on-coordinates},
     it suffices to check
     \begin{equation*}
      \left( \epsilon \circ \sigma(P)\right) |_{\mathcal{O}\left(\TT_{k_{0}}\right)^{\circ}}
      = j_{L_{r}} \left( P \right) |_{\mathcal{O}\left(\TT_{k_{0}}\right)^{\circ}}.
     \end{equation*}
     Fix the $O^{\circ}$-basis $p^{r}\partial_{1},\dots,p^{r}\partial_{d}$
     for $L_{r}$ and write $P=\sum_{i=1}^{d}f_{i}p^{r}\partial_{i}$
     with $f_{i}\in O^{\circ}$. Both
     $\left(\epsilon \circ \sigma( - )\right) |_{\mathcal{O}\left(\TT_{k_{0}}\right)^{\circ}}$
     and $j_{L_{r}} \left( - \right) |_{\mathcal{O}\left(\TT_{k_{0}}\right)^{\circ}}$
     are $O^{\circ}$-linear, thus it suffices to check
     \begin{equation*}
     \left(\epsilon \circ \sigma(p^{r}\partial_{l})\right) |_{\mathcal{O}\left(\TT_{k_{0}}\right)^{\circ}}
      =j_{L} \left( p^{r}\partial_{l} \right) |_{\mathcal{O}\left(\TT_{k_{0}}\right)^{\circ}}
     \end{equation*}
     for every $l=1,\dots,d$. Fix such an $l$. Now computing
     \begin{equation}\label{eq:DtBlam-Dtalgebraobject-3-3}
     \left(\epsilon \circ \sigma(p^{r}\partial_{l})\right) \left(T_{j}\right)
     = j_{L_{t}} \left( p^{r}\partial_{l} \right) \left(j_{O^{\circ}}\left(T_{j}\right)\right)
     \end{equation}
     for every $j=1,\dots,d$ suffices. 
     This is a direct computation:
     \begin{equation*}
       \left(\epsilon \circ \sigma(p^{r}\partial_{l})\right) \left(T_{j}\right)
       =\epsilon\left( p^{r}\frac{d}{dT_{l}}\left(T_{j}\right) \right)
       =p^{r}\delta_{lj},
     \end{equation*}
     where $\delta_{ij}$ denotes the Kronecker delta, and
     \begin{equation*}
       j_{L_{r}} \left( p^{r}\partial_{l} \right) \left(\epsilon\left(T_{j}\right)\right)
       =p^{r}\frac{d}{dZ_{l}}\left( \left[T_{j}^{\flat}\right] + Z_{j} \right)
       =p^{r}\delta_{lj}.
     \end{equation*}     
     The last step uses that $\left[T_{j}^{\flat}\right]\in B^{q,+}$ is a constant and
     $d/dZ_{l}$ is an $B^{q,+}$-linear derivation. This gives the
     identity~\ref{eq:DtBlam-Dtalgebraobject-3-3}, and thus finishes the
     proof of Lemma~\ref{lem:DtBlam-Dtalgebraobject-3}.
     \end{proof}
     
     We verify Condition~\ref{condition:UP-envelopingalgebra}(b) in the following.
     
     \begin{lem}\label{lem:DtBlam-Dtalgebraobject-4}
       The morphism $j_{L_{r}}$ is an $O^{\circ}$-Lie algebra
       homomorphism.
     \end{lem}
     
     \begin{proof}
       $j_{L_{r}}$ is $k^{\circ}$-linear. It remains to show that it preserves the
       Lie bracket:
       \begin{equation}\label{eq:DtBlam-Dtalgebraobject-4-1}
         \left[ j_{L_{r}}(P) , j_{L_{t}}(Q) \right]
         =j_{L_{r}}\left( [P,Q] \right),
       \end{equation}
       for any $P,Q\in L_{r}$. For all $j=1,\dots,d$,
       it suffices to check
       \begin{equation*}
         \left[ j_{L_{r}}(P) , j_{L_{r}}(Q) \right]\left( Z_{j} \right)
         =j_{L_{r}}\left( [P,Q] \right)\left( Z_{j} \right)
       \end{equation*}
       because both sides are derivations.
       But then we would have to compute
       \begin{equation}\label{eq:DtBlam-Dtalgebraobject-4-2}
         \left[ j_{L_{r}}(P) , j_{L_{r}}(Q) \right]\circ\epsilon
         =j_{L_{r}}\left( [P,Q] \right)\circ\epsilon.
       \end{equation}
       Indeed, if~(\ref{eq:DtBlam-Dtalgebraobject-4-2}) is true, we could compute
       \begin{align*}
         \left[ j_{L_{r}}(P) , j_{L_{r}}(Q) \right]\left( Z_{j} \right)
         &=\left[ j_{L_{r}}(P) , j_{L_{r}}(Q) \right]\left( \left[T_{j}^{\flat}\right] + Z_{j} \right) \\
         &=\left( \left[ j_{L_{r}}(P) , j_{L_{r}}(Q) \right]\circ\epsilon\right)\left( T_{j} \right) \\
         &\stackrel{\text{(\ref{eq:DtBlam-Dtalgebraobject-4-2})}}{=}
         \left(j_{L_{r}}\left( [P,Q] \right)\circ\epsilon\right)\left( T_{j} \right) \\
         &= j_{L_{r}}\left( [P,Q] \right)\left( \left[T_{j}^{\flat}\right] + Z_{j} \right) \\
         &= j_{L_{r}}\left( [P,Q] \right)\left( Z_{j} \right).
       \end{align*}
       In particular,~(\ref{eq:DtBlam-Dtalgebraobject-4-2}) implies~(\ref{eq:DtBlam-Dtalgebraobject-4-1}).
       The following computation
       \begin{align*}
         &\left[ j_{L_{r}}(P) , j_{L_{r}}(Q) \right]\circ\epsilon \\
         &=\left(j_{L_{r}}(P) \circ j_{L_{r}}(Q) - j_{L_{r}}(Q) \circ j_{L_{r}}(P)\right) \circ \epsilon \\
         &=j_{L_{r}}(P) \circ \left(j_{L_{r}}(Q) \circ \epsilon\right)
            - j_{L_{r}}(Q) \circ \left(j_{L_{r}}(P) \circ \epsilon\right) \\
         &\stackrel{\text{\ref{lem:DtBlam-Dtalgebraobject-3}}}{=}
          \left(j_{L_{r}}(P) \circ \epsilon\right) \circ \sigma(Q)
            - \left(j_{L_{r}}(Q) \circ \epsilon\right) \circ \sigma(P)\\
         &\stackrel{\text{\ref{lem:DtBlam-Dtalgebraobject-3}}}{=}
          \left(\epsilon \circ \sigma(P) \right) \circ \sigma(Q)
            - \left(\epsilon \circ \sigma(Q) \right) \circ \sigma(P)\\
         &\stackrel{\text{\ref{lem:DtBlam-Dtalgebraobject-3}}}{=}
          \epsilon \circ \left(\sigma(P) \circ \sigma(Q)
            - \sigma(Q) \circ \sigma(P)\right)\\
         &=\epsilon\circ\left( \sigma(P) \circ\sigma(Q)
         -\sigma(Q)\circ\sigma(P)\right) \\
         &=\epsilon \circ \left[\sigma(P),\sigma(Q)\right] \\
         &=\epsilon \circ \sigma\left(\left[ P, Q \right]\right) \\
         &\stackrel{\text{\ref{lem:DtBlam-Dtalgebraobject-3}}}{=}
          j_{L_{r}}\left(\left[ P, Q \right] \right)\circ\epsilon
       \end{align*}
       checks~(\ref{eq:DtBlam-Dtalgebraobject-4-2}).
     \end{proof}
     
     We verify Condition~\ref{condition:UP-envelopingalgebra}(d) in the following.
     
      \begin{lem}\label{lem:DtBlam-Dtalgebraobject-5}
        For all $f\in O$ and $P\in L_{r}$,
       $\left[j_{L_{r}}(P),j_{O^{\circ}}(f)\right]=j_{O^{\circ}}\left(\sigma(P)(f)\right)$.
     \end{lem}    
     
     \begin{proof}
     Write $P=\sum_{l=1}^{d}f_{l}p^{r}\partial_{l}$.
     Compute for every $h\in O\kern -1,00pt B^{q,+}$
     and $l=1,\dots,d$,
     \begin{align*}
       \left[j_{O^{\circ}}\left(f_{l}\right)p^{r}\frac{d}{dZ_{l}},j_{O^{\circ}}(f)\right](h)
       &=j_{O^{\circ}}\left(f_{l}\right)\left[p^{r}\frac{d}{dZ_{l}},j_{O^{\circ}}\left(f\right)\right](h) \\
       &=\epsilon\left(f_{l}\right)\left(p^{r}\frac{d}{dZ_{l}}\left(\epsilon\left(f\right)\right)\right)h \\
       &=\left(j_{O^{\circ}}\left(f_{l}\right)p^{r}\frac{d}{dZ_{l}}\right)\left(\epsilon\left(f\right)\right)h.
      \end{align*}
      This establishes the identity
      at the left-hand side here:
      \begin{equation*}
        \left[j_{L_{r}}(P),j_{O^{\circ}}(f)\right](h)
        =\left(j_{L_{r}}(P)\circ\epsilon\right)(f)h
        \stackrel{\text{\ref{lem:DtBlam-Dtalgebraobject-3}}}{=}
        \left(\epsilon\circ\sigma(P)\right)(f)h
        =j_{O^{\circ}}\left(\sigma(P)(f)\right)h.
      \end{equation*}
      The identity at the right-hand side comes from the
      definition of $j_{O^{\circ}}$.
   \end{proof} 

   We have thus verified Condition~\ref{condition:UP-envelopingalgebra}.
   This gives the morphism~(\ref{eq:Ulr-toHomAqOAqOAq})
   and thus the
   $\widehat{U\left( \mathcal{B} \otimes_{\mathcal{A}}\mathcal{L}_{r}\right)}$-$\BB_{\dR}^{q,+}(U)$-bimodule
   structure on $\BB_{\dR}^{q,+}\left( U \right)\left\< \frac{Z_{1},\dots,Z_{d}}{p^{q}}\right\>$.
   Invert $p$ via $-\widehat{\otimes}_{W(\kappa)}k_{0}$.
   Lemma~\ref{lem:bimodulestructure-limcolim} applies because
   Condition~\ref{condition:gluethemaps-along-limcolim} is satisfied, giving
   the desired $\Dcap(V)$-$\BB_{\dR}^{q,+}(U)$-bimodule
   structure on $\OB_{\dR}^{q,+}\left(U\right)$. It is
   functorial in both $V$ and $U$, giving $\OB_{\dR}^{q,+}$ the structure
   of a $\nu^{-1}\Dcap$-$\BB_{\dR}^{q,+}$-bimodule object.
      
    We discussed above that uniqueness is immediate
    once we have checked that the canonical
    morphism~(\ref{eq:Otensorbla-to-OBla-epi}) is a morphism
    of sheaves of $\nu^{-1}\Dcap$-$\BB_{\dR}^{q,+}$-bimodule objects.
    Given $V$ and $U$ as above, we may check that
    \begin{equation}\label{eq:OVotimesBlaplus-to-OBlaplusU-onebeforethat}
      \mathcal{O}(V) \widehat{\otimes}_{k_{0}}\BB_{\dR}^{q,+}\left( U\right)
      \to\OB_{\dR}^{q,+}\left(U\right)
    \end{equation}
    is a morphism of $\Dcap(V)$-$\BB_{\dR}^{q,+}(U)$-bimodule objects.
    It is $\BB_{\dR}^{q}(U)$-linear and $\mathcal{O}(V)$-linear.
    Therefore, it remains to compare the actions of the differentials
    $\partial_{1},\dots,\partial_{d}$.
    But we have étale coordinates on $\mathcal{O}(V)$, and it suffices to compare
    the actions on these. First, compute for all $T_{j}\in\mathcal{O}(V)$,
    \begin{equation*}
      \partial_{l}\cdot T_{j} = \delta_{lj} .
    \end{equation*}
    Second, the map~(\ref{eq:OVotimesBlaplus-to-OBlaplusU-onebeforethat})
    sends $T_{j}$ to $\left[ T^{\flat}_{j}\right] + Z_{j}$. This implies
    \begin{equation*}
      \partial_{l} \cdot \left(\left[ T^{\flat}_{j}\right] + Z_{j}\right) = \delta_{lj},
    \end{equation*}
    see the proof of Lemma~\ref{lem:DtBlam-Dtalgebraobject-3}.
    Therefore, the actions of the differentials
    $\partial_{1},\dots,\partial_{d}$ coincide.
\end{proof}

Furthermore, the proof of
Theorem~\ref{thm:OBlaplus-bimodule} implies the following
Lemma~\ref{lem:OBlaqplus-bimodule-veryq}
and~\ref{lem:OBlaqplus-bimodule-varyq-unitmaps}.

\begin{lem}\label{lem:OBlaqplus-bimodule-veryq}
  Let $X$ be affinoid. Fix an étale morphism $g\colon X\to\TT^{d}$
  and $q\geq q_{g}$ (see~\cite[Notation~\ref*{notation:qgbound}]{WiersigPeriods}).
  Consider the diagram
  \begin{equation}\label{cd:lem:OBlaqplus-bimodule-veryq}
  \begin{tikzcd}
    \left(\nu^{-1}\Dcap \widehat{\otimes}_{k_{0}} \BB_{\dR}^{q+1,+}\right)
      \widehat{\otimes}_{k_{0}} \OB_{\dR}^{q+1,+}
      \arrow{r} &
    \OB_{\dR}^{q+1,+} \\
    \left(\nu^{-1}\Dcap \widehat{\otimes}_{k_{0}} \BB_{\dR}^{q,+}\right)
      \widehat{\otimes}_{k_{0}} \OB_{\dR}^{q,+}
      \arrow{r}\arrow{u} &
    \OB_{\dR}^{q,+}, \arrow{u}
  \end{tikzcd}
  \end{equation}
  where the horizontal maps are the bimodule action maps
  coming from Theorem~\ref{thm:OBlaplus-bimodule},
  and the vertical maps are the canonical ones.
  This diagram commutes~\ref{lem:OBlaqplus-bimodule-veryq}.
\end{lem}

\begin{lem}\label{lem:OBlaqplus-bimodule-varyq-unitmaps}
  Let $X$ be affinoid. Fix an étale morphism $g\colon X\to\TT^{d}$
  and $q\geq q_{g}$ (see~\cite[Notation~\ref*{notation:qgbound}]{WiersigPeriods}).
  Consider the diagram  \begin{equation}\label{cf:OBlaqplus-bimodule-varyq-unitmaps}
  \begin{tikzcd}
    k_{0}
      \arrow{r} &
    \OB_{\dR}^{q+1,+} \\
    k_{0}
      \arrow{r}\arrow{u} &
    \OB_{\dR}^{q,+}, \arrow{u}
  \end{tikzcd}
  \end{equation}
  where the horizontal maps are the unit maps
  of the bimodule structures coming from Theorem~\ref{thm:OBlaplus-bimodule},
  and the vertical maps are the canonical ones.
  This diagram~(\ref{cf:OBlaqplus-bimodule-varyq-unitmaps}) commutes.
\end{lem}

\begin{cor}\label{cor:OBdRdagplus-bimodule}
  There exists a unique
  $\nu^{-1}\Dcap$-$\BB_{\dR}^{\dag,+}$-bimodule structure on $\OB_{\dR}^{\dag,+}$
  such that the canonical morphism
  \begin{equation}\label{eq:Otensorbla-to-OBdag-epi}
    \nu^{-1}\mathcal{O}\widehat{\otimes}_{k_{0}}\BB_{\dR}^{\dag,+}
    \to\OB_{\dR}^{\dag,+}
  \end{equation}  
  is a morphism of $\nu^{-1}\Dcap$-$\BB_{\dR}^{\dag,+}$-bimodule objects.
\end{cor}

\begin{proof}
  Locally on $X$, we construct the bimodule structure from the ones
  in Theorem~\ref{thm:OBlaplus-bimodule}.
  One does this via Lemma~\ref{lem:bimodulestructure-colim},
  which applies thanks to Lemma~\ref{lem:OBlaqplus-bimodule-veryq}
  and~\ref{lem:OBlaqplus-bimodule-varyq-unitmaps}.
  These glue bimodule structures glue, because they are unique.
  
  It remains to check that the resulting $\nu^{-1}\Dcap$-$\BB_{\dR}^{\dag,+}$-bimodule
  structure on $\OB_{\dR}^{\dag,+}$ is unique.
  By Corollary~\ref{cor:epi-uniquebimoduleaction}, it suffices
  to show that~(\ref{eq:Otensorbla-to-OBdag-epi}) is an epimorphism.
  This follows from Lemma~\ref{lem:Otensorbla-to-OBla-epi},
  as~(\ref{eq:Otensorbla-to-OBdag-epi}) is the colimit of the map
  (\ref{eq:Otensorbla-to-OBla-epi}) along $q\to\infty$.  
\end{proof}

Repeating the previous arguments, one easily deduces:

\begin{cor}\label{thm:OBpdRq-bimodule}
  Let $X$ be affinoid. Fix an étale morphism $g\colon X\to\TT^{d}$
  and $q\geq q_{g}$ (see~\cite[Notation~\ref*{notation:qgbound}]{WiersigPeriods}).
  Then there exists a unique
  $\nu^{-1}\Dcap$-$\BB_{\pdR}^{q}$-bimodule structure on $\OB_{\pdR}^{q}$
  such that the canonical morphism
  \begin{equation}\label{eq:Otensorbla-to-OBla-epi}
    \nu^{-1}\mathcal{O}\widehat{\otimes}_{k_{0}}\BB_{\dR}^{q}
    \to\OB_{\pdR}^{q}
  \end{equation}  
  is a morphism of $\nu^{-1}\Dcap$-$\BB_{\pdR}^{q}$-bimodule objects.
\end{cor}

\begin{cor}\label{cor:OBpdRdag-bimodule}
  There exists a unique
  $\nu^{-1}\Dcap$-$\BB_{\pdR}^{\dag}$-bimodule structure on $\OB_{\pdR}^{\dag}$
  such that the canonical morphism
  \begin{equation*}
    \nu^{-1}\mathcal{O}\widehat{\otimes}_{k_{0}}\BB_{\pdR}^{\dag}
    \to\OB_{\pdR}^{\dag}
  \end{equation*}  
  is a morphism of $\nu^{-1}\Dcap$-$\BB_{\pdR}^{\dag}$-bimodule objects.
\end{cor}


\subsection{Definitions of the solution and de Rham functors}
\label{subsec:solutionderHamfunctorsDcapmod-solutionpaper}

In an ideal world, we could define the solution functor to be
\begin{equation}\label{eq:idealworld-sol}
  \D\left( \Dcap \right)^{\op} \to \D\left( \BB_{\pdR}^{\dag}\right),
  \mathcal{M}^{\bullet} \mapsto
    \varinjlim_{q\gg0}\R\shHom_{\nu^{-1}\Dcap}\left( \nu^{-1}\cal{M}^{\bullet} , \OB_{\pdR}^{q} \right).
\end{equation}

\begin{prob}\label{prob:idealworld-sol}
  With~(\ref{eq:idealworld-sol}), we run into the following two difficulties:
  \begin{itemize}
    \item[(i)] We do not know how to derive $\shHom_{\nu^{-1}\Dcap}$. The issue is that we do not know
      whether $\IndBan_{k}$ has enough injectives. In fact, we do not expect this to be the case.
    \item[(ii)] Given a sheaf $\cal{M}$ of $\Dcap$-modules,
      $\varinjlim_{q\gg0}\shHom_{\nu^{-1}\Dcap}\left( \nu^{-1}\cal{M}^{\bullet} , \OB_{\pdR}^{q} \right)$
      becomes a sheaf of $\BB_{\pdR}^{\dag}$-ind-Banach modules.
      However, given
      a complex $\cal{M}^{\bullet}$ of sheaves of $\Dcap$-ind-Banach modules,
      and even if we could derive the $\shHom$-functor, it is not clear to us how to make
      $\varinjlim_{q\gg0}\R\shHom_{\nu^{-1}\Dcap}\left( \nu^{-1}\cal{M}^{\bullet} , \OB_{\pdR}^{q} \right)$
      an object of $\D\left(\BB_{\pdR}^{\dag}\right)$.
  \end{itemize}
\end{prob}

In order to address Problem~\ref{prob:idealworld-sol}(ii), we aim to replace
the colimit in~(\ref{eq:idealworld-sol}) by a pullback functor. Therefore,
we introduce the following notation.

\begin{defn}
  $X_{\proet}\times\NN_{\gg0}^{\op}$ denotes the following site.
  \begin{itemize}
    \item Its objects are pairs $(U,q)$, where $U\in X_{\proet}$ and $q\in\NN_{\geq2}$
      such that there exists a $\nu^{-1}\Dcap$-$\BB_{\pdR}^{q}$-bimodule structure on $\OB_{\pdR}^{q}$
      such that the canonical morphism
      \begin{equation*}
        \nu^{-1}\mathcal{O}(U)\widehat{\otimes}_{k_{0}}\BB_{\dR}^{q}(U)
        \to\OB_{\pdR}^{q}(U)
      \end{equation*}
      is a morphism of $\nu^{-1}\Dcap(U)$-$\BB_{\pdR}^{q}(U)$-bimodules.
      (If such a bimodule structure exists, then it is necessarily unique.
      This is explained in \S\ref{subsec:OBla-bimodule}.)
    \item Its morphisms are
      \begin{equation*}
        \Hom\left(\left(U,q^{\prime}\right),\left(V,q\right)\right)
        =\begin{cases}
          \Hom\left(U,V\right) &\text{ if $q\leq q^{\prime}$, and} \\
          \emptyset &\text{ otherwise}.
        \end{cases}
      \end{equation*}
    \item We equip this category with the Grothendieck topology whose coverings are
      \begin{equation*}
        \left\{ \left( U_{i},q\right) \to \left(U,q\right)\right\}_{i\in I}
      \end{equation*}
      for all coverings $\left\{U_{i}\to U\right\}_{i\in I}$ in $X_{\proet}$ and $q\gg0$.
  \end{itemize}
\end{defn}

Let $X_{w}$ denote the category
whose objects are the affinoid subdomains of $X$ and whose morphisms
are the inclusions, carrying the weak Grothendieck topology.
Then $\lambda\colon X_{\proet}\times\NN_{\gg0}^{\op}\to X_{w}$ is the morphism
of sites $\left(\nu^{-1}\left(U\right),q\right)\mapsfrom U$, where $q$ is minimal with respect
to the property that $\left(\nu^{-1}\left(U\right),q\right)\in X_{\proet}\times\NN_{\gg0}^{\op}$.
This is well-defined. Indeed, any affinoid $U\in X_{w}$ can be covered by finitely many affinoids $V$
admitting étale coordinates $U_{i}\to\TT^{d_{i}}$ for which we then get the
constant $q_{g_{i}}$ (see~\cite[Notation~\ref*{notation:qgbound}]{WiersigPeriods}). Now define $q:=\min_{i}q_{g_{i}}$.

On $X_{\proet}\times\NN_{\gg0}^{\op}$, we introduce the following two sheaves
\begin{equation*}
\begin{split}
  \OB_{\pdR}^{*}\colon \left(U,q\right)&\mapsto\OB_{\pdR}^{q}(U), \text{ and} \\
  \BB_{\pdR}^{*}\colon \left(U,q\right)&\mapsto\BB_{\pdR}^{q}(U)
\end{split}
\end{equation*}
of $k_{0}$-ind-Banach algebras. There is canonical morphism
\begin{equation}\label{eq:OtensorbdRstarplus-to-OBdrstarplus-epi}
  \lambda^{-1}\mathcal{O}\widehat{\otimes}_{k_{0}}\BB_{\pdR}^{*}
  \to\OB_{\pdR}^{*}
\end{equation}  
of $k_{0}$-ind-Banach algebras on $X_{\proet}\times\NN_{\gg0}^{\op}$.

\begin{lem}\label{lem:OBdRstarplus-bimodule}
  There exists a
  $\lambda^{-1}\Dcap$-$\BB_{\pdR}^{*}$-bimodule structure on $\OB_{\pdR}^{*}$
  such that the canonical morphism~(\ref{eq:OtensorbdRstarplus-to-OBdrstarplus-epi})
  is a morphism of $\lambda^{-1}\Dcap$-$\BB_{\pdR}^{*}$-bimodule objects.
  It is unique.
\end{lem}

\begin{proof}
  This follows easily from the discussion in \S\ref{subsec:OBla-bimodule}.
\end{proof}

We have the projection $\sigma\colon X_{\proet}\to X_{\proet}\times\NN_{\gg0}^{\op}$,
$U\mapsfrom(U,q)$. There is a canonical morphism $\sigma^{-1}\BB_{\pdR}^{*}\to\BB_{\pdR}^{\dag}$
of sheaves of $k_{0}$-ind-Banach algebras making $\sigma$ a morphism of ringed sites
\begin{equation*}
  \left(X_{\proet},\BB_{\pdR}^{\dag}\right)
  \to \left(X_{\proet}\times\NN_{\gg0}^{\op},\BB_{\pdR}^{*}\right).
\end{equation*}
The pullback $\sigma^{*}$ is supposed to replace the colimit in~(\ref{eq:idealworld-sol}),
cf. Definition~\ref{defn:sol} below.

Next, we address Problem~\ref{prob:idealworld-sol}(i) by working with sheaves valued in the left heart
$\IndBan_{\I\left(k_{0}\right)}=\LH\left(\IndBan_{k_{0}}\right)$. This is a Grothendieck abelian category
by~\cite[Lemma 4.20 and a remark following Definition 3.2]{Bo21},
therefore it has enough injectives. Given an arbitrary site $Y$,
$\Sh\left(Y,\IndBan_{\I\left(k_{0}\right)}\right)$ admits a closed symmetric monoidal structure
by~\cite[\S 3.4, which applies by Corollary 4.24]{Bo21}.
By~\cite[Corollary 1.2.28]{Sch99} and~\cite[Lemma 3.33]{Bo21}, the functor
\begin{equation*}
  \I\colon\Sh\left(Y,\IndBan_{k_{0}}\right) \to \Sh\left(Y,\IndBan_{\I\left(k_{0}\right)}\right),
\end{equation*}
is fully faithful and strongly monoidal. As a consequence, Lemma~\ref{lem:OBdRstarplus-bimodule}
gives a $\lambda^{-1}\I\left(\Dcap\right)$-$\I\left(\BB_{\pdR}^{*}\right)$-bimodule structure on $\I\left(\OB_{\pdR}^{*}\right)$
and $\sigma$ gives rise to a morphism of ringed sites
\begin{equation*}
  \left(X_{\proet},\I\left(\BB_{\pdR}^{\dag}\right)\right)
  \to \left(X_{\proet}\times\NN_{\gg0}^{\op},\I\left(\BB_{\pdR}^{*}\right)\right).
\end{equation*}
Here, we have also used that pullback $\lambda^{-1}$ commutes with $\I$,
cf.~\cite[Lemma 3.39(ii)]{Bo21}.

In the following definition, we can derive
$\sigma^{*}$ because it is a left adjoint by the discussion in~\cite[\S 3.7]{Bo21}.
One derives the internal homomorphisms via \emph{loc. cit.} Proposition 3.35.

\begin{defn}\label{defn:sol}
  The \emph{solution functor} is
  \begin{equation*}
    \SolB_{\pdR}^{\dag}\colon
      \D\left(\I\left(\Dcap\right)\right)
      \to
      \D\left(\I\left(\BB_{\pdR}^{\dag}\right)\right),
      \cal{M}^{\bullet}\mapsto
        \rL\sigma^{*}\R\shHom_{\lambda^{-1}\I\left(\Dcap\right)}
        \left(
          \lambda^{-1}\cal{M}^{\bullet},
          \I\left(\OB_{\pdR}^{*}\right)
        \right).
  \end{equation*}
\end{defn}

The duality functor
$\DD\colon\D\left( \Dcap \right)
  \to\D\left( \Dcap \right)^{\op}$
has been introduced in~\cite[\S 7.2]{Bo21}.
Since $\D\left( \Dcap \right)\simeq\D\left( \I\left(\Dcap\right) \right)$,
we can make the following Definition~\ref{defn:de-Rham-Dcapmodules},
where $d$ is the dimension of $X$.
    
\begin{defn}\label{defn:de-Rham-Dcapmodules}
Following the classical~\cite[Proposition 4.2.1]{HTTDmodules},
the \emph{de Rham functor} is
\begin{equation*}
  \dRfunctorB_{\pdR}^{\dag}\colon
  \D\left( \I\left(\Dcap\right) \right) \to \D\left( \I\left(\BB_{\pdR}^{\dag}\right)\right),
  \mathcal{M}^{\bullet} \mapsto \Sol\left(\DD\left(\mathcal{M}^{\bullet}\right)\right)\left[d\right].
\end{equation*}
\end{defn}

The previous discussion also goes through with coefficients in $\BB_{\dR}^{\dag,+}$:

\begin{defn}
 Define the \emph{solution} and \emph{de Rham functors}
  \begin{align*}
    \Sol:=\SolB_{\dR}^{\dag,+}&\colon
      \D\left(\I\left(\Dcap\right)\right)
      \to
      \D\left(\I\left(\BB_{\dR}^{\dag,+}\right)\right),
      \cal{M}^{\bullet}\mapsto
        \rL\sigma^{*}\R\shHom_{\lambda^{-1}\I\left(\Dcap\right)}
        \left(
          \lambda^{-1}\cal{M}^{\bullet},
          \I\left(\OB_{\dR}^{*,+}\right)
        \right) \\
    \dRfunctor:=\dRfunctorB_{\dR}^{\dag,+}&\colon
    \D\left( \I\left(\Dcap\right) \right) \to \D\left( \I\left(\BB_{\pdR}^{\dag}\right)\right),
    \mathcal{M}^{\bullet} \mapsto \Sol\left(\DD\left(\mathcal{M}^{\bullet}\right)\right)\left[d\right].
  \end{align*}
  where $\sigma^{*}=\BB_{\dR}^{\dag,+}\widehat{\otimes}_{\sigma^{-1}\BB_{\dR}^{*,+}}-$.
\end{defn}


\section{Reconstruction theorems and fully faithfulness}
\label{ch:recthms-recpaper}

Throughout \S\ref{ch:recthms-recpaper},
we fix a smooth rigid-analytic $k$-variety $X$ of dimension $d$.


\subsection{$\cal{C}$-complexes}
\label{subsec:C-recpaper}

We fix the notation introduced in \S\ref{subsec:Dcap-modules}.
In particular, given an étale morphism $g\colon X\to\TT^{d}$ to the $d$-dimensional
torus,
we obtain the certain constant $r_{g}\in\NN$ such that for all $r\geq r_{g}$,
the sites $X_{g,r}$ of $\cal{L}_{r}$-admissible affinoid subdomains
of $X$ are defined.

See~\cite[\S 8 and \S 9]{AW19} for the definition of the category of
coadmissible $\Dcap$-modules as a full subcategory
of abstract $\Dcap$-modules. \cite{Bo21} realised it as a full
subcategory of the category of sheaves of $\Dcap$-ind-Banach modules.
We are interested in a derived analog: the category
\begin{equation*}
  \D_{\cal{C}}\left(\Dcap\right)\subseteq\D\left(\Dcap\right)
\end{equation*}
of $\cal{C}$-complexes,
which is a full subcategory of the derived category of $\Dcap$-ind-Banach modules.
Here, we do not recall the precise definition as in~\cite[\S 8.1]{Bo21}.
Instead, we refer to the characterisation as in Lemma~\ref{lem:C-complexdefn-reconstructionpaper}.
For its formulation, we work in the derived category of
$\I\left(\Dcap\right)$-module objects instead, where $\I$ denotes the morphism into the left heart:
\begin{equation*}
  \I\colon\Sh\left(X,\IndBan_{k}\right) \to 
  \LH\left(\Sh\left(X,\IndBan_{k}\right)\right)
  \simeq\Sh\left(X,\IndBan_{\I\left(k\right)}\right).
\end{equation*}
We refer the reader to~\cite[\S 3, which applies because of Corollary 4.24]{Bo21}
for a detailed discussion of the relevant category theoretical background. Note that
passing to the left heart is allowed, as $\I$ induces the equivalence of
categories $\D\left(\Dcap\right)\simeq\D\left(\I\left(\Dcap\right)\right)$,
cf.~\cite[Proposition 3.8]{Bo21}.

\begin{defn}\label{defn:perfectcomplex-recpaper}
  Consider a closed symmetric monoidal quasi-abelian category $\E$, containing a ring object $R$.
  Then the category of $R$-modules is again quasi-abelian, cf.~\cite[Proposition 1.5.1]{Sch99},
  thus it admits a derived category $\D\left(R\right)$. The \emph{category
  $\D_{\perf}^{\bd}\left(R\right)$ of bounded perfect complexes
  of $R$-modules} is the smallest full triangulated subcategory of
  $\D\left(R\right)$ which contains $R$, and is closed under
  direct summands and isomorphisms. A \emph{bounded perfect
  complex of $R$-modules} is an object of $\D_{\perf}^{\bd}\left(R\right)$.
\end{defn}

\begin{lem}\label{lem:C-complexdefn-reconstructionpaper}
  Suppose $X$ is affinoid and equipped with an étale morphism $g\colon X\to\TT^{d}$.
  An object $\mathcal{M}^{\bullet}$ in the derived category of
  of $\I\left(\Dcap\right)$-module objects on $X$ is a $\mathcal{C}$-complex if
  and only if
  \begin{itemize}
    \item[(i)] $\mathcal{M}_{r}^{\bullet}
    :=\I\left(\mathcal{D}_{r}\right)\widehat{\otimes}_{\I\left(\Dcap\right)|_{X_{r}}}^{\rL}\mathcal{M}^{\bullet}|_{X_{r}}$
    is a bounded perfect complex of sheaves of
    $\I\left(\mathcal{D}_{r}\right)$-module objects for every $r\geq r_{g_{i}}$, and
    \item[(ii)]
    $\Ho^{i}\left( \mathcal{M}^{\bullet}\right)\stackrel{\cong}{\longrightarrow}
    \varprojlim_{r\geq r_{g}}\Ho^{i}\left(\mathcal{M}_{r}^{\bullet}\right)$
    for every $i\in\ZZ$.
  \end{itemize}
\end{lem}

\begin{proof}
  See~\cite[\S 8.1]{Bo21} for the definition of $\cal{C}$-complexes.
  \emph{Loc. cit.} works with the left heart of the category of bornological $\Dcap$-modules.
  This is equivalent to our setup by~\cite[Proposition 4.22(v)]{Bo21}.
  The only difference is that we require $\mathcal{M}_{r}^{\bullet}$
  to be a bounded perfect complex, where \emph{loc. cit.} requires it
  to be a bounded complex with coherent cohomology. These notions
  agree by~\cite[Theorem 1.1]{Bode2024Auslanderregularitydiffoper}.
\end{proof}


\subsection{Theorems}\label{subsec:thms-recpaper}

In \S\ref{subsec:solutionderHamfunctorsDcapmod-solutionpaper}, we defined the solution functor
\begin{equation*}
  \Sol\colon
    \D\left(\I\left(\Dcap\right)\right)^{\op}
    \to
    \D\left(\I\left(\BB_{\dR}^{\dag,+}\right)\right),
    \cal{M}^{\bullet}\mapsto
      \rL\sigma^{*}\R\shHom_{\lambda^{-1}\I\left(\Dcap\right)}
      \left(
        \lambda^{-1}\cal{M}^{\bullet},
        \I\left(\OB_{\dR}^{*,+}\right)
      \right).
\end{equation*}

By~\cite[Lemma 3.33]{Bo21},
Corollary~\ref{cor:OBpdRdag-bimodule} makes $\I\left(\OB_{\pdR}^{\dag}\right)$
a $\I\left(\nu^{-1}\Dcap\right)$-$\I\left(\BB_{\pdR}^{\dag}\right)$-bimodule object.
In the following Definition~\ref{defn:rec-reconstructionpaper}, we can derive
$\nu_{*}$ to a functor between the unbounded derived categories as explained
in~\cite[\S 3.7]{Bo21}.
One derives the internal homomorphisms via \emph{loc. cit.} Proposition 3.35.
We also use $\I\left(\nu^{-1}\Dcap\right)\cong\nu^{-1}\I\left(\Dcap\right)$,
cf. \emph{loc. cit.} Lemma 3.39.

\begin{defn}\label{defn:rec-reconstructionpaper}
  The \emph{reconstruction functor} is
  \begin{equation*}
    \Rec\colon\D\left( \I\left(\BB_{\dR}^{\dag,+}\right) \right) \to \D\left( \I\left(\Dcap\right) \right)^{\op},
    \mathcal{F}^{\bullet}
    \mapsto \R\nu_{*}\R\shHom_{\I\left(\BB_{\dR}^{\dag,+}\right)}\left( \mathcal{F}^{\bullet} , \I\left(\OB_{\pdR}^{\dag}\right) \right).
  \end{equation*}
\end{defn}

Now we define the morphisms which relate $\Sol$ and $\Rec$.

\begin{construction}\label{construct:rhocalMbullet-reconstructionpaper}
  In this Construction~\ref{construct:rhocalMbullet-reconstructionpaper},
  we always work in left hearts. Therefore, we omit $\I$ throughout from the notation.
  Furthermore, we consider the following morphisms of ringed sites.
  \begin{itemize}
    \item[(i)] As before, $\sigma$ and the canonical map
      $\sigma^{-1}\BB_{\dR}^{*,+}\to\BB_{\dR}^{\dag,+}$ induces a morphism of ringed sites
      \begin{equation*}
        \sigma\colon\left(X_{\proet},\BB_{\dR}^{\dag,+}\right)
        \to \left(X_{\proet}\times\NN_{\gg0}^{\op},\BB_{\dR}^{*,+}\right).
      \end{equation*}
    \item[(ii)] $\nu$ together with the canonical map
      $\nu^{-1}\Dcap
        \cong\nu^{-1}\Dcap\widehat{\otimes}_{k_{0}}k_{0}
        \to\nu^{-1}\Dcap\widehat{\otimes}_{k_{0}}\BB_{\dR}^{\dag,+}$
      give
      \begin{equation*}
        \nu\colon\left(X_{\proet},\nu^{-1}\Dcap\widehat{\otimes}_{k_{0}}\BB_{\dR}^{\dag,+}\right)
        \to \left(X,\Dcap\right).
      \end{equation*}
    \item[(iii)] Finally, $\lambda$ induces a morphism of ringed sites
      \begin{equation*}
        \lambda\colon\left(X_{\proet}\times\NN_{\gg0}^{\op},\lambda^{-1}\Dcap\widehat{\otimes}_{k_{0}}\BB_{\dR}^{*,+}\right)
        \to \left(X,\Dcap\right)
      \end{equation*}
      via the canonical map
        $\lambda^{-1}\Dcap
          \cong\lambda^{-1}\Dcap\widehat{\otimes}_{k_{0}}k_{0}
          \to\lambda^{-1}\Dcap\widehat{\otimes}_{k_{0}}\BB_{\dR}^{*,+}$.
  \end{itemize}

  Fix a complex $\cal{M}^{\bullet}$ of $\Dcap$-module objects.
  We have the canonical isomorphism
  \begin{equation}\label{eq:construct:rhocalMbullet-someiso-reconstructionpaper}
  \begin{split}
    \rL\sigma^{*}\rL\lambda^{*}\cal{M}^{\bullet}
    &=\BB_{\dR}^{\dag,+}\widehat{\otimes}_{\sigma^{-1}\BB_{\dR}^{*,+}}
      \sigma^{-1}\left(
        \left(\lambda^{-1}\Dcap\widehat{\otimes}_{k_{0}}\BB_{\dR}^{*,+}\right)
        \widehat{\otimes}_{\lambda^{-1}\Dcap}^{\rL}\lambda^{-1}\cal{M}^{\bullet}
        \right) \\
    &\cong\sigma^{-1}\left(
        \left(\lambda^{-1}\Dcap\widehat{\otimes}_{k_{0}}\BB_{\dR}^{*,+}\right)
        \widehat{\otimes}_{\lambda^{-1}\Dcap}^{\rL}\lambda^{-1}\cal{M}^{\bullet}
        \right) \\
    &=\left(\nu^{-1}\Dcap\widehat{\otimes}_{k_{0}}\BB_{\dR}^{\dag,+}\right)
      \widehat{\otimes}_{\lambda^{-1}\Dcap}^{\rL}\nu^{-1}\cal{M}^{\bullet} \\
    &\cong\rL\nu^{*}\cal{M}^{\bullet}.
  \end{split}
  \end{equation}
  Here, we used that $\sigma^{-1}\BB_{\dR}^{*,+}\to\BB_{\dR}^{\dag,+}$
  is an isomorphism,~\cite[Lemma 3.32]{Bo21}
  and that $\sigma^{-1}$ is strongly monoidal~\cite[Lemma 2.28]{Bo21}.
  Next, consider the evaluation morphism
  \begin{equation*}
    \rL\lambda^{*}\mathcal{M}^{\bullet}
    \widehat{\otimes}_{\BB_{\dR}^{*,+}}^{\rL}
    \R\shHom_{\nu^{-1}\Dcap\widehat{\otimes}_{k_{0}}\BB_{\dR}^{*,+}}\left(
      \rL\lambda^{*}\mathcal{M},\OB_{\dR}^{*,+}\right)
    \to\OB_{\dR}^{*,+}.
  \end{equation*}
  Now apply $\rL\sigma^{*}$ and consider the composition
  \begin{equation}\label{eq:construct:rhocalMbullet-somemap-reconstructionpaper}
  \begin{split}
    &\rL\sigma^{*}\rL\lambda^{*}\mathcal{M}^{\bullet}
    \widehat{\otimes}_{\rL\sigma^{*}\BB_{\dR}^{*,+}}^{\rL}
    \rL\sigma^{*}
    \R\shHom_{\nu^{-1}\Dcap\widehat{\otimes}_{k_{0}}\BB_{\dR}^{*,+}}\left(
      \rL\lambda^{*}\mathcal{M},\OB_{\dR}^{*,+}\right) \\
    &\to
    \rL\sigma^{*}\left(\rL\lambda^{*}\mathcal{M}^{\bullet}
    \widehat{\otimes}_{\BB_{\dR}^{*,+}}^{\rL}
    \R\shHom_{\nu^{-1}\Dcap\widehat{\otimes}_{k_{0}}\BB_{\dR}^{*,+}}\left(
      \rL\lambda^{*}\mathcal{M},\OB_{\dR}^{*,+}\right)\right) \\    
    &\to\rL\sigma^{*}\OB_{\dR}^{*,+} \\
    &\to\OB_{\pdR}^{\dag}
  \end{split}
  \end{equation}
  Thanks to~(\ref{eq:construct:rhocalMbullet-someiso-reconstructionpaper})
  and since $\rL\sigma^{*}\BB_{\dR}^{*,+}\cong\BB_{\dR}^{\dag,+}$,
  this map~(\ref{eq:construct:rhocalMbullet-somemap-reconstructionpaper}) may be written as
  \begin{equation*}
    \rL\nu^{*}\mathcal{M}^{\bullet}
    \widehat{\otimes}_{\BB_{\dR}^{\dag,+}}^{\rL}
    \rL\sigma^{*}
    \R\shHom_{\nu^{-1}\Dcap\widehat{\otimes}_{k_{0}}\BB_{\dR}^{*,+}}\left(
      \rL\lambda^{*}\mathcal{M},\OB_{\dR}^{*,+}\right) \\
    \to\OB_{\pdR}^{\dag}.
  \end{equation*}
  Produce
  \begin{equation*}
    \rL\nu^{*}\mathcal{M}^{\bullet}
    \to
    \R\shHom_{\BB_{\dR}^{\dag,+}}\left(
    \rL\sigma^{*}
    \R\shHom_{\nu^{-1}\Dcap\widehat{\otimes}_{k_{0}}\BB_{\dR}^{*,+}}\left(
      \rL\lambda^{*}\mathcal{M},\OB_{\dR}^{*,+}\right),
    \OB_{\pdR}^{\dag}\right)
  \end{equation*}
  via the tensor-hom adjunction.
  \begin{align*}
    \mathcal{M}^{\bullet}
    &\to
    \rL\nu^{*}\R\shHom_{\BB_{\dR}^{\dag,+}}\left(
    \rL\sigma^{*}
    \R\shHom_{\nu^{-1}\Dcap\widehat{\otimes}_{k_{0}}\BB_{\dR}^{*,+}}\left(
      \rL\lambda^{*}\mathcal{M},\OB_{\dR}^{*,+}\right),
    \OB_{\pdR}^{\dag}\right) \\
    &\to
    \rL\nu^{*}\R\shHom_{\BB_{\dR}^{\dag,+}}\left(
    \rL\sigma^{*}
    \R\shHom_{\lambda^{-1}\Dcap}\left(
      \lambda^{-1}\mathcal{M},\OB_{\dR}^{*,+}\right),
    \OB_{\pdR}^{\dag}\right) \\
    &=\Rec\left( \Sol\left(\mathcal{M}^{\bullet}\right)\right).
  \end{align*}
  is its adjunct. It is by construction canonical and functorial in $\cal{M}^{\bullet}$.
\end{construction}

\begin{notation}
  For every complex $\cal{M}^{\bullet}$ of $\I\left(\Dcap\right)$-module objects,
  denote $\cal{M}^{\bullet}\to\Rec\left(\Sol\left(\cal{M}^{\bullet}\right)\right)$
  as in Construction~\ref{construct:rhocalMbullet-reconstructionpaper}
  by $\rho_{\cal{M}^{\bullet}}$.
\end{notation}

\begin{thm}\label{thm:reconstruction-thm-solfunctor-reconstructionpaper}
  For every $\mathcal{C}$-complex $\mathcal{M}^{\bullet}$,
  the canonical morphism
  \begin{equation*}
    \rho_{\cal{M}^{\bullet}}\colon\mathcal{M}^{\bullet}
    \stackrel{\cong}{\longrightarrow} \Rec\left(\Sol\left( \mathcal{M}^{\bullet} \right) \right)
  \end{equation*}
  is an isomorphism in the derived category of $\I\left(\Dcap\right)$ module objects.
\end{thm}

\begin{remark}\label{remark:fullyfaithful-is-generalised-rectheorem-recpaper}
  We even prove a generalisation
  of Theorem~\ref{thm:reconstruction-thm-solfunctor-reconstructionpaper},
  namely Theorem~\ref{thm:solfunctorfullyfaithfull-reconstructionpaper-localversion}.
  Indeed, Theorem~\ref{thm:solfunctorfullyfaithfull-reconstructionpaper-localversion}
  for $\cal{N}^{\bullet}=\I\left(\Dcap\right)$ recovers Theorem~\ref{thm:reconstruction-thm-solfunctor-reconstructionpaper}.
  However, we decided to include a separate
  proof for Theorem~\ref{thm:reconstruction-thm-solfunctor-reconstructionpaper},
  as we believe that this makes our arguments more readable.
\end{remark}

See \S\ref{subsec:proofsofmainthms-reconstructionpaper}
for the proof of Theorem~\ref{thm:reconstruction-thm-solfunctor-reconstructionpaper}.

\begin{cor}\label{cor:padic-ishimuraprosmansschneiders}
  The $\I\left(\nu^{-1}\Dcap\right)$-$\I\left(\BB_{\pdR}^{\dag}\right)$-bimodule structure
  on $\I\left(\OB_{\pdR}^{\dag}\right)$ induces
  \begin{equation}\label{eq:derivedIshimura-Dcap-reconstructionpaper}
    \I\left(\Dcap\right)\isomap
    \R\nu_{*}\R\shHom_{\I\left(\BB_{\pdR}^{\dag}\right)}\left(
    \I\left(\OB_{\pdR}^{\dag}\right),\I\left(\OB_{\pdR}^{\dag}\right)\right).
  \end{equation}
\end{cor}

We present the proof of Corollary~\ref{cor:padic-ishimuraprosmansschneiders}
on page~\pageref{proof:cor:padic-ishimuraprosmansschneiders} below.
It uses Lemma~\ref{lem:2--cor:padic-ishimuraprosmansschneiders}
and~\ref{lem:1--cor:padic-ishimuraprosmansschneiders}.

\begin{lem}\label{lem:2--cor:padic-ishimuraprosmansschneiders}
  We have the canonical isomorphism
  \begin{equation*}
    \rL\sigma^{*}\I\left(\OB_{\dR}^{*,+}\right)\isomap\I\left(\OB_{\dR}^{\dag,+}\right).
  \end{equation*}
\end{lem}

\begin{proof}
  Recall that $\sigma$ is the morphism of ringed sites
  \begin{equation*}
    \sigma\colon\left(X_{\proet},\I\left(\BB_{\dR}^{\dag,+}\right)\right)
    \to \left(X_{\proet}\times\NN_{\gg0}^{\op},\I\left(\BB_{\dR}^{*,+}\right)\right).
  \end{equation*}
  The structural map $\sigma^{-1}\I\left(\BB_{\dR}^{*,+}\right)\to\I\left(\BB_{\dR}^{\dag,+}\right)$ is an isomorphism,
  cf. the proof of Lemma~\ref{lem:solexplicitdescription-MIC}.
  $\rL\sigma^{*}\simeq\sigma^{-1}$ follows. Together with~\cite[Lemma 3.39(ii)]{Bo21},
  we find that it suffices to check that the canonical map
  $\sigma^{-1}\OB_{\dR}^{*,+}\to\OB_{\dR}^{\dag,+}$ is an isomorphism.
  This follows from the computation
  \begin{equation*}
    \varinjlim_{\sigma^{-1}(V)\leftarrow U}\OB_{\dR}^{*,+}(V)
    =\varinjlim_{q\in\NN_{\geq2}}\OB_{\dR}^{*,+}\left(\left(U,q\right)\right)
    =\varinjlim_{q\in\NN_{\geq2}}\OB_{\dR}^{q,+}\left(U\right)
    =\OB_{\dR}^{\dag,+}\left(U\right)
  \end{equation*}
  for all $U\in X_{\proet}$; the last step in this computation furthermore
  requires that $U$ is affinoid perfectoid. Indeed, by~\cite[Proposition 3.12(i)]{Sch13pAdicHodge},
  we can compute the colimit on the sections over such $U$.
\end{proof}

\begin{lem}\label{lem:1--cor:padic-ishimuraprosmansschneiders}
  We have the canonical isomorphism of complexes of
  $\I\left(\BB_{\pdR}^{\dag}\right)$-module objects
  \begin{equation*}
    \I\left(\BB_{\pdR}^{\dag}\right)
      \widehat{\otimes}_{\I\left(\BB_{\dR}^{\dag,+}\right)}^{\rL}
      \I\left(\OB_{\dR}^{\dag,+}\right)
    \isomap\I\left(\OB_{\pdR}^{\dag}\right).
  \end{equation*}
\end{lem}

\begin{proof}
  It suffices to check that the restriction of the canonical map to the covering $\widetilde{X}\to X$ as in
  \cite[Corollary~\ref*{cor:localdescription--OBpdRdag-sheafy-reconstructionpaper}]{WiersigPeriods}
  is an isomorphism. Using~\cite[Lemma 3.33]{Bo21}, we compute
  in degree zero
  \begin{align*}
    \I\left(\OB_{\pdR}^{\dag}\right)|_{\widetilde{X}}
    &=\I\left(\OB_{\pdR}^{\dag}|_{\widetilde{X}}\right) \\
    &\stackrel{\text{\ref{cor:localdescription--OBpdRdag-sheafy-reconstructionpaper}}}{\cong}
    \I\left(
      \BB_{\pdR}^{\dag}|_{\widetilde{X}}\widehat{\otimes}_{k_{0}}k_{0}\left\<\frac{Z_{1},\dots,Z_{d}}{p^{q}}\right\>
    \right) \\
    &\cong
    \I\left(\BB_{\pdR}^{\dag}|_{\widetilde{X}}\right)
    \widehat{\otimes}_{\I\left(k_{0}\right)}\I\left(k_{0}\left\<\frac{Z_{1},\dots,Z_{d}}{p^{q}}\right\>\right)
      \\
    &\cong
    \I\left(\BB_{\pdR}^{\dag}|_{\widetilde{X}}\right)
    \widehat{\otimes}_{\I\left(\BB_{\dR}^{\dag,+}\right)|_{\widetilde{X}}}
    \I\left(\BB_{\dR}^{\dag,+}\right)|_{\widetilde{X}}
    \widehat{\otimes}_{\I\left(k_{0}\right)}\I\left(k_{0}\left\<\frac{Z_{1},\dots,Z_{d}}{p^{q}}\right\>\right)\\
    &\cong
    \I\left(\BB_{\pdR}^{\dag}\right)|_{\widetilde{X}}
    \widehat{\otimes}_{\I\left(\BB_{\dR}^{\dag,+}\right)}|_{\widetilde{X}}
    \I\left(\BB_{\dR}^{\dag,+}|_{\widetilde{X}}\right)\widehat{\otimes}_{\I\left(k_{0}\right)}
    \I\left(k_{0}\left\<\frac{Z_{1},\dots,Z_{d}}{p^{q}}\right\>\right) \\
    &\cong
    \I\left(\BB_{\pdR}^{\dag}\right)|_{\widetilde{X}}
    \widehat{\otimes}_{\I\left(\BB_{\dR}^{\dag,+}\right)}|_{\widetilde{X}}
    \I\left(\BB_{\dR}^{\dag,+}|_{\widetilde{X}}\widehat{\otimes}_{k_{0}}
    k_{0}\left\<\frac{Z_{1},\dots,Z_{d}}{p^{q}}\right\>\right) \\
    &\stackrel{\text{\ref{prop:localdescription-of-subsections-OBla-reconstructionpaper}}}{\cong}
    \I\left(\BB_{\pdR}^{\dag}\right)|_{\widetilde{X}}
    \widehat{\otimes}_{\I\left(\BB_{\dR}^{\dag,+}\right)}|_{\widetilde{X}}
    \I\left(\OB_{\dR}^{\dag,+}\right)|_{\widetilde{X}}.
  \end{align*}
  It remains to check that
  $\I\left(\OB_{\pdR}^{\dag}\right)|_{\widetilde{X}}$ is
  flat as an $\I\left(\BB_{\pdR}^{\dag}\right)|_{\widetilde{X}}$-module object.
  We find
  \begin{align*}
    \I\left(\OB_{\pdR}^{\dag,+}\right)|_{\widetilde{X}}
    &=\I\left(\OB_{\pdR}^{\dag}|_{\widetilde{X}}\right) \\
    &\stackrel{\text{\ref{cor:localdescription--OBpdRdag-sheafy-reconstructionpaper}}}{\cong}
    \I\left(
      \BB_{\pdR}^{\dag}|_{\widetilde{X}}\widehat{\otimes}_{k_{0}}k_{0}\left\<\frac{Z_{1},\dots,Z_{d}}{p^{q}}\right\>
    \right) \\
    &\cong
    \I\left(\BB_{\pdR}^{\dag}\right)|_{\widetilde{X}}
    \widehat{\otimes}_{\I\left(k_{0}\right)}\I\left(k_{0}\left\<\frac{Z_{1},\dots,Z_{d}}{p^{q}}\right\>\right),
  \end{align*}
  again using~\cite[Lemma 3.33]{Bo21}.
  Thus Lemma~\ref{lem:constantsheaf-Tatealgebra-LH-exact} implies the flatness.
\end{proof}

\begin{proof}[Proof of Corollary~\ref{cor:padic-ishimuraprosmansschneiders}]
  \label{proof:cor:padic-ishimuraprosmansschneiders}
  Using~\cite[Proposition 3.36]{Bo21}, we compute
  \begin{equation}\label{eq:maincomputation---cor:padic-ishimuraprosmansschneiders}
  \begin{split}
    \Rec\left(\Sol\left(\I\left(\Dcap\right)\right)\right)
    &\cong\Rec\left(\rL\sigma^{*}\I\left(\OB_{\dR}^{*,+}\right)\right) \\
    &\stackrel{\text{\ref{lem:2--cor:padic-ishimuraprosmansschneiders}}}{\cong}
      \Rec\left(\I\left(\OB_{\dR}^{\dag,+}\right)\right) \\
    &\cong\R\nu_{*}\R\shHom_{\I\left(\BB_{\dR}^{\dag,+}\right)}\left(
      \I\left(\OB_{\dR}^{\dag,+}\right),\I\left(\OB_{\pdR}^{\dag}\right)\right) \\
    &\stackrel{\text{\ref{lem:1--cor:padic-ishimuraprosmansschneiders}}}{\cong}
      \R\nu_{*}\R\shHom_{\I\left(\BB_{\pdR}^{\dag}\right)}\left(
      \I\left(\OB_{\dR}^{\dag,+}\right)\widehat{\otimes}_{\I\left(\BB_{\dR}^{\dag,+}\right)}^{\rL}\I\left(\BB_{\pdR}^{\dag}\right)
      \I\left(\OB_{\pdR}^{\dag}\right)\right) \\
    &\cong\R\nu_{*}\R\shHom_{\I\left(\BB_{\pdR}^{\dag}\right)}\left(
      \I\left(\OB_{\pdR}^{\dag}\right),
      \I\left(\OB_{\pdR}^{\dag}\right)\right).
  \end{split}
  \end{equation}
  From~(\ref{eq:maincomputation---cor:padic-ishimuraprosmansschneiders}),
  we observe that~(\ref{eq:derivedIshimura-Dcap-reconstructionpaper}) coincides with $\rho_{\Dcap}$.
  Thus Theorem~\ref{thm:reconstruction-thm-solfunctor-reconstructionpaper} gives the result.
\end{proof}

We state another consequence of Theorem~\ref{thm:reconstruction-thm-solfunctor-reconstructionpaper}.
Here, we use $\D\left(\BB_{\dR}^{\dag,+}\right)\simeq\D\left(\I\left(\BB_{\dR}^{\dag,+}\right)\right)$.

\begin{cor}\label{cor:solfunctor-fullyfaithful-reconstructionpaper}
  The restriction of the solution functor to the category of
  $\mathcal{C}$-complexes
  \begin{equation*}
    \D_{\mathcal{C}}\left( \I\left(\Dcap\right) \right)^{\op}
    \hookrightarrow \D\left( \I\left( \BB_{\dR}^{\dag,+} \right) \right),
    \;
    \cal{M}^{\bullet}\mapsto\Sol\left(\cal{M}^{\bullet}\right)
  \end{equation*}
  is faithful but not full in general.
\end{cor}

\begin{proof}
  The faithfulness is an immediate consequence of Theorem~\ref{thm:reconstruction-thm-solfunctor-reconstructionpaper}.
  We will show that $\Sol$ is not full when $X=\Sp k$ is a point. In this case,
  $\Dcap=k$, and
  \begin{equation*}
    \Hom_{\I(k)}\left(\I(k),\I(k)[1]\right)=\Ext_{\I(k)}^{1}\left(\I(k),\I(k)\right)=0.
  \end{equation*}
  On the other hand,
  \begin{equation*}
  \begin{split}
    \Hom_{\I\left( \BB_{\dR}^{\dag,+}\right)}\left(\Sol\left(\I(k)[1]\right),\Sol\left(\I(k)\right)\right)
    &=\Hom_{\I\left( \BB_{\dR}^{\dag,+}\right)}\left(\I\left( \BB_{\dR}^{\dag,+}\right)[-1],\I\left( \BB_{\dR}^{\dag,+}\right)\right) \\
    &=\Ext_{\I\left( \BB_{\dR}^{\dag,+}\right)}^{1}\left(\I\left( \BB_{\dR}^{\dag,+}\right),\I\left( \BB_{\dR}^{\dag,+}\right)\right) \\
    &=\Ho^{1}\left( X , \I\left( \BB_{\dR}^{\dag,+}\right)\right).
  \end{split}
  \end{equation*}
  Using Theorem~\ref{thm:galois-cohomology-of-solidBdRdaggerplus-born},
  one can easily deduce the non-vanishing of $\Ho^{1}\left( X , \I\left( \BB_{\dR}^{\dag,+}\right)\right)$.
  Therefore
  \begin{equation*}
    \Hom_{\I(k)}\left(\I(k),\I(k)[1]\right)
    \to\Hom_{\I\left( \BB_{\dR}^{\dag,+}\right)}\left(\Sol\left(\I(k)[1]\right),\Sol\left(\I(k)\right)\right)
  \end{equation*}
  is not surjective.
\end{proof}

We still get a fully faithful functor as follows:

\begin{thm}\label{thm:pdRdagsolfunctor-fullyfaithful-reconstructionpaper}
  The following functor on the category of $\mathcal{C}$-complexes
  \begin{equation*}
    \D_{\mathcal{C}}\left( \I\left( \Dcap \right) \right)^{\op} \hookrightarrow \D\left( \I\left( \BB_{\pdR}^{\dag} \right) \right),
    \;
    \cal{M}^{\bullet}\mapsto\SolB_{\pdR}^{\dag}\left(\cal{M}^{\bullet}\right)
  \end{equation*}
  is a fully faithful embedding of triangulated categories.
\end{thm}

For the proof of Theorem~\ref{thm:pdRdagsolfunctor-fullyfaithful-reconstructionpaper},
we refer the reader to \S\ref{subsec:proofofsolfunctorfullyfaithfull-reconstructionpaper},
page~\pageref{proof:--thm:pdRdagsolfunctor-fullyfaithful-reconstructionpaper}.

In the remainder of \S\ref{subsec:thms-recpaper}, we discuss covariant versions.
Recall that $d$ is the dimension of $X$.

\begin{defn}
  The \emph{Riemann-Hilbert functor} is
  \begin{equation*}
    \RH\colon\D\left( \I\left(\BB_{\dR}^{\dag,+}\right)\right) \to \D\left( \I\left(\Dcap\right) \right),
    \mathcal{F}^{\bullet} \mapsto \DD\left(\Rec\left(\mathcal{F}^{\bullet}\right)\right)\left[-d\right].
  \end{equation*}
\end{defn}

\begin{cor}\label{cor:reconstruction-thm-dRfunctor-reconstructionpaper}
  For every $\mathcal{C}$-complex $\mathcal{M}^{\bullet}$,
  $\rho_{\cal{M}^{\bullet}}$ induces a canonical morphism
  \begin{equation*}
    \mathcal{M}^{\bullet}
    \stackrel{\cong}{\longrightarrow} \RH\left(\dRfunctor\left( \mathcal{M}^{\bullet} \right) \right)
  \end{equation*}
  which is an isomorphism in the derived category of $\I\left(\Dcap\right)$-module objects.
\end{cor}

\begin{proof}
  This follows from
  Theorem~\ref{thm:reconstruction-thm-solfunctor-reconstructionpaper}
  as the canonical maps $\mathcal{M}^{\bullet}
    \to\DD^{2}\left(\mathcal{M}^{\bullet}\right)$
  are isomorphisms by~\cite[Theorem 9.17]{Bo21}.
\end{proof}

\begin{cor}\label{cor:deRhamfunctor-fullyfaithful-reconstructionpaper}
  The restriction of the de Rham functor to the category of
  $\mathcal{C}$-complexes
  \begin{equation*}
    \D_{\mathcal{C}}\left( \I\left(\Dcap\right) \right)
    \to \D\left( \I\left(\BB_{\dR}^{\dag,+}\right) \right),\;
    \cal{M}^{\bullet}\mapsto\dRfunctor\left(\cal{M}^{\bullet}\right)
  \end{equation*}
  is faithful but not full in general.
\end{cor}

\begin{proof}
  The duality functor $\DD$ is fully faithful on the category of
  $\cal{C}$-complexes by~\cite[Theorem 9.17]{Bo21}.
  Therefore, Corollary~\ref{cor:deRhamfunctor-fullyfaithful-reconstructionpaper}
  follows from the definition of the de Rham functor and
  Theorem~\ref{cor:solfunctor-fullyfaithful-reconstructionpaper}.
\end{proof}

\begin{cor}\label{cor:pdRdagdRfunctor-fullyfaithful-reconstructionpaper}
  The following functor on the category of $\mathcal{C}$-complexes
  \begin{equation*}
    \D_{\mathcal{C}}\left( \I\left( \Dcap \right) \right) \hookrightarrow \D\left( \I\left( \BB_{\pdR}^{\dag} \right) \right),
    \;
    \cal{M}^{\bullet}\mapsto\dRfunctorB_{\pdR}^{\dag}\left(\cal{M}^{\bullet}\right)
  \end{equation*}
  is a fully faithful embedding of triangulated categories.
\end{cor}

\begin{proof}
  The duality functor $\DD$ is fully faithful on the category of
  $\cal{C}$-complexes by~\cite[Theorem 9.17]{Bo21}.
  Therefore, Theorem~\ref{cor:pdRdagdRfunctor-fullyfaithful-reconstructionpaper}
  follows from the definition of the de Rham functor and
  Theorem~\ref{thm:pdRdagsolfunctor-fullyfaithful-reconstructionpaper}.
\end{proof}


\subsection{Proof of the Reconstruction Theorem}
\label{subsec:proofsofmainthms-reconstructionpaper}

Now we prove Theorem~\ref{thm:reconstruction-thm-solfunctor-reconstructionpaper}.

Theorem~\ref{thm:reconstruction-thm-solfunctor-reconstructionpaper} is a local statement. Therefore,
we can safely assume that $X$ is affinoid and equipped with an étale morphism $g\colon X\to\TT^{d}$.

We give an overview of \S\ref{subsec:proofsofmainthms-reconstructionpaper}.
Our proof follows Prosmans-Schneiders' ideas in complex geometry~\cite{PS00}.
At the heart of their arguments is the isomorphism
\begin{equation*}
  \rho=\rho_{\cal{D}^{\infty}}\colon\cal{D}^{\infty}\isomap\R\shHom_{\CC}\left(\cal{O},\cal{O}\right),
\end{equation*}
from which they deduce their
reconstruction theorem for bounded perfect complex of sheaves of $\cal{D}^{\infty}$-modules.
However, $\cal{C}$-complexes are not bounded perfect complex of sheaves of $\Dcap$-modules.
Instead, at least locally on $X$, we may think of them as homotopy limits of bounded perfect
complexes of sheaves of $\cal{D}_{r}$-modules, cf. the notation as in Lemma~\ref{lem:C-complexdefn-reconstructionpaper}.
Therefore, it does not suffice to establish one analog of $\rho$. Instead,
we need countably many $p$-adic versions which we denote by $\varrho_{r}$.

Interestingly, $\varrho_{r}$ can be parametrised in two different ways: firstly, we have $r$
as in Lemma~\ref{lem:C-complexdefn-reconstructionpaper}. On the other hand, we could parametrise
it by $q$ as in \S\ref{sec:Recap-Periods-Rec}. For technical reasons, we choose
to write $\varrho_{q}$ for the $p$-adic versions of Prosmans-Schneiders' $\rho$.

The first instance of $\varrho_{q}$ appears
as~(\ref{eq:nu-1DrPDXX-BBdRgplus-bimoduleaction-on-OBdRqt-reconstructionpaper})
in \S\ref{subsubsec:subsec:proofsofmainthms-reconstructionpaper-DrPD}.
In \S\ref{subsubsec:subsec:proofsofmainthms-reconstructionpaper-DrPD-to-E},
we check that the $\varrho_{q}$ are isomorphisms in cohomological degree zero.
It remains to check the vanishing of the higher ext, which we do in
\S\ref{subsubsec:subsec:proofsofmainthms-reconstructionpaper-vanishingExt}.
This completes the proof that the $\varrho_{q}$ are isomorphisms.

It remains to deduce Theorem~\ref{subsec:proofsofmainthms-reconstructionpaper}.
Since the $\varrho_{q}$ are isomorphisms, we get reconstruction theorems
for bounded perfect complexes of sheaves of $\cal{D}_{q}$-modules, locally on $X$,
cf. Lemma~\ref{prop:construction-thm-finite-level}.
In \S\ref{subsubsec:subsec:proofsofmainthms-reconstructionpaper-recgloablsection},
we explain how to glue these reconstruction theorems to deduce a
version of Theorem~\ref{thm:reconstruction-thm-solfunctor-reconstructionpaper}
on the level of global sections; this is Theorem~\ref{thm:recthm-globalsections-recpaper}.
Here, we still assume that $X$ is affinoid and equipped with
an étale morphism $X\to\TT^{d}$. Working on the level
of sections is necessary, as this allows us to invoke~\cite[Lemma 8.11]{Bo21}.
There is no analog of \emph{loc. cit.} on the level of sheaves, because the
$\cal{D}_{q}$ all live on different sites $X_{q}$.

$\cal{C}$-complexes are determined by their global sections~\cite[Theorem 8.12]{Bo21},
at least locally on $X$. We are thus able to deduce
Theorem~\ref{subsec:proofsofmainthms-reconstructionpaper} in
\S\ref{subsubsec:subsec:proofsofmainthms-reconstructionpaper-endgame}
from the aforementioned computations.


\subsubsection{Rings of differential operators and divided powers}\label{subsubsec:subsec:proofsofmainthms-reconstructionpaper-DrPD}

This \S\ref{subsubsec:subsec:proofsofmainthms-reconstructionpaper-DrPD}
follows the discussion in~\cite[\S 3.1]{AW2024_GlobalsectionsoneqlinebundlesonOmega}.
Fix the notation from \S\ref{subsec:Dcap-modules}.
We continue to assume that $X$ is affinoid and equipped with an étale morphism
$g\colon X\to\TT^{d}$. Fix the notation $\cal{O}\left(\TT^{d}\right)=k\left\<T_{1}^{\pm},\dots,T_{d}^{\pm}\right\>$.
The \emph{ring of algebraic differential operators} on $X$ is
\begin{equation*}
  \cal{D}(X):=U\left(\Der_{k}\cal{O}(X)\right)\cong\cal{O}(X)\left[\partial_{1},\dots,\partial_{d}\right]
\end{equation*}
Here, $\partial_{i}$ is the lift of the canonical vector field $d/dT^{i}\in\Der_{k}\left(\TT^{d}\right)$.

\begin{notation}
  For $\alpha=\left(\alpha_{1},\dots,\alpha_{d}\right)\in\NN^{d}$, write
  $|\alpha|:=\alpha_{1}+\dots+\alpha_{d}$,
  $\alpha!:=\alpha_{1}!\cdots\alpha_{d}!$,
  $\partial^{\alpha}:=\partial^{\alpha_{1}}\cdots\partial^{\alpha_{d}}$, and
  $\partial^{[\alpha]}=\partial^{\alpha}/\alpha!$.
\end{notation}

The following computation is classical, thus we omit its proof.

\begin{lem}\label{lem:DX-PD-relations-reconstructionpaper}
  The following relations hold in $\cal{D}(X)$:
  \begin{itemize}
    \item[(i)] For all $\alpha,\alpha^{\prime}\in\NN^{d}$,
      \begin{equation*}
        \partial^{\left[\alpha\right]}\partial^{\left[\alpha^{\prime}\right]}
        =\binom{\alpha+\alpha^{\prime}}{\alpha}\partial^{\left[\alpha+\alpha^{\prime}\right]}.
      \end{equation*}
    \item[(ii)] For all $\alpha\in\NN^{d}$ and $f\in\cal{O}(X)$,
      \begin{equation*}
        \partial^{\left[\alpha\right]}f
        =\sum_{\alpha^{\prime}+\alpha^{\prime\prime}=\alpha}
          \partial^{\left[\alpha^{\prime}\right]}(f)\partial^{\left[\alpha^{\prime\prime}\right]}.
      \end{equation*}
  \end{itemize}
\end{lem}


\begin{defn}
  For all $r\in\NN$,
  $\cal{O}(X)^{\circ}\left[p^{r}\partial_{1},\dots,p^{r}\partial_{d}\right]^{\PD}$
  is the $\cal{O}(X)^{\circ}$-subalgebra of $\cal{D}(X)$ generated by the
  set $\left\{ p^{|\alpha|r}\partial^{[\alpha]} \colon \alpha\in\NN^{d} \right\}$.
\end{defn}

See \S\ref{subsec:Dcap-modules}
for the definition of $r_{g}\in\NN$.

\begin{lem}\label{lem:basis-algebraPDdiff-reconstructionpaper}
  If $r\geq r_{g}+1$, then $\left\{ p^{|\alpha|r}\partial^{[\alpha]} \colon \alpha\in\NN^{d} \right\}$
  is a free generating set for
  the left $\cal{O}(X)^{\circ}$-module
  $\cal{O}(X)^{\circ}\left[p^{r}\partial_{1},\dots,p^{r}\partial_{d}\right]^{\PD}$.
\end{lem}

\begin{proof}
  The $\cal{O}(X)^{\circ}$-linear independence is clear. It remains to check that
  the $\cal{O}(X)^{\circ}$-submodule of $\cal{D}(X)$ generated by
  $\left\{ p^{|\alpha|r}\partial^{[\alpha]} \colon \alpha\in\NN^{d} \right\}$ is a subring.
  In light of Lemma~\ref{lem:DX-PD-relations-reconstructionpaper}, it remains to check the following:
  for all $f\in\cal{O}(X)$ and $\alpha^{\prime}\in\NN^{d}$,
  \begin{equation*}
    p^{r|\alpha^{\prime}|}\partial^{\left[\alpha^{\prime}\right]}(f)\in\cal{O}(X)^{\circ}.
  \end{equation*}
  It suffices to check this for $r=r_{g}+1$, where the following computation gives the result:
  \begin{equation*}
    \|p^{(r_{g}+1)|\alpha^{\prime}|}\partial^{\left[\alpha^{\prime}\right]}(f)\|
    =\|p^{r_{g}|\alpha^{\prime}|}\partial^{\alpha^{\prime}}(f)\|
      |p^{|\alpha^{\prime}|}/\alpha!|
    \leq 1\cdot 1=1.
  \end{equation*}
  We used $p^{r_{g}|\alpha^{\prime}|}\partial^{\alpha^{\prime}}(f)\in\cal{O}(X)^{\circ}$,
  which follows from the definition of $r_{g}$
  in~\S\ref{subsec:Dcap-modules}.
\end{proof}

\begin{defn}
  For all $r\geq r_{g}+1$, define the $k$-Banach algebra
  \begin{equation*}
    \cal{D}_{r}^{\PD}(X):=
    \widehat{\cal{O}(X)^{\circ}\left[p^{r}\partial_{1},\dots,p^{r}\partial_{d}\right]^{\PD}}[1/p].
  \end{equation*}
  Here, the completion is the $p$-adic one.
\end{defn}

\begin{lem}\label{lem:descriptioncalDrPDX-reconstructionpaper}
  If $r\geq r_{g}+1$, then
  \begin{equation*}
    \cal{D}_{r}^{\PD}(X)=\left\{
      \sum_{\alpha\in\NN^{d}}f_{\alpha}\partial^{[\alpha]}\in\prod_{\alpha\in\NN^{d}}\cal{O}(X)\partial^{[\alpha]}
        \colon
        p^{-r|\alpha|}f_{\alpha}\to0\text{ for }|\alpha|\to\infty
    \right\},
  \end{equation*}
  equipped with the norm
  $\|\sum_{\alpha\in\NN^{d}}f_{\alpha}\partial^{[\alpha]}\|
    =\sup_{\alpha\in\NN^{d}}p^{r|\alpha|}\|f_{\alpha}\|$.
\end{lem}

\begin{proof}
  This is a direct consequence of Lemma~\ref{lem:basis-algebraPDdiff-reconstructionpaper}.
\end{proof}

\begin{remark}\label{rem:DcapX-to-calDrPDX-reconstructionpaper}
  From Lemma~\ref{lem:descriptioncalDrPDX-reconstructionpaper}
  and~\cite[Lemma 3.4]{KBB18}, we find a canonical morphism
  \begin{equation}\label{eq:DcapX-to-calDrPDX-reconstructionpaper}
    \Dcap(X)\to\cal{D}_{r}^{\PD}(X)
  \end{equation}  
  of $k$-Fréchet spaces. This induces a canonical morphism between
  the associated $k$-ind-Banach spaces.
\end{remark}

\begin{notation}\label{notation:constantsheaf-reconstructionpaper}
  Given a $k$-ind-Banach space $V$, denote the associated constant sheaf on $X$ by $V_{X}$.
\end{notation}

Fix $q\gg0$ large enough so that we have the $\nu^{-1}\Dcap$-$\BB_{\dR}^{q,+}$-bimodule
structure on $\OB_{\dR}^{q,+}$ as constructed in Theorem~\ref{thm:OBlaplus-bimodule}.
This makes
$\OB_{\dR}^{q,+}$ a $\nu^{-1}\Dcap(X)_{X}$-$\BB_{\dR}^{q,+}$-bimodule object, which is equivalent
to a map
\begin{equation}\label{eq:nu-1DcapXX-BBdRgplus-bimoduleaction-on-OBdRqt-reconstructionpaper}
  \nu^{-1}\Dcap(X)_{X}
  \to\shHom_{\BB_{\dR}^{q,+}}\left(\OB_{\dR}^{q,+},\OB_{\dR}^{q,+}\right)
\end{equation}
of sheaves of $k$-in-Banach algebras.

\begin{lem}\label{lem:nu-1DrPDXX-BBdRgplus-bimoduleaction-on-OBdRqt-reconstructionpaper}
  Given natural numbers $r\geq r_{g}+1$, $r\geq q\gg0$,
  (\ref{eq:nu-1DcapXX-BBdRgplus-bimoduleaction-on-OBdRqt-reconstructionpaper})
  factors uniquely through a morphism
  \begin{equation}\label{eq:nu-1DrPDXX-BBdRgplus-bimoduleaction-on-OBdRqt-reconstructionpaper}
    \nu^{-1}\cal{D}_{r}^{\PD}(X)_{X}
    \to\shHom_{\BB_{\dR}^{q,+}}\left(\OB_{\dR}^{q,+},\OB_{\dR}^{q,+}\right)
  \end{equation}
  of sheaves of $k$-ind-Banach algebras.
\end{lem}

\begin{proof}
  We start with the existence.
  Recall that $X$ is affinoid and equipped with an étale morphism
  $g\colon X\to\TT^{d}$. We may therefore consider the affinoid perfectoid pro-étale covering
  $\widetilde{X}\to X$. For any $U\in X_{\proet}/\widetilde{X}$,
  (\ref{eq:nu-1DcapXX-BBdRgplus-bimoduleaction-on-OBdRqt-reconstructionpaper})
  induces a morphism
  \begin{equation}
    \Dcap(X)
    \to\shHom_{\BB_{\dR}^{q,+}}\left(\OB_{\dR}^{q,+},\OB_{\dR}^{q,+}\right)(U).
  \end{equation}
  We consider its restriction
  \begin{equation*}
    \varrho\colon
    \cal{O}(X)^{\circ}\left[p^{r}\partial_{1},\dots,p^{r}\partial_{d}\right]^{\PD}
    \to\shHom_{\BB_{\dR}^{q,+}}\left(\OB_{\dR}^{q,+},\OB_{\dR}^{q,+}\right)(U).
  \end{equation*}
  By Remark~\ref{remark:OBdRqplus-bimodulestructureproof-overview}, this is a map of
  $\cal{O}(X)^{\circ}$-algebras determined by $\partial_{i}\mapsto d/dZ_{i}$
  for all $i=1,\dots,d$.
  It suffices to check that it extends by continuity to a morphism
  \begin{equation*}
    \widehat{\cal{O}^{\circ}\left[p^{r}\partial_{1},\dots,p^{r}\partial_{d}\right]^{\PD}}
    \to\shHom_{\BB_{\dR}^{q,+}}\left(\OB_{\dR}^{q,+},\OB_{\dR}^{q,+}\right)(U)
  \end{equation*}
  of $k$-ind-Banach algebras. To do this, we compute
  \begin{equation}\label{eq:nu-1DrPDXX-BBdRgplus-bimoduleaction-on-OBdRqt-somecomputation-reconstructionpaper}
  \begin{split}
    \shHom_{\BB_{\dR}^{q,+}}\left(\OB_{\dR}^{q,+},\OB_{\dR}^{q,+}\right)(U)
    &\cong
      \shHom_{\BB_{\dR}^{q,+}}\left(\BB_{\dR}^{q,+}\left\<\frac{Z_{1},\dots,Z_{d}}{p^{q}}\right\>,\OB_{\dR}^{q,+}\right)(U) \\
    &\stackrel{\clubsuit}{\cong}
      \shHom_{k}\left(k\left\<\frac{Z_{1},\dots,Z_{d}}{p^{q}}\right\>_{X},\OB_{\dR}^{q,+}\right)(U) \\
    &\cong
      \intHom_{k}\left(k\left\<\frac{Z_{1},\dots,Z_{d}}{p^{q}}\right\>,\OB_{\dR}^{q,+}(U)\right) \\
    &\stackrel{\clubsuit}{\cong}
      \intHom_{k}\left(k\left\<\frac{Z_{1},\dots,Z_{d}}{p^{q}}\right\>,\BB_{\dR}^{q,+}(U)\left\<\frac{Z_{1},\dots,Z_{d}}{p^{q}}\right\>\right).
  \end{split}
  \end{equation}
  Here we used the~\cite[equation before Proposition 2.2.10]{Sch99},
  and $\clubsuit$ follows
  from~\cite[Proposition~\ref*{prop:localdescription-of-subsections-OBla-reconstructionpaper}]{WiersigPeriods}
  $\shHom_{\BB_{\dR}^{q,+}}\left(\OB_{\dR}^{q,+},\OB_{\dR}^{q,+}\right)(U)$
  is thus a Banach space. Therefore, we may check that the image of
  $\varrho$ lies in its unit ball. In symbols:
  \begin{equation*}
    \varrho\left(P\right)\in\shHom_{\BB_{\dR}^{q,+}}\left(\OB_{\dR}^{q,+},\OB_{\dR}^{q,+}\right)(U)^{\circ}
  \end{equation*}
  for every $P\in\cal{O}^{\circ}\left[p^{r}\partial_{1},\dots,p^{r}\partial_{d}\right]^{\PD}$.
  With~(\ref{eq:nu-1DrPDXX-BBdRgplus-bimoduleaction-on-OBdRqt-somecomputation-reconstructionpaper}),
  this translates into the following:
  \begin{equation}\label{lem:nu-1DrPDXX-BBdRgplus-bimoduleaction-on-OBdRqt-todoforexistence-reconstructionpaper}
    \|\varrho(P)(f)\|\leq\|f\| \text{ for all } f\in k\left\<\frac{Z_{1},\dots,Z_{d}}{p^{q}}\right\>.
  \end{equation}
  This is what we check. Without loss of generality, $P=p^{r|\alpha|}\partial^{[\alpha]}$ for $\alpha\in\NN^{d}$.
  Then use $r\geq q$ to compute for every $f=\sum_{\beta\in\NN^{d}}\lambda_{\beta}Z^{\beta}\in k\left\<\frac{Z_{1},\dots,Z_{d}}{p^{q}}\right\>$,
  \begin{equation*}
    \|\varrho(P)(f)\|
    =\|p^{r|\alpha|}\sum_{\beta\geq\alpha}\lambda_{\alpha}\binom{\beta}{\alpha}Z^{\beta-\alpha}\|
    =|p^{r|\alpha|}|\sup_{\beta\geq\alpha}|p^{-q|\beta-\alpha|}|\|\lambda_{\alpha}\binom{\beta}{\alpha}\|
    \leq \sup_{\beta\in\NN^{d}}|p^{q|\beta|}|\|\lambda_{\alpha}\|
    =\|f\|.
  \end{equation*}
  We have thus
  checked~(\ref{lem:nu-1DrPDXX-BBdRgplus-bimoduleaction-on-OBdRqt-todoforexistence-reconstructionpaper}).
  As explained above, this gives the existence
  of~(\ref{eq:nu-1DrPDXX-BBdRgplus-bimoduleaction-on-OBdRqt-reconstructionpaper}).
  
  It remains to check uniqueness. This would follow immediately once we have checked that
  the canonical morphism
  \begin{equation*}
    \nu^{-1}\Dcap(X)_{X}\to\nu^{-1}\cal{D}_{r}^{\PD}(X)_{X},
  \end{equation*}
  induced by~(\ref{eq:DcapX-to-calDrPDX-reconstructionpaper}),
  is an epimorphism. Since $\nu^{-1}$ is a left adjoint and because
  sheafification is strongly exact, cf.~\cite[Lemma 2.11]{Sch99},
  it suffices to check that
  $\Dcap(X)\to\cal{D}_{r}^{\PD}(X)$ is an epimorphism.
  This follows because its image is dense, which is implied
  by Lemma~\ref{lem:descriptioncalDrPDX-reconstructionpaper}.
\end{proof}

\begin{cor}\label{cor:OBdRqplus-bimodulestructureover-calDrPDXXBdRqplus}
  Given natural numbers $r\geq r_{g}+1$, $r\geq q\gg0$,
  $\OB_{\dR}^{q,+}$ admits a unqiue $\nu^{-1}\cal{D}_{r}^{\PD}(X)_{X}$-$\BB_{\dR}^{q,+}$-bimodule
  structure, which is compatible with the
  bimodule structure
  in Corollary~\ref{cor:OBdRdagplus-bimodule}.
\end{cor}

\begin{proof}
  This is a consequence of Lemma~\ref{lem:nu-1DrPDXX-BBdRgplus-bimoduleaction-on-OBdRqt-reconstructionpaper}.
\end{proof}


\subsubsection{Differential operators as bounded linear endomorphisms}\label{subsubsec:subsec:proofsofmainthms-reconstructionpaper-DrPD-to-E}

We continue to fix $q\gg0$ as in~\S\ref{subsubsec:subsec:proofsofmainthms-reconstructionpaper-DrPD}.
Set $r:=q$.
The global sections of $\shHom_{\BB_{\dR}^{q,+}}\left(\OB_{\dR}^{q,+},\OB_{\dR}^{q,+}\right)$ are
\begin{equation*}
  \intHom_{\BB_{\dR}^{q,+}}\left(\OB_{\dR}^{q,+},\OB_{\dR}^{q,+}\right).
\end{equation*}
Invoking abstract nonsense,~(\ref{eq:nu-1DrPDXX-BBdRgplus-bimoduleaction-on-OBdRqt-reconstructionpaper})
induces a morphism
\begin{equation*}
  \varrho_{q}^{0}\colon\cal{D}_{q}^{\PD}(X)
  \to\intHom_{\BB_{\dR}^{q,+}}\left(\OB_{\dR}^{q,+},\OB_{\dR}^{q,+}\right).
\end{equation*}
In this \S\ref{subsubsec:subsec:proofsofmainthms-reconstructionpaper-DrPD-to-E},
we prove the following result.

\begin{prop}\label{prop:nu-1DrPDXX-BBdRgplus-bimoduleaction-on-OBdRqt-globalsections-iso-reconstructionpaper}
  $\varrho_{q}^{0}$ is an isomorphism.
\end{prop}

We complete the proof of
Proposition~\ref{prop:nu-1DrPDXX-BBdRgplus-bimoduleaction-on-OBdRqt-globalsections-iso-reconstructionpaper} on
page~\pageref{proof:prop:nu-1DrPDXX-BBdRgplus-bimoduleaction-on-OBdRqt-globalsections-iso-reconstructionpaper}
and start with some notation.

\begin{notation}
  Let $\iota$ denote the map
  $\cal{O}\left(X\right)\to\OB_{\dR}^{q,+}\left(\widetilde{X}_{C}\right)$ induced
  by the canonical $\nu^{-1}\cal{O}\to\OB_{\dR}^{q,+}$.
  On the other hand, we use~\cite[Corollary~\ref*{cor:localdescription-of-subsections-OBqplus}]{WiersigPeriods}
  to write down the morphism
  \begin{equation*}
    \Ovartheta\colon\OB_{\dR}^{q,+}\left(\widetilde{X}_{C}\right)
    \cong
    \BB_{\dR}^{q,+}\left(\widetilde{X}_{C}\right)\left\<\frac{Z_{1},\dots,Z_{d}}{p^{q}}\right\>
    \to\widehat{\cal{O}}\left(\widetilde{X}_{C}\right),
      \sum_{\beta\in\NN^{d}}b_{\beta}Z^{\beta}\mapsto\vartheta_{\dR}^{q,+}\left(b_{0}\right).
  \end{equation*}
  We refer the reader to~\cite[\S\ref*{subsubsection:invertingp-solutionpaper}]{WiersigPeriods}
  for the definition of $\vartheta_{\dR}^{q,+}$.
\end{notation}

\begin{lem}\label{lem:cd:nu-1DrPDXX-BBdRgplus-bimoduleaction-on-OBdRqt-globalsections-injective-reconstructionpaper}
  We have the commutative diagram
  \begin{equation}\label{cd:nu-1DrPDXX-BBdRgplus-bimoduleaction-on-OBdRqt-globalsections-injective-reconstructionpaper}
  \begin{tikzcd}
    \cal{O}(X)
      \arrow{r}{\iota}\arrow{rd} &
    \OB_{\dR}^{q,+}\left(\widetilde{X}_{C}\right)
      \arrow{d}{\Ovartheta} \\
    \empty &
    \widehat{\cal{O}}\left(\widetilde{X}_{C}\right).
  \end{tikzcd}
  \end{equation}
\end{lem}

\begin{proof}
  $\iota$
  comes from $\nu^{-1}\cal{O}\to\OB_{\dR}^{q,+}$, which we constructed
  via~\cite[Lemma~\ref*{lem:imagetildenu-Ainf}]{WiersigPeriods}.
  The second half of \emph{loc. cit.} gives the result.
\end{proof}

\begin{lem}\label{lem:varrhoqiso-iotainjective-reconstructionpaper}
  $\iota\colon\cal{O}\left(X\right)\to\OB_{\dR}^{q,+}\left(\widetilde{X}_{C}\right)$
  is injective.
\end{lem}

\begin{proof}
  Consider~(\ref{cd:nu-1DrPDXX-BBdRgplus-bimoduleaction-on-OBdRqt-globalsections-injective-reconstructionpaper}).
  The diagonal map is injective, as it is
  $\cal{O}(X)\cong\widehat{\cal{O}}(X)\to\widehat{\cal{O}}\left(\widetilde{X}_{C}\right)$,
  cf.~\cite[Corollary 6.19]{Sch13pAdicHodge}, and
  $\widetilde{X}_{C}\to X$ is a covering.
  Lemma~\ref{lem:cd:nu-1DrPDXX-BBdRgplus-bimoduleaction-on-OBdRqt-globalsections-injective-reconstructionpaper}
  implies that $\iota$ is injective.
\end{proof}

Next, we note that the codomain of $\varrho_{q}$
is more concrete than it first appears.

\begin{lem}\label{lem:intHoBBdRqplusOBdRqplusOBdRqplus-Banach-reconstructionpaper}
  $\intHom_{\BB_{\dR}^{q,+}}\left(\OB_{\dR}^{q,+},\OB_{\dR}^{q,+}\right)$
  is a $k$-Banach space.
\end{lem}

\begin{proof}
  Use the~\cite[equation before Proposition 2.2.10]{Sch99}
  and~\cite[Proposition~\ref*{prop:localdescription-of-subsections-OBla-reconstructionpaper}]{WiersigPeriods}:
  \begin{equation}\label{eq:lem:intHoBBdRqplusOBdRqplusOBdRqplus-Banach-reconstructionpaper}
  \begin{split}
    \intHom_{\BB_{\dR}^{q,+}}\left(\OB_{\dR}^{q,+},\OB_{\dR}^{q,+}\right)
    &=\shHom_{\BB_{\dR}^{q,+}}\left(\OB_{\dR}^{q,+},\OB_{\dR}^{q,+}\right)(X) \\
    &=\ker\left(
      \begin{array}{l}
        \shHom_{\BB_{\dR}^{q,+}}\left(\OB_{\dR}^{q,+},\OB_{\dR}^{q,+}\right)\left(\widetilde{X}_{C}\right) \\
        \to\shHom_{\BB_{\dR}^{q,+}}\left(\OB_{\dR}^{q,+},\OB_{\dR}^{q,+}\right)\left(\widetilde{X}_{C}\times_{X}\widetilde{X}_{C}\right)
      \end{array}
    \right) \\
    &\cong
      \ker\left(
      \begin{array}{l}
        \shHom_{k_{0}}\left(k_{0}\left\<\frac{Z_{1},\dots,Z_{d}}{p^{q}}\right\>,\OB_{\dR}^{q,+}\right)\left(\widetilde{X}_{C}\right) \\
        \to\shHom_{k_{0}}\left(k_{0}\left\<\frac{Z_{1},\dots,Z_{d}}{p^{q}}\right\>,\OB_{\dR}^{q,+}\right)\left(\widetilde{X}_{C}\times_{X}\widetilde{X}_{C}\right)
      \end{array}
    \right) \\
    &\cong
      \ker\left(
      \begin{array}{l}
        \intHom_{k_{0}}\left(k_{0}\left\<\frac{Z_{1},\dots,Z_{d}}{p^{q}}\right\>,\OB_{\dR}^{q,+}\left(\widetilde{X}_{C}\right)\right) \\
        \to\intHom_{k_{0}}\left(k_{0}\left\<\frac{Z_{1},\dots,Z_{d}}{p^{q}}\right\>,\OB_{\dR}^{q,+}\left(\widetilde{X}_{C}\times_{X}\widetilde{X}_{C}\right)\right)
      \end{array}
    \right).
  \end{split}
  \end{equation}
  Lemma~\ref{lem:intHoBBdRqplusOBdRqplusOBdRqplus-Banach-reconstructionpaper} follows because
  kernels of morphisms of Banach spaces are Banach spaces.
\end{proof}

In the following lemma, we consider the adjunct of the canonical
morphism $\nu^{-1}\cal{O}\to\OB_{\dR}^{q,+}$.

\begin{lem}\label{lem:globalsections-OBdRqplus-reconstructionpaper}
  The canonical map $\cal{O}\to\nu_{*}\OB_{\dR}^{q,+}$ induces
  the isomorphism
  $\cal{O}(X)\isomap\OB_{\dR}^{q,+}(X)$.
\end{lem}

\begin{proof}
  Consider the commutative diagram
  \begin{equation*}
  \begin{tikzcd}
    \cal{O}(X) \arrow{r}\arrow{rd} &
    \OB_{\dR}^{q,+}(X) \arrow{d} \\
    \empty &
    \OB_{\dR}^{\dag}(X).
  \end{tikzcd}
  \end{equation*}
  The diagonal morphism is an isomorphism by Theorem~\ref{thm:borncohomologyOBdRdagX-reconstructionpaper}.
  Thanks to~\cite[Proposition 1.1.8]{Sch99}, it remains to check that
  $\OB_{\dR}^{q,+}(X)\to\OB_{\dR}^{\dag}(X)$ is a monomorphism.
  This follows because $\OB_{\dR}^{q,+}\to\OB_{\dR}^{\dag}$ is a monomorphism,
  which can be checked locally and is thus implied
  by~\cite[Corollary~\ref*{cor:localdescription-of-subsections-OBqplus}
     Theorem~\ref*{thm:cohomology-OBdRdag-over-affperfd-reconstructionpaper}]{WiersigPeriods},
  and because
  $\BB_{\dR}^{q,+}\left(R,R^{+}\right)\to\BB_{\dR}^{\dag}\left(R,R^{+}\right)$
  is a monomorphism for any affinoid perfectoid algebra $\left(R,R^{+}\right)$,
  cf.~\cite[the proof of Theorem~\ref*{thm:BdRdagRRplus-bornological}
  on page~\pageref*{proof:thm:BdRdagRRplus-bornological}]{WiersigPeriods}.
\end{proof}

We start with the construction of a two-sided inverse of $\varrho_{q}$.
Let $\phi\mapsto\phi(X)$ denote
\begin{align*}
  \intHom_{\BB_{\dR}^{q,+}}\left(\OB_{\dR}^{q,+},\OB_{\dR}^{q,+}\right)
  &\to\intHom_{\BB_{\dR}^{q,+}\left(X\right)}
    \left(\OB_{\dR}^{q,+}\left(X\right),\OB_{\dR}^{q,+}\left(X\right)\right)\\
  \empty &\quad\stackrel{\text{\ref{lem:globalsections-OBdRqplus-reconstructionpaper}}}{\cong}
    \intHom_{\BB_{\dR}^{q,+}\left(X\right)}
    \left(\cal{O}\left(X\right),\cal{O}\left(X\right)\right),
\end{align*}
the canonical projection coming from Lemma~\ref{lem:sheafinternalhom-over-monoid}.
The following is inspired by~\cite[Definition 4.5]{KBB18}.

\begin{defn}\label{defn:etaalpha-reconstructionpaper}
  For $\phi\in\intHom_{\BB_{\dR}^{q,+}}\left(\OB_{\dR}^{q,+},\OB_{\dR}^{q,+}\right)$
  (see Lemma~\ref{lem:intHoBBdRqplusOBdRqplusOBdRqplus-Banach-reconstructionpaper})
  and $\alpha\in\NN^{d}$,
  \begin{equation*}
    \eta_{[\alpha]}\left(\phi\right)
    :=\sum_{\beta\leq\alpha}\phi(X)\left(T^{\beta}\right)\binom{\alpha}{\beta}\left(-T\right)^{\alpha-\beta}
      \in\cal{O}(X).
  \end{equation*}
\end{defn}

Here is our version of~\cite[Theorem 4.9]{KBB18}:

\begin{lem}\label{lem:alpha!etaalpha-converges-reconstructionpaper}
  Given $\phi\in\intHom_{\BB_{\dR}^{q,+}}\left(\OB_{\dR}^{q,+},\OB_{\dR}^{q,+}\right)$,
  $p^{q|\alpha|}\eta_{[\alpha]}\left(\phi\right)\to0$ for $\alpha\to\infty$.
\end{lem}

\begin{proof}
  Again from Lemma~\ref{lem:sheafinternalhom-over-monoid}, we get the canonical projection
  \begin{equation*}
    \phi\left(\widetilde{X}_{C}\right)\in
    \intHom_{\BB_{\dR}^{q,+}\left(\widetilde{X}_{C}\right)}
      \left(\OB_{\dR}^{q,+}\left(\widetilde{X}_{C}\right),\OB_{\dR}^{q,+}\left(\widetilde{X}_{C}\right)\right) 
  \end{equation*}
  of $\phi$.
  \emph{Loc. cit.} also gives $\phi\left(\widetilde{X}_{C}\right)\circ\iota=\iota\circ\phi(X)$.
  Compse with $\Ovartheta$ on both sides to get
  \begin{equation}\label{eq:OvarthetaphitXiota-is-phi-reconstructionpaper}
    \Ovartheta\circ\phi\left(\widetilde{X}_{C}\right)\circ\iota
    =\left(\Ovartheta\circ\iota\right)\circ\phi(X)
    \stackrel{\text{\ref{lem:cd:nu-1DrPDXX-BBdRgplus-bimoduleaction-on-OBdRqt-globalsections-injective-reconstructionpaper}}}{=}\phi(X).
  \end{equation}
  Next, we define for all $\alpha\in\NN^{d}$,
  \begin{align*}
    \widetilde{\eta}_{[\alpha]}\left(\phi\right)
    :=\phi\left(\widetilde{X}_{C}\right)\left(Z^{\alpha}\right)
      \in\OB_{\dR}^{q,+}\left(\widetilde{X}_{C}\right).
  \end{align*}
  Using $Z_{i}=\iota\left(T_{i}\right)-\left[T_{i}^{\flat}\right]$
  and $\Ovartheta\left(\left[T_{i}^{\flat}\right]\right)=\vartheta_{\dR}^{q,+}\left(\left[T_{i}^{\flat}\right]\right)=T_{i}$
  for all $i=1,\dots,d$, one computes
  \begin{equation*}
  \begin{split}
    \Ovartheta\left(\widetilde{\eta}_{[\alpha]}\left(\phi\right)\right)
    &=\Ovartheta\left(\phi\left(\widetilde{X}_{C}\right)\left(Z^{\alpha}\right)\right) \\
    &=\Ovartheta\left(\phi\left(\widetilde{X}_{C}\right)\left(
      \sum_{\beta\leq\alpha}
        \iota\left( T \right)^{\beta}
        \binom{\alpha}{\beta}
        \left(-\left[T^{\flat}\right]\right)^{\alpha-\beta}\right)
      \right) \\
    &=\Ovartheta\left(
      \sum_{\beta\leq\alpha}
        \phi\left(\widetilde{X}_{C}\right)\left(\iota\left( T \right)^{\beta}\right)
        \binom{\alpha}{\beta}
        \left(-\left[T^{\flat}\right]\right)^{\alpha-\beta}
      \right) \\
    &=
      \sum_{\beta\leq\alpha}
        \Ovartheta\left(\phi\left(\widetilde{X}_{C}\right)\left(\iota\left( T \right)^{\beta}\right)\right)
        \binom{\alpha}{\beta}
        \left(-T\right)^{\alpha-\beta} \\
    &\stackrel{\text{(\ref{eq:OvarthetaphitXiota-is-phi-reconstructionpaper})}}{=}
      \sum_{\beta\leq\alpha}
        \phi\left(X\right)\left(T^{\beta}\right)
        \binom{\alpha}{\beta}
        \left(-T\right)^{\alpha-\beta}.
  \end{split}
  \end{equation*}
  That is, for all $\alpha\in\NN^{d}$,
  $\Ovartheta\left(\widetilde{\eta}_{[\alpha]}\left(\phi\right)\right)
    =\eta_{[\alpha]}\left(\phi\right)$.
  This implies
  \begin{equation*}
    p^{q|\alpha|}\eta_{[\alpha]}\left(\phi(X)\right)
    =\Ovartheta\left(\phi\left(\widetilde{X}_{C}\right)\left(p^{q|\alpha|}Z^{\alpha}\right)\right)
    \to0
  \end{equation*}
  for $\alpha\to\infty$, because
  $p^{q|\alpha|}Z^{\alpha}\to0$ for $\alpha\to\infty$ in $\OB_{\dR}^{q,+}\left(\widetilde{X}_{C}\right)$.
\end{proof}

Lemma~\ref{lem:descriptioncalDrPDX-reconstructionpaper}
and~\ref{lem:alpha!etaalpha-converges-reconstructionpaper}
allow us to define
\begin{equation*}
  \eta^{q}\colon
  \intHom_{\BB_{\dR}^{q,+}}\left(\OB_{\dR}^{q,+},\OB_{\dR}^{q,+}\right)
  \to
  \cal{D}_{r}^{\PD}(X),
  \phi\mapsto\sum_{\alpha\in\NN^{d}}\eta_{[\alpha]}\left(\phi(X)\right)\partial^{[\alpha]}.
\end{equation*}

\begin{lem}\label{lem:etaqcircrhoq-identity-reconstructionpaper}
  $\eta^{q}\circ\varrho_{q}^{0}$ is the identity on $\cal{D}_{r}^{\PD}(X)$.
\end{lem}

\begin{proof}
  Given $P=\sum_{\alpha\in\NN^{d}}f_{\alpha}\partial^{[\alpha]}\in\cal{D}_{r}^{\PD}(X)$,
  set $\phi:=\varrho_{q}^{0}(P)$. We have to check
  \begin{equation*}
    \eta_{[\alpha]}\left(\phi(X)\right) = f_{\alpha}
  \end{equation*}
  for all $\alpha\in\NN^{d}$.
  Proceed by following the arguments in
  the~\cite[proof of Lemma 4.11]{KBB18}.
\end{proof}

\begin{proof}[Proof of Proposition~\ref{prop:nu-1DrPDXX-BBdRgplus-bimoduleaction-on-OBdRqt-globalsections-iso-reconstructionpaper}]
\label{proof:prop:nu-1DrPDXX-BBdRgplus-bimoduleaction-on-OBdRqt-globalsections-iso-reconstructionpaper}
  Thanks to Lemma~\ref{lem:intHoBBdRqplusOBdRqplusOBdRqplus-Banach-reconstructionpaper},
  we have the open mapping theorem at our disposal. Therefore,
  it suffices to check that $\varrho_{q}^{0}$ is bijective. It is injective by
  Lemma~\ref{lem:etaqcircrhoq-identity-reconstructionpaper}.
  Furthermore,
  \begin{equation*}
    \eta^{q}\circ\left(\varrho_{q}^{0}\circ\eta^{q}-\id\right)
    =\left(\eta^{q}\circ\varrho_{q}^{0}-\id\right)\circ\eta^{q}=0,
  \end{equation*}
  again by Lemma~\ref{lem:etaqcircrhoq-identity-reconstructionpaper}.
  Therefore, it suffices to check that $\eta^{q}$ is injective. To show this,
  consider $\phi\in\ker\eta^{q}$. Following the arguments in
  the~\cite[first paragraph of the proof of Proposition 4.13]{KBB18},
  we find $\phi(X)\left(T^{\beta}\right)=0$ for all $\beta\in\NN^{d}$.
  This implies, for all $\beta\in\NN^{d}$,
  \begin{equation*}
    \phi\left(\widetilde{X}_{C}\right)\left( Z^{\beta} \right)
    =\sum_{\alpha\leq\beta}
      \phi\left(\widetilde{X}_{C}\right)\left( \iota\left(T^{\alpha}\right) \right)\left(-\left[T^{\flat}\right]\right)^{\beta-\alpha}
    =\sum_{\alpha\leq\beta}
      \iota\left(\phi\left(X\right)\left( T^{\alpha}\right) \right)\left(-\left[T^{\flat}\right]\right)^{\beta-\alpha}
    =0
  \end{equation*}
  From~\cite[Theorem~\ref*{thm:localdescription-of-OBqplus}]{WiersigPeriods},
  we deduce $\phi|_{\widetilde{X}_{C}}=0$.
  Because $\widetilde{X}_{C}\to X$ is a covering, $\phi=0$.
\end{proof}


\subsubsection{Vanishing Ext}\label{subsubsec:subsec:proofsofmainthms-reconstructionpaper-vanishingExt}

We follow the notation in \S\ref{subsubsec:subsec:proofsofmainthms-reconstructionpaper-DrPD-to-E} and
keep $q\gg0$, $r=q$ fixed. We need a derived version of
Proposition~\ref{prop:nu-1DrPDXX-BBdRgplus-bimoduleaction-on-OBdRqt-globalsections-iso-reconstructionpaper}.
To ensure the homological algebra runs smoothly, we therefore work
in the left heart $\IndBan_{\I\left(k\right)}:=\LH\left(\IndBan_{k}\right)$.
As before, $\I$ denotes the canonical embedding $\IndBan_{k}\to\IndBan_{\I\left(k\right)}$.
$\I$ is lax monoidal by~\cite[Proposition 2.12]{2024BambozziKremnizerSheafynessproperty}.
Therefore, Corollary~\ref{cor:OBdRqplus-bimodulestructureover-calDrPDXXBdRqplus}
equippes $\I\left(\OB_{\dR}^{q,+}\right)$ with the structure of a
$\I\left(\nu^{-1}\cal{D}_{q}^{\PD}(X)_{X}\right)$-$\I\left(\BB_{\dR}^{q,+}\right)$-bimodule object.
This is equivalent to the data of a morphism
\begin{equation*}
  \I\left(\nu^{-1}\cal{D}_{q}^{\PD}(X)_{X}\right)
  \to\shHom_{\I\left(\BB_{\dR}^{q,+}\right)}\left(\I\left(\OB_{\dR}^{q,+}\right),\I\left(\OB_{\dR}^{q,+}\right)\right)
\end{equation*}
of $\IndBan_{\left(k\right)}$-sheaves of rings. Taking global sections induces a morphism
\begin{equation}\label{eq:prop:derivedRq-iso-constructionmap-reconstructionpaper}
  \I\left(\cal{D}_{q}^{\PD}(X)\right)
  \to\intHom_{\I\left(\BB_{\dR}^{q,+}\right)}\left(\I\left(\OB_{\dR}^{q,+}\right),\I\left(\OB_{\dR}^{q,+}\right)\right).
\end{equation}
Derive the internal homomorphisms as in~\cite[\S 4.2]{Bo21} and
compose~(\ref{eq:prop:derivedRq-iso-constructionmap-reconstructionpaper}) with
\begin{equation*}
  \intHom_{\I\left(\BB_{\dR}^{q,+}\right)}\left(\I\left(\OB_{\dR}^{q,+}\right),\I\left(\OB_{\dR}^{q,+}\right)\right)
    \to\R\intHom_{\I\left(\BB_{\dR}^{q,+}\right)}\left(\I\left(\OB_{\dR}^{q,+}\right),\I\left(\OB_{\dR}^{q,+}\right)\right)
\end{equation*}
to get the morphism
\begin{equation}\label{eq:almostvarrho-but-not-iso-recpaper}
  \I\left(\cal{D}_{q}^{\PD}(X)\right)
  \to\R\intHom_{\I\left(\BB_{\dR}^{q,+}\right)}\left(\I\left(\OB_{\dR}^{q,+}\right),\I\left(\OB_{\dR}^{q,+}\right)\right).
\end{equation}
We do not expect~(\ref{eq:almostvarrho-but-not-iso-recpaper}) to be
an isomorphism.
Indeed, Theorem~\ref{thm:galois-cohomology-of-solidBdRdaggerplus-born}
suggests that $\I\left(\OB_{\dR}^{q,+}\right)$ has non-vanishing higher cohomology.
Therefore, we compose~(\ref{eq:almostvarrho-but-not-iso-recpaper}) with the canonical morphism
\begin{equation*}
  \R\intHom_{\I\left(\BB_{\dR}^{q,+}\right)}\left(\I\left(\OB_{\dR}^{q,+}\right),\I\left(\OB_{\dR}^{q,+}\right)\right)
  \to\R\intHom_{\I\left(\BB_{\dR}^{q,+}\right)}\left(\I\left(\OB_{\dR}^{q,+}\right),\I\left(\OB_{\pdR}^{\dag}\right)\right)
\end{equation*}
in the following Proposition~\ref{prop:derivedRq-iso-reconstructionpaper}.

\begin{prop}\label{prop:derivedRq-iso-reconstructionpaper}
  The morphism
  \begin{equation*}
    \varrho_{q}\colon
    \I\left(\cal{D}_{q}^{\PD}(X)\right)
    \isomap\R\intHom_{\I\left(\BB_{\dR}^{q,+}\right)}\left(\I\left(\OB_{\dR}^{q,+}\right),\I\left(\OB_{\pdR}^{\dag}\right)\right)
  \end{equation*}
  constructed above is an isomorphism.
\end{prop}

We complete the proof of
Proposition~\ref{prop:derivedRq-iso-reconstructionpaper} on
page~\pageref{proof:prop:derivedRq-iso-reconstructionpaper}
and start with computations in cohomological degree zero.

\begin{lem}\label{lem1:lem:prop:derivedRq-iso-degree0-reconstructionpaper}
  Given any $\I\left(\BB_{\dR}^{q,+}\right)$-module object $\cal{G}$,
  we have the isomorphism
  \begin{equation*}
    \intHom_{\I\left(\BB_{\dR}^{q,+}\right)}\left(\I\left(\OB_{\dR}^{q,+}\right),\cal{G}\right)
    \cong\intHom_{\I\left(k_{0}\right)}\left(\I\left(k_{0}\left\<\frac{Z_{1},\dots,Z_{d}}{p^{q}}\right\>\right),\cal{G}(X)\right)
  \end{equation*}  
  which is functorial in $\cal{G}$.
\end{lem}

\begin{proof}
  Recall Definition~\ref{defn:constantsheaf} and compute
  \begin{equation}\label{eq1:lem1:lem:prop:derivedRq-iso-degree0-reconstructionpaper}
  \begin{split}
    \I\left(\OB_{\dR}^{q,+}\right)|_{\widetilde{X}}
    &\cong
    \I\left(\BB_{\dR}^{q,+}|_{\widetilde{X}}\right)
      \widehat{\otimes}_{\I\left(k_{0}\right)}
      \I\left(k_{0}\left\<\frac{Z_{1},\dots,Z_{d}}{p^{q}}\right\>_{X_{\proet}/\widetilde{X}}\right) \\
    &\stackrel{\text{\ref{lem:I-commutes-constantsheaf-recpaper}}}{\cong}
      \I\left(\BB_{\dR}^{q,+}|_{\widetilde{X}}\right)
      \widehat{\otimes}_{\I\left(k_{0}\right)}
      \I\left(k_{0}\left\<\frac{Z_{1},\dots,Z_{d}}{p^{q}}\right\>\right)_{X_{\proet}/\widetilde{X}},
  \end{split}
  \end{equation}
  where we
  used~\cite[Proposition~\ref*{prop:localdescription-of-subsections-OBla-reconstructionpaper}]{WiersigPeriods}
  and~\cite[Lemma 3.33]{Bo21}.
  Now we find
  \begin{equation*}
  \begin{split}
    &\intHom_{\I\left(\BB_{\dR}^{q,+}\right)}\left(\I\left(\OB_{\dR}^{q,+}\right),\cal{G}\right) \\
    &=\shHom_{\I\left(\BB_{\dR}^{q,+}\right)}\left(\I\left(\OB_{\dR}^{q,+}\right),\cal{G}\right)(X) \\
    &=\ker\left(
      \begin{array}{l}
        \shHom_{\I\left(\BB_{\dR}^{q,+}\right)}\left(\I\left(\OB_{\dR}^{q,+}\right),\cal{G}\right)\left(\widetilde{X}\right) \\
        \to\shHom_{\I\left(\BB_{\dR}^{q,+}\right)}\left(\I\left(\OB_{\dR}^{q,+}\right),\cal{G}\right)\left(\widetilde{X}\times_{X}\widetilde{X}\right)
      \end{array}
    \right) \\
    &\stackrel{\text{(\ref{eq1:lem1:lem:prop:derivedRq-iso-degree0-reconstructionpaper})}}{\cong}
      \ker\left(
      \begin{array}{l}
        \shHom_{\I\left(\BB_{\dR}^{q,+}\right)}\left(
          \I\left(\BB_{\dR}^{q,+}|_{\widetilde{X}}\right)
            \widehat{\otimes}_{\I\left(k_{0}\right)}
            \I\left(k_{0}\left\<\frac{Z_{1},\dots,Z_{d}}{p^{q}}\right\>\right)_{X_{\proet}/\widetilde{X}},
          \cal{G}\right)\left(\widetilde{X}\right) \\
        \to
        \shHom_{\I\left(\BB_{\dR}^{q,+}\right)}\left(
          \I\left(\BB_{\dR}^{q,+}|_{\widetilde{X}}\right)
            \widehat{\otimes}_{\I\left(k_{0}\right)}
            \I\left(k_{0}\left\<\frac{Z_{1},\dots,Z_{d}}{p^{q}}\right\>\right)_{X_{\proet}/\widetilde{X}},
          \cal{G}\right)\left(\widetilde{X}\times_{X}\widetilde{X}\right)
      \end{array}
    \right) \\
    &\cong
      \ker\left(
      \begin{array}{l}
        \shHom_{\I\left(k_{0}\right)}\left(
            \I\left(k_{0}\left\<\frac{Z_{1},\dots,Z_{d}}{p^{q}}\right\>\right)_{X_{\proet}/\widetilde{X}},
          \cal{G}\right)\left(\widetilde{X}\right) \\
        \to
        \shHom_{\I\left(k_{0}\right)}\left(
            \I\left(k_{0}\left\<\frac{Z_{1},\dots,Z_{d}}{p^{q}}\right\>\right)_{X_{\proet}/\widetilde{X}},
          \cal{G}\right)\left(\widetilde{X}\times_{X}\widetilde{X}\right)
      \end{array}
    \right) \\
    &\stackrel{\text{($*$)}}{\cong}
      \ker\left(
      \begin{array}{l}
        \intHom_{\I\left(k_{0}\right)}\left(
            \I\left(k_{0}\left\<\frac{Z_{1},\dots,Z_{d}}{p^{q}}\right\>\right),
          \cal{G}\left(\widetilde{X}\right)\right) \\
        \to
        \intHom_{\I\left(k_{0}\right)}\left(
            \I\left(k_{0}\left\<\frac{Z_{1},\dots,Z_{d}}{p^{q}}\right\>\right),
          \cal{G}\left(\widetilde{X}\times_{X}\widetilde{X}\right)\right)
      \end{array}
    \right) \\
    &\cong\intHom_{\I\left(k_{0}\right)}\left(
      \I\left(k_{0}\left\<\frac{Z_{1},\dots,Z_{d}}{p^{q}}\right\>\right),
      \ker\left(
          \cal{G}\left(\widetilde{X}\right)
          \to
          \cal{G}\left(\widetilde{X}\times_{X}\widetilde{X}\right)
          \right)
        \right) \\
    &\cong\intHom_{\I\left(k_{0}\right)}\left(
      \I\left(k_{0}\left\<\frac{Z_{1},\dots,Z_{d}}{p^{q}}\right\>\right),
        \cal{G}\left(X\right)
        \right),
  \end{split}
  \end{equation*}
  as desired. Here, the isomorphism ($*$) comes from the~\cite[equation before Proposition 2.2.10]{Sch99}.
\end{proof}

\begin{lem}\label{lem2:lem:prop:derivedRq-iso-degree0-reconstructionpaper}
  Given any sheaf $\cal{F}$ of $\BB_{\dR}^{q,+}$-ind-Banach modules,
  we have the isomorphism
  \begin{equation*}
    \intHom_{\I\left(\BB_{\dR}^{q,+}\right)}\left(\I\left(\OB_{\dR}^{q,+}\right),\I\left(\cal{F}\right)\right)
    \cong\I\left(\intHom_{k_{0}}\left(k_{0}\left\<\frac{Z_{1},\dots,Z_{d}}{p^{q}}\right\>,\cal{F}(X)\right)\right)
  \end{equation*}  
  which is functorial in $\cal{F}$.
\end{lem}

\begin{proof}
  Recall that the sections of $\I\left(\cal{F}\right)$ over a $U\in X_{\proet}$
  are given by $\I\left(\cal{F}(U)\right)$.
  \begin{equation*}
  \begin{split}
    \intHom_{\I\left(\BB_{\dR}^{q,+}\right)}\left(\I\left(\OB_{\dR}^{q,+}\right),\I\left(\cal{F}\right)\right)  
    &\cong\intHom_{\I\left(k_{0}\right)}\left(
      \I\left(k_{0}\left\<\frac{Z_{1},\dots,Z_{d}}{p^{q}}\right\>\right),
      \I\left(\cal{F}\left(X\right)\right)
      \right) \\
    &\cong\I\left(\intHom_{k_{0}}\left(k_{0}\left\<\frac{Z_{1},\dots,Z_{d}}{p^{q}}\right\>,\cal{F}(X)\right)\right),
  \end{split}
  \end{equation*}  
  follows with Lemma~\ref{lem1:lem:prop:derivedRq-iso-degree0-reconstructionpaper}
  and~\ref{lem:IHom-Tatealgebra-is-HomITatealgebra-reconstructionpaper}.
\end{proof}

\begin{lem}\label{lem:prop:derivedRq-iso-degree0-reconstructionpaper}
  We have the isomorphism
  \begin{equation*}
    \Ho^{0}\left(\varrho_{q}\right)\colon
    \I\left(\cal{D}_{q}^{\PD}(X)\right)
    \isomap\intHom_{\I\left(\BB_{\dR}^{q,+}\right)}\left(\I\left(\OB_{\dR}^{q,+}\right),\I\left(\OB_{\pdR}^{\dag}\right)\right).
  \end{equation*}
\end{lem}

\begin{proof}
  By construction, $\Ho^{0}\left(\varrho_{q}\right)$ fits into a commutative diagram
  \begin{equation*}
  \begin{tikzcd}
    \I\left(\cal{D}_{q}^{\PD}(X)\right)
      \arrow[swap]{rd}{\Ho^{0}\left(\varrho_{q}\right)} \arrow{r}{\phi} &
    \intHom_{\I\left(\BB_{\dR}^{q,+}\right)}\left(\I\left(\OB_{\dR}^{q,+}\right),\I\left(\OB_{\dR}^{q,+}\right)\right)
      \arrow{d}{\psi} \\
    \empty &
    \intHom_{\I\left(\BB_{\dR}^{q,+}\right)}\left(\I\left(\OB_{\dR}^{q,+}\right),\I\left(\OB_{\pdR}^{\dag}\right)\right).
  \end{tikzcd}
  \end{equation*}
  Lemma~\ref{lem:globalsections-OBdRqplus-reconstructionpaper},
  \cite[Theorem~\ref*{thm:borncohomologyOBpdRdagX-reconstructionpaper}]{WiersigPeriods},
  and Lemma~\ref{lem1:lem:prop:derivedRq-iso-degree0-reconstructionpaper}
  imply that $\psi$ is an isomorphism. It remains to check that $\phi$ is an isomorphism.
  $\phi$ is~(\ref{eq:prop:derivedRq-iso-constructionmap-reconstructionpaper}),
  which we constructed via the $\nu^{-1}\cal{D}_{r}^{\PD}(X)_{X}$-$\BB_{\dR}^{q,+}$-bimodule
  structure on $\OB_{\dR}^{q,+}$. This was similar to the construction of $\varrho_{q}$,
  before
  Proposition~\ref{prop:nu-1DrPDXX-BBdRgplus-bimoduleaction-on-OBdRqt-globalsections-iso-reconstructionpaper}.
  \begin{equation*}
  \begin{split}
    &\intHom_{\I\left(\BB_{\dR}^{q,+}\right)}\left(
      \I\left(\OB_{\dR}^{q,+}\right),
      \I\left(\OB_{\pdR}^{\dag}\right)\right) \\
    &\stackrel{\text{\ref{lem2:lem:prop:derivedRq-iso-degree0-reconstructionpaper}}}{\cong}
      \I\left(\intHom_{k_{0}}\left(k_{0}\left\<\frac{Z_{1},\dots,Z_{d}}{p^{q}}\right\>,\OB_{\dR}^{q,+}(X)\right)\right) \\
    &\cong
      \I\left(\ker\left( 
      \begin{array}{l}
        \intHom_{k_{0}}\left(k_{0}\left\<\frac{Z_{1},\dots,Z_{d}}{p^{q}}\right\>,\OB_{\dR}^{q,+}\left(\widetilde{X}\right)\right) \\
        \to\intHom_{k_{0}}\left(k_{0}\left\<\frac{Z_{1},\dots,Z_{d}}{p^{q}}\right\>,\OB_{\dR}^{q,+}\left(\widetilde{X}\times_{X}\widetilde{X}\right)\right)
       \end{array}\right)\right) \\
     &\stackrel{\text{(\ref{eq:lem:intHoBBdRqplusOBdRqplusOBdRqplus-Banach-reconstructionpaper})}}{\cong}      
       \I\left(\intHom_{\BB_{\dR}^{q,+}}
         \left(
         \OB_{\dR}^{q,+},
         \OB_{\pdR}^{\dag}\right)\right),
  \end{split}
  \end{equation*}
  thus implies $\phi=\I\left(\varrho_{q}^{0}\right)$,
  cf. Proposition~\ref{prop:nu-1DrPDXX-BBdRgplus-bimoduleaction-on-OBdRqt-globalsections-iso-reconstructionpaper}.
  \emph{Loc. cit.} implies that $\phi$ is an isomorphism.
\end{proof}

We continue with computations in higher cohomological degree.

\begin{lem}\label{lem:prop:derivedRq-iso-positivedegree-reconstructionpaper}
  For all $i>0$,
  \begin{equation*}
    \R^{i}\intHom_{\I\left(\BB_{\dR}^{q,+}\right)}\left(\I\left(\OB_{\dR}^{q,+}\right),\I\left(\OB_{\pdR}^{\dag}\right)\right)=0.
  \end{equation*}
\end{lem}

\begin{proof}
  Pick an injective resolution $\I\left(\OB_{\pdR}^{\dag}\right)\to\cal{I}^{\bullet}$.
  We find
  \begin{align*}
    \Ho^{i}\left(\intHom_{\I\left(\BB_{\dR}^{q,+}\right)}\left(\I\left(\OB_{\dR}^{q,+}\right),\cal{I}^{\bullet}\right)\right)
    &\stackrel{\text{\ref{lem1:lem:prop:derivedRq-iso-degree0-reconstructionpaper}}}{\cong}
    \Ho^{i}\left(\intHom_{\I\left(k_{0}\right)}\left(
      \I\left(k_{0}\left\<\frac{Z_{1},\dots,Z_{d}}{p^{q}}\right\>\right),\Gamma\left(X,\cal{I}^{\bullet}\right)\right)\right) \\
    &\stackrel{\text{\ref{lem:Tate-algebra-internally-projective-LHIndBanF-reconstructionpaper}}}{\cong}
    \intHom_{\I\left(k_{0}\right)}\left(
      \I\left(k_{0}\left\<\frac{Z_{1},\dots,Z_{d}}{p^{q}}\right\>\right),\Ho^{i}\left(\Gamma\left(X,\cal{I}^{\bullet}\right)\right)\right) \\
    &\stackrel{\clubsuit}{\cong}
    \intHom_{\I\left(k_{0}\right)}\left(
      \I\left(k_{0}\left\<\frac{Z_{1},\dots,Z_{d}}{p^{q}}\right\>\right),0\right) \\
    &=0
  \end{align*}
  for all $i>0$. The isomorphism $\clubsuit$
  is~\cite[Theorem~\ref*{thm:borncohomologyOBpdRdagX-reconstructionpaper}]{WiersigPeriods}.
\end{proof}

\begin{proof}[Proof of Proposition~\ref{prop:derivedRq-iso-reconstructionpaper}]
  \label{proof:prop:derivedRq-iso-reconstructionpaper}
  Put Lemma~\ref{lem:prop:derivedRq-iso-degree0-reconstructionpaper}
  and~\ref{lem:prop:derivedRq-iso-positivedegree-reconstructionpaper}
  together.
\end{proof}


\subsubsection{A local reconstruction theorem}\label{subsubsec:subsec:proofsofmainthms-reconstructionpaper-recgloablsection}

We continue to assume that $X$ is affinoid and equipped with an étale morphism $g\colon X\to\TT^{d}$.
In \S\ref{subsubsec:subsec:proofsofmainthms-reconstructionpaper-recgloablsection},
we establish a local version of Theorem~\ref{thm:reconstruction-thm-solfunctor-reconstructionpaper}.
It is based on a corresponding adaptation of the notion of a $\cal{C}$-complex.
We do not recall the~\cite[Definition 8.10]{Bo21}, but use the
characterisation as in Lemma~\ref{lem:sections-C-complexdefn-reconstructionpaper} below.
Recall Definition~\ref{defn:perfectcomplex-recpaper}.

\begin{lem}\label{lem:sections-C-complexdefn-reconstructionpaper}
  $M^{\bullet}\in\D\left(\I\left(\Dcap(X)\right)\right)$ is a $\mathcal{C}$-complex if and only if
  \begin{itemize}
    \item[(i)] $M_{r}^{\bullet}
    :=\I\left(\mathcal{D}_{r}(X)\right)\widehat{\otimes}_{\I\left(\Dcap(X)\right)}^{\rL}M^{\bullet}$
    is a bounded perfect complex of
    $\I\left(\mathcal{D}_{r}(X)\right)$-module objects for every $r\geq r_{g}$, and
    \item[(ii)]
    $\Ho^{j}\left( M^{\bullet}\right)\stackrel{\cong}{\longrightarrow}
    \varprojlim_{r\geq r_{g}}\Ho^{j}\left(M_{r}^{\bullet}\right)$
    for every $j\in\ZZ$.
  \end{itemize}
\end{lem}

\begin{proof}
  Proceed as in the proof of Lemma~\ref{lem:C-complexdefn-reconstructionpaper}.
\end{proof}

Recall Notation~\ref{notation:constantsheaf-reconstructionpaper}.

\begin{defn}\label{defn:solfrefunctor-globalsections-recpaper}
  Denote the following functor by $\nSol$:
  \begin{equation*}
    \D\left(\I\left(\Dcap(X)\right)\right)^{\op}
    \to
    \D\left(\I\left(\BB_{\dR}^{\dag,+}\right)\right),
    M^{\bullet}
    \mapsto
    \rL\sigma^{*}\R\shHom_{\I\left(\Dcap(X)\right)}\left(
      M^{\bullet}_{X_{\proet}\times\NN_{\gg0}^{\op}},\I\left(\OB_{\dR}^{*,+}\right)
    \right).
  \end{equation*}
  Denote the following functor by $\nRec$:
  \begin{equation*}
    \D\left(\I\left(\BB_{\dR}^{\dag,+}\right)\right)\to\D\left(\I\left(\Dcap(X)\right)\right)^{\op},
    \cal{F}^{\bullet}
    \mapsto
    \R\Gamma\left(
      X,
      \R\shHom_{\I\left(\BB_{\dR}^{\dag,+}\right)}\left(\cal{F}^{\bullet},\I\left(\OB_{\pdR}^{\dag}\right)\right)
    \right).
  \end{equation*}
\end{defn}

For every complex $M^{\bullet}$ of $\I\left(\Dcap(X)\right)$-module objects,
we follow Construction~\ref{construct:rhocalMbullet-reconstructionpaper}
to find morphisms
\begin{equation*}
  \varrho_{M^{\bullet}}\colon M^{\bullet}\to \nRec\left( \nSol \left( M^{\bullet} \right) \right)
\end{equation*}
which are canonical and functorial in $M^{\bullet}$.
Here is the main result of \S\ref{subsubsec:subsec:proofsofmainthms-reconstructionpaper-recgloablsection}:

\begin{thm}\label{thm:recthm-globalsections-recpaper}
  $\varrho_{M^{\bullet}}\colon M^{\bullet}\to\nRec\left( \nSol \left( M^{\bullet} \right) \right)$
  is an isomorphism when $M^{\bullet}$ is a $\cal{C}$-complex.
\end{thm}

We complete the proof of Theorem~\ref{thm:recthm-globalsections-recpaper}
on page~\pageref{proof:thm:recthm-globalsections-recpaper}.
For now,  we introduce more variants of the solution functor,
cf. Definition~\ref{defn:solfunctordagq-globalsections-recpaper}
and~\ref{defn:solfunctorq-globalsections-recpaper}
below. We start with notation.

\begin{defn}
  For every $q\in\NN$, $\OB_{\dR}^{\dag,q,+}:=\BB_{\dR}^{\dag,+}\widehat{\otimes}_{\BB_{\dR}^{q,+}}\OB_{\dR}^{q,+}$.
\end{defn}

Let $q\gg0$.
We consider the $\Dcap(X)$-$\BB_{\dR}^{q,+}$-bimodule
structure on $\OB_{\dR}^{q,+}$ coming from Theorem~\ref{thm:OBlaplus-bimodule}.
It induces a $\Dcap(X)_{X_{\proet}}$-$\BB_{\dR}^{\dag,+}$-bimodule
structure on $\OB_{\dR}^{\dag,q,+}$. Since $\I$
is strongly monoidal, cf.~\cite[Lemma 3.33]{Bo21}, and thanks
to Lemma~\ref{lem:I-commutes-constantsheaf-recpaper},
this gives a $\I\left(\Dcap(X)\right)_{X_{\proet}}$-$\I\left(\BB_{\dR}^{\dag,+}\right)$-bimodule
structure on $\I\left(\OB_{\dR}^{\dag,q,+}\right)$, which we fix.
It plays a crucial role in the following Definition~\ref{defn:solfunctordagq-globalsections-recpaper}.

\begin{defn}\label{defn:solfunctordagq-globalsections-recpaper}
  For every $q\gg0$,
  denote the following functor by $\nSol^{\dag,q}$:
  \begin{equation*}
    \D\left(\I\left(\Dcap(X)\right)\right)^{\op}
    \to
    \D\left(\I\left(\BB_{\dR}^{\dag,+}\right)\right),
    M^{\bullet}
    \mapsto
    \R\shHom_{\I\left(\Dcap(X)\right)}\left(
      M^{\bullet}_{X_{\proet}},\I\left(\OB_{\dR}^{\dag,q,+}\right)
    \right).
  \end{equation*}
\end{defn}

\begin{construction}\label{construct:map--prop:construction-thm-finite-level}
  In this Construction~\ref{construct:map--prop:construction-thm-finite-level},
  we work in left hearts and omit the symbol $\I$ throughout.
  This notation does not cause imprecisions
  by Lemma~\ref{lem:I-commutes-constantsheaf-recpaper}
  and~\cite[Lemma 3.39]{Bo21}.
  
We continue to fix $q\gg0$ and use
Notation~\ref{notation:constantsheaf-reconstructionpaper}.
Equip $X_{\proet}$ with $\nu^{-1}\Dcap(X)_{X}\widehat{\otimes}_{k_{0}}\BB_{\dR}^{\dag,+}$ and
$X$ with $\Dcap(X)$. The canonical map
\begin{equation*}
  \nu^{-1}\Dcap(X)_{X}
  \cong\nu^{-1}\Dcap(X)_{X}\widehat{\otimes}_{k_{0}}k_{0}
  \to\nu^{-1}\Dcap(X)_{X}\widehat{\otimes}_{k_{0}}\BB_{\dR}^{\dag,+}
\end{equation*}
makes $\nu\colon X_{\proet}\to X$ a morphism of ringed sites. 
For every complex $M^{\bullet}$ of $\Dcap(X)$-ind-Banach modules,
let $M_{X}^{\bullet}$ denote the corresponding complex of sheaves of
$\Dcap(X)_{X}$-ind-Banach modules.
Furthermore, the
$\nu^{-1}\Dcap$-$\BB_{\dR}^{q,+}$-bimodule structure
on $\OB_{\dR}^{q,+}$ as in Corollary~\ref{cor:OBdRqplus-bimodulestructureover-calDrPDXXBdRqplus} induces a
$\nu^{-1}\Dcap(X)_{X}$-$\BB_{\dR}^{\dag,+}$-bimodule structure on $\OB_{\dR}^{\dag,q,+}$.
We use this to write down the following evaluation morphism
\begin{equation*}
  \rL\nu^{*}M_{X}^{\bullet}
  \widehat{\otimes}_{\BB_{\dR}^{\dag,+}}^{\rL}\R\shHom_{\nu^{-1}\Dcap(X)_{X}\widehat{\otimes}_{k_{0}}\BB_{\dR}^{\dag,+}}\left(\rL\nu^{*}M_{X}^{\bullet},\OB_{\dR}^{\dag,q,+}\right)
  \to\OB_{\dR}^{\dag,q,+}.
\end{equation*}
Produce
\begin{equation*}
  \rL\nu^{*}M_{X}^{\bullet}
  \to
  \R\shHom_{\BB_{\dR}^{\dag,+}}\left(\R\shHom_{\nu^{-1}\Dcap(X)_{X}\widehat{\otimes}_{k_{0}}\BB_{\dR}^{\dag,+}}\left(\rL\nu^{*}M_{X}^{\bullet},\OB_{\dR}^{\dag,q,+}\right), \OB_{\dR}^{\dag,q,+}\right)
\end{equation*}
via the tensor-hom adjunction.
This is by construction a morphism of
$\nu^{-1}\Dcap(X)_{X}\widehat{\otimes}_{k_{0}}\BB_{\dR}^{\dag,+}$-ind-Banach modules,
where $\nu^{-1}\Dcap(X)_{X}$ acts on the right-most $\OB_{\dR}^{\dag,q,+}$. It induces
\begin{equation}\label{eq:construction-thm-finite-level-somemaptoconstructthemap}
\begin{split}
  M^{\bullet}
  \to
  &\R\Gamma\left(X,\R\shHom_{\BB_{\dR}^{\dag,+}}\left(\R\shHom_{\nu^{-1}\Dcap(X)_{X}\widehat{\otimes}_{k_{0}}\BB_{\dR}^{\dag,+}}\left(\rL\nu^{*}M_{X}^{\bullet},\OB_{\dR}^{\dag,q,+}\right), \OB_{\dR}^{\dag,q,+}\right)\right) \\
  &\cong\R\Gamma\left(X,\R\shHom_{\BB_{\dR}^{\dag,+}}\left(\R\shHom_{\Dcap(X)}\left(\nu^{-1}M_{X}^{\bullet},\OB_{\dR}^{\dag,q,+}\right), \OB_{\dR}^{\dag,q,+}\right)\right) \\
  &\stackrel{\text{\ref{lem:constantsheaf-commuteswith-inverseimage-recpaper}}}{\cong}\R\Gamma\left(X,\R\shHom_{\BB_{\dR}^{\dag,+}}\left(\R\shHom_{\Dcap(X)}\left(M_{X_{\proet}}^{\bullet},\OB_{\dR}^{\dag,q,+}\right), \OB_{\dR}^{\dag,q,+}\right)\right).
\end{split}
\end{equation}
As $\OB_{\dR}^{\dag,q,+}$ is
a $\nu^{-1}\cal{D}_{q}^{\PD}(X)_{X}\widehat{\otimes}_{k_{0}}\BB_{\dR}^{\dag,+}$-bimodule object,
the morphism~(\ref{eq:construction-thm-finite-level-somemaptoconstructthemap}) lifts to
\begin{equation*}
\begin{split}
  \cal{D}_{q}^{\PD}(X)\widehat{\otimes}_{\cal{D}_{q}(X)}^{\rL}M^{\bullet}
  \to
  \R\Gamma\left(X,\R\shHom_{\BB_{\dR}^{\dag,+}}\left(\R\shHom_{\Dcap(X)}\left(M_{X_{\proet}}^{\bullet},\OB_{\dR}^{\dag,q,+}\right), \OB_{\dR}^{\dag,q,+}\right)\right).
\end{split}
\end{equation*}
Compose this map with the following map, which is induced by the canonical
$\OB_{\dR}^{\dag,q,+}\to\OB_{\pdR}^{\dag}$,
\begin{align*}
    \R\Gamma\left(X,
      \RshHom_{\BB_{\dR}^{\dag,+}}\left(\R\shHom_{\Dcap(X)}\left(
        M_{X_{\proet}}^{\bullet},\OB_{\dR}^{\dag,q,+}\right),
        \OB_{\dR}^{\dag,q,+}\right)\right)
    \to\nRec\left(\nSol^{\dag,q}\left(M^{\bullet}\right)\right)
\end{align*}
to get $\varrho_{M^{\bullet},q}^{\dag}$
in the following Proposition~\ref{prop:construction-thm-finite-level}.
\end{construction}

\begin{prop}\label{prop:construction-thm-finite-level}
  Fix $q\gg0$.
  Then, for every $\cal{C}$-complex $M^{\bullet}\in\D\left(\I\left(\Dcap(X)\right)\right)$,
  \begin{equation*}
    \varrho_{M^{\bullet},q}^{\dag}\colon
    \I\left(\cal{D}_{q}^{\PD}(X)\right)\widehat{\otimes}_{\I\left(\Dcap(X)\right)}^{\rL}M^{\bullet}
    \isomap\nRec\left(\nSol^{\dag,q}\left(M^{\bullet}\right)\right).
  \end{equation*}
\end{prop}

We present the proof of Proposition~\ref{prop:construction-thm-finite-level}
below the following two Lemma~\ref{lem:prop:construction-thm-finite-level}
and~\ref{lem:isoBdRdagBdRqOBdRdagq-LH--lem:prop:construction-thm-finite-level}.

\begin{lem}\label{lem:prop:construction-thm-finite-level}
  For any $\IndBan_{\I\left(k_{0}\right)}$-sheaf $\cal{F}$, we have the equivalence of functors
  \begin{equation*}
    \R\Gamma(X,-)\circ\R\shHom_{\I\left(\BB_{\dR}^{q,+}\right)}\left(\I\left(\OB_{\dR}^{q,+}\right),-\right)
    \simeq\R\intHom_{\I\left(\BB_{\dR}^{q,+}\right)}\left(\I\left(\OB_{\dR}^{q,+}\right),-\right).
  \end{equation*}
\end{lem}

\begin{proof}
  Lemma~\ref{lem:Tate-algebra-internally-projective-LHIndBanF-reconstructionpaper}
  and~\ref{lem1:lem:prop:derivedRq-iso-degree0-reconstructionpaper} imply
  that $\shHom_{\I\left(\BB_{\dR}^{q,+}\right)}\left(\I\left(\OB_{\dR}^{q,+}\right),-\right)$
  sends injective objects to flasque sheaves. Now Lemma~\ref{lem:prop:construction-thm-finite-level}
  follows from
  \begin{equation*}
    \Gamma(X,-)\circ\shHom_{\I\left(\BB_{\dR}^{q,+}\right)}\left(\I\left(\OB_{\dR}^{q,+}\right),-\right)
    =\intHom_{\I\left(\BB_{\dR}^{q,+}\right)}\left(\I\left(\OB_{\dR}^{q,+}\right),-\right),
  \end{equation*}
  together with the associated Grothendieck spectral sequence.
\end{proof}

\begin{lem}\label{lem:isoBdRdagBdRqOBdRdagq-LH--lem:prop:construction-thm-finite-level}
  We have the canonical isomorphism of
  $\I\left(\BB_{\dR}^{q,+}\right)$-module objects
  \begin{equation*}
    \I\left(\BB_{\dR}^{\dag,+}\right)
      \widehat{\otimes}_{\I\left(\BB_{\dR}^{q,+}\right)}^{\rL}
      \I\left(\OB_{\dR}^{q,+}\right)
    \isomap\I\left(\OB_{\dR}^{\dag,q,+}\right).
  \end{equation*}
\end{lem}

\begin{proof}
  It suffices to check that the restriction of the canonical map to the covering $\widetilde{X}\to X$ as in
  \cite[Proposition~\ref*{prop:localdescription-of-subsections-OBla-reconstructionpaper}]{WiersigPeriods}
  is an isomorphism.
  Using \emph{loc. cit.}
  and~\cite[Lemma 3.33]{Bo21}, we compute
  \begin{align*}
    \I\left(\OB_{\dR}^{\dag,q,+}\right)|_{\widetilde{X}}
    &=\I\left(\OB_{\dR}^{\dag,q,+}|_{\widetilde{X}}\right) \\
    &\cong
    \I\left(
      \BB_{\dR}^{\dag,+}|_{\widetilde{X}}\widehat{\otimes}_{k_{0}}k_{0}\left\<\frac{Z_{1},\dots,Z_{d}}{p^{q}}\right\>
    \right) \\
    &\cong
    \I\left(\BB_{\dR}^{\dag,+}|_{\widetilde{X}}\right)
    \widehat{\otimes}_{\I\left(k_{0}\right)}\I\left(k_{0}\left\<\frac{Z_{1},\dots,Z_{d}}{p^{q}}\right\>\right)
      \\
    &\cong
    \I\left(\BB_{\dR}^{\dag,+}|_{\widetilde{X}}\right)
    \widehat{\otimes}_{\I\left(\BB_{\dR}^{q,+}\right)|_{\widetilde{X}}}
    \I\left(\BB_{\dR}^{q,+}\right)|_{\widetilde{X}}
    \widehat{\otimes}_{\I\left(k_{0}\right)}\I\left(k_{0}\left\<\frac{Z_{1},\dots,Z_{d}}{p^{q}}\right\>\right)\\
    &\cong
    \I\left(\BB_{\dR}^{\dag,+}\right)|_{\widetilde{X}}
    \widehat{\otimes}_{\I\left(\BB_{\dR}^{q,+}\right)}|_{\widetilde{X}}
    \I\left(\BB_{\dR}^{q,+}|_{\widetilde{X}}\right)\widehat{\otimes}_{\I\left(k_{0}\right)}
    \I\left(k_{0}\left\<\frac{Z_{1},\dots,Z_{d}}{p^{q}}\right\>\right) \\
    &\cong
    \I\left(\BB_{\dR}^{\dag,+}\right)|_{\widetilde{X}}
    \widehat{\otimes}_{\I\left(\BB_{\dR}^{q,+}\right)}|_{\widetilde{X}}
    \I\left(\BB_{\dR}^{q,+}|_{\widetilde{X}}\widehat{\otimes}_{k_{0}}
    k_{0}\left\<\frac{Z_{1},\dots,Z_{d}}{p^{q}}\right\>\right) \\
    &\cong
    \I\left(\BB_{\dR}^{\dag,+}\right)|_{\widetilde{X}}
    \widehat{\otimes}_{\I\left(\BB_{\dR}^{q,+}\right)}|_{\widetilde{X}}
    \I\left(\OB_{\dR}^{q,+}\right)|_{\widetilde{X}}
  \end{align*}
  in degree zero. It therefore remains to check that
  $\I\left(\OB_{\dR}^{q,+}\right)|_{\widetilde{X}}$ is flat as an $\I\left(\BB_{\dR}^{q,+}\right)|_{\widetilde{X}}$-module object.
  Using again~\cite[Proposition~\ref*{prop:localdescription-of-subsections-OBla-reconstructionpaper}]{WiersigPeriods}
  and~\cite[Lemma 3.33]{Bo21}, we compute
  \begin{align*}
    \I\left(\OB_{\dR}^{q,+}\right)|_{\widetilde{X}}
    &=\I\left(\OB_{\dR}^{q,+}|_{\widetilde{X}}\right) \\
    &\cong
    \I\left(
      \BB_{\dR}^{q,+}|_{\widetilde{X}}\widehat{\otimes}_{k_{0}}k_{0}\left\<\frac{Z_{1},\dots,Z_{d}}{p^{q}}\right\>
    \right) \\
    &\cong
    \I\left(\BB_{\dR}^{q,+}\right)|_{\widetilde{X}}
    \widehat{\otimes}_{\I\left(k_{0}\right)}\I\left(k_{0}\left\<\frac{Z_{1},\dots,Z_{d}}{p^{q}}\right\>\right).
  \end{align*}
  Thus Lemma~\ref{lem:constantsheaf-Tatealgebra-LH-exact} implies the flatness.
\end{proof}

\begin{proof}[Proof of Proposition~\ref{prop:construction-thm-finite-level}]
  Write $M_{q}^{\bullet}:=\I\left(\cal{D}_{q}(X)\right)\widehat{\otimes}_{\I\left(\Dcap(X)\right)}^{\rL}M^{\bullet}$.
  We have
  \begin{equation*}
    \I\left(\cal{D}_{q}^{\PD}(X)\right)\widehat{\otimes}_{\I\left(\Dcap(X)\right)}^{\rL}M^{\bullet}
    \isomap\I\left(\cal{D}_{q}^{\PD}(X)\right)\widehat{\otimes}_{\I\left(\cal{D}_{q}(X)\right)}^{\rL}M_{q}^{\bullet}
  \end{equation*}
  and
  \begin{equation}\label{eq:Soldagq-dependson-Mq--prop:construction-thm-finite-level}
  \begin{split}
    \nSol^{\dag,q}\left(M^{\bullet}\right)
    &\cong
      \R\shHom_{\I\left(\Dcap(X)\right)}\left(M_{X_{\proet}}^{\bullet},\I\left(\OB_{\dR}^{\dag,q,+}\right)\right), \\ 
    &\cong
          \R\shHom_{\I\left(\cal{D}_{q}(X)\right)}\left(
            \I\left(\cal{D}_{q}(X)\right)_{X_{\proet}}\widehat{\otimes}_{\Dcap(X)}^{\rL}M_{X_{\proet}}^{\bullet},
              \I\left(\OB_{\dR}^{\dag,q,+}\right)\right) \\
    &\stackrel{\text{\ref{lem:constantsheaf-commutes-compeltedtensorproduct-arbitrarybase-recpaper}}}{\cong}
          \R\shHom_{\I\left(\cal{D}_{q}(X)\right)}\left(M_{q,X_{\proet}}^{\bullet},\I\left(\OB_{\dR}^{\dag,q,+}\right)\right).
  \end{split}
  \end{equation}
  Henceforth,
  \begin{align*}
    \nRec\left(\nSol^{\dag,q}\left(M^{\bullet}\right)\right)
    \cong\R\Gamma\left(X,
      \R\shHom_{\I\left(\BB_{\dR}^{\dag,+}\right)}\left(
          \R\shHom_{\I\left(\cal{D}_{q}(X)\right)}\left(M_{q,X_{\proet}}^{\bullet},\I\left(\OB_{\dR}^{\dag,q,+}\right)\right),
          \I\left(\OB_{\pdR}^{\dag}\right)
        \right)
    \right).
  \end{align*}
  That is $\varrho_{M^{\bullet},q}^{\dag}$ only depends on $M_{q}^{\bullet}$.
  Since $M_{q}^{\bullet}$ is a bounded perfect complex of
  $\I\left(\cal{D}_{q}(X)\right)$-module objects, cf. Lemma~\ref{lem:sections-C-complexdefn-reconstructionpaper},
  we may assume $M_{q}^{\bullet}=\I\left(\cal{D}_{q}(X)\right)$.
  Then $\varrho_{M^{\bullet},q}^{\dag}$ becomes
  \begin{equation}\label{eq:construction-thm-finite-level-themap-Mbullet-isDqX-LH}
  \begin{split}
    \I\left(\cal{D}_{q}^{\PD}(X)\right)
    \to
    \R\Gamma\left(X,
      \R\shHom_{\I\left(\BB_{\dR}^{\dag,+}\right)}\left(
          \I\left(\OB_{\dR}^{\dag,q,+}\right),
          \I\left(\OB_{\pdR}^{\dag}\right)
        \right)
    \right) \\
    \stackrel{\text{\ref{lem:isoBdRdagBdRqOBdRdagq-LH--lem:prop:construction-thm-finite-level}}}{\cong}
      \R\Gamma\left(X,\R\shHom_{\I\left(\BB_{\dR}^{q,+}\right)}\left(
        \I\left(\OB_{\dR}^{q,+}\right),
        \I\left(\OB_{\pdR}^{\dag}\right)\right)\right).
  \end{split}
  \end{equation}  
  From Construction~\ref{construct:map--prop:construction-thm-finite-level}
  and Lemma~\ref{lem:prop:construction-thm-finite-level}, we find that the
  map~(\ref{eq:construction-thm-finite-level-themap-Mbullet-isDqX-LH}) coincides with
  $\varrho_{q}$ in Proposition~\ref{prop:derivedRq-iso-reconstructionpaper}.
  \emph{Loc. cit.} thus implies Proposition~\ref{prop:construction-thm-finite-level}.
\end{proof}

The following discussion makes frequent use of the notion of homotopy
limits in triangulated categories. We follow the terminology as
in~\cite[\href{https://stacks.math.columbia.edu/tag/08TB}{Tag 08TB}]{stacks-project}.

\begin{cor}\label{cor:construction-thm-finite-level-themap-fromDcapmoduleis-homlim}
  For every $\cal{C}$-complex $M^{\bullet}\in\D\left(\I\left(\Dcap(X)\right)\right)$,
  the maps
  \begin{equation*}
    \varrho_{M^{\bullet},q}^{\dag}\colon M^{\bullet}\to\nRec\left(\nSol^{\dag,q}\left(M^{\bullet}\right)\right)
  \end{equation*}
  exhibit $M^{\bullet}$ as the homotopy limit of the $\nRec\left(\nSol^{\dag,q}\left(M^{\bullet}\right)\right)$
  along the canonical maps
  \begin{equation*}
    \nRec\left(\nSol^{\dag,q+1}\left(M^{\bullet}\right)\right)
    \to\nRec\left(\nSol^{\dag,q}\left(M^{\bullet}\right)\right).
  \end{equation*}
  In symbols,
  $M^{\bullet}\cong\holim_{q}\nRec\left(\nSol^{\dag,q}\left(M^{\bullet}\right)\right)$.
\end{cor}

\begin{proof}
  Set $M_{q}^{\PD,\bullet}:=\I\left(\cal{D}_{q}^{\PD}(X)\right)\widehat{\otimes}_{\I\left(\Dcap(X)\right)}^{\rL}M^{\bullet}$.
  By Proposition~\ref{prop:construction-thm-finite-level}, we have to check
  that the canonical maps $M^{\bullet}\to M_{q}^{\PD,\bullet}$
  exhibit $M^{\bullet}$ as a homotopy limit of the $M_{q}^{\PD,\bullet}\in\D\left(\Dcap(X)\right)$,
  where the transition maps are the canonical ones. In symbols,
  $M^{\bullet}\cong\holim_{q}M_{q}^{\PD,\bullet}$.
  But we know $M^{\bullet}\cong\holim_{q}M_{q}^{\bullet}$
  for $M_{q}^{\bullet}:=\I\left(\cal{D}_{q}(X)\right)\widehat{\otimes}_{\I\left(\Dcap(X)\right)}^{\rL}M^{\bullet}$,
  cf.~\cite[Lemma 8.11]{Bo21}.
  Thus it remains to find a suitably functorial isomorphism
  $\holim_{q}M_{q}^{\PD,\bullet}\cong\holim_{q}M_{q}^{\bullet}$.
  We get one from the dual of~\cite[\href{https://stacks.math.columbia.edu/tag/0CRJ}{Tag 0CRJ}]{stacks-project},
  because both $q\mapsto M_{q}^{\PD,\bullet}$
  and $q\mapsto M_{q}^{\bullet}$ are cofinal in the diagram
  \begin{equation*}
    M_{q_{0}}^{\bullet} \to M_{q_{0}}^{\PD,\bullet} \to
    M_{q_{0}+1}^{\bullet} \to M_{q_{0}+1}^{\PD,\bullet} \to \dots
  \end{equation*}
  for a fixed large enough $q_{0}\gg0$.
  Here, all morphisms are the canonical ones.
\end{proof}

\begin{cor}\label{cor:construction-thm-finite-level-themap-fromDcapmoduleis-homlim-Hoi}
  For every $\cal{C}$-complex $M^{\bullet}\in\D\left(\I\left(\Dcap(X)\right)\right)$,
  the maps
  \begin{equation*}
    \varrho_{M^{\bullet},q}^{\dag}\colon M^{\bullet}\to\nRec\left(\nSol^{\dag,q}\left(M^{\bullet}\right)\right)
  \end{equation*}
  induce the following isomorphisms for all $i\in\ZZ$:
  \begin{equation*}
    \Ho^{i}\left(M^{\bullet}\right)\isomap
    \varprojlim_{q}\Ho^{i}\left(\nRec\left(\nSol^{\dag,q}\left(M^{\bullet}\right)\right)\right).
  \end{equation*}
\end{cor}

\begin{proof}
  By Corollary~\ref{cor:construction-thm-finite-level-themap-fromDcapmoduleis-homlim}
  and the short exact sequence in~\cite[Lemma 3.3]{Bo21}, we have to check that
  \begin{equation*}
    \R^{1}\varprojlim_{q}\Ho^{i}\left(\nRec\left(\nSol^{\dag,q}\left(M^{\bullet}\right)\right)\right)
    \stackrel{\text{\ref{prop:construction-thm-finite-level}}}{\cong}
    \R^{1}\varprojlim_{q}\Ho^{i}\left(M_{q}^{\PD,\bullet}\right)
  \end{equation*}
  vanishes for every $i\in\ZZ$ where
  $M_{q}^{\PD,\bullet}:=\I\left(\cal{D}_{q}^{\PD}(X)\right)\widehat{\otimes}_{\I\left(\Dcap(X)\right)}^{\rL}M^{\bullet}$.
  Using similar cofinality arguments as in the proof of
  Corollary~\ref{cor:construction-thm-finite-level-themap-fromDcapmoduleis-homlim},
  this follows from
  \begin{equation*}
    0=\R^{1}\varprojlim_{q}\Ho^{i}\left(M_{q}^{\bullet}\right)\isomap\R^{1}\varprojlim_{q}\Ho^{i}\left(M_{q}^{\PD,\bullet}\right)
  \end{equation*}
  where $M_{q}^{\bullet}:=\I\left(\cal{D}_{q}(X)\right)\widehat{\otimes}_{\I\left(\Dcap(X)\right)}^{\rL}M^{\bullet}$,
  cf. the~\cite[proof of Lemma 8.11]{Bo21}.
\end{proof}

Now we relate $\nRec\left(\nSol\left(M^{\bullet}\right)\right)$ and
$\holim_{q}\nRec\left(\nSol^{\dag,q}\left(M^{\bullet}\right)\right)$.
With Corollary~\ref{cor:construction-thm-finite-level-themap-fromDcapmoduleis-homlim},
this will allow us to relate $M^{\bullet}$ to $\nRec\left(\nSol\left(M^{\bullet}\right)\right)$
and prove Theorem~\ref{thm:recthm-globalsections-recpaper}.

Homotopy colimits are dual to homotopy limits.
We follow the formalism as
in~\cite[\href{https://stacks.math.columbia.edu/tag/0A5K}{Tag 0A5K}]{stacks-project}.

\begin{prop}\label{prop:relateSolq-andSol-recpaper}
  For every $\cal{C}$-complex $M^{\bullet}\in\D\left(\I\left(\Dcap(X)\right)\right)$
  and $q\gg0$, there exist maps
  \begin{equation*}
    \omega_{M^{\bullet}}^{\dag,q}\colon\nSol^{\dag,q}\left(M^{\bullet}\right)\to\nSol\left(M^{\bullet}\right)
  \end{equation*}
  of complexes of $\I\left(\BB_{\dR}^{\dag,+}\right)$-module objects satisfying
  the following two conditions:
  \begin{itemize}
    \item[(i)] 
    If $\iota_{M^{\bullet}}^{\dag,q}$ denotes the canonical maps, then following diagrams commute for all large $q$:
    \begin{equation}\label{cd1:prop:relateSolq-andSol-recpaper}
      \begin{tikzcd}
        \empty & \nSol\left(M^{\bullet}\right) \\
        \nSol^{\dag,q}\left(M^{\bullet}\right) \arrow{r}{\iota_{M^{\bullet}}^{\dag,q}}\arrow{ru}{\omega_{M^{\bullet}}^{\dag,q}} &
        \nSol^{\dag,q+1}\left(M^{\bullet}\right). \arrow[swap]{u}{\omega_{M^{\bullet}}^{\dag,q+1}}
      \end{tikzcd}
    \end{equation}  
    Furthermore, any induced morphism
    $\holim_{q}\nSol^{\dag,q}\left(M^{\bullet}\right)\isomap\nSol\left(M^{\bullet}\right)$
    is an isomorphism.
    \item[(ii)] The following
    diagrams commute for large $q$:
    \begin{equation}\label{cd2:prop:relateSolq-andSol-recpaper}
      \begin{tikzcd}
        M^{\bullet}\arrow{r}{\varrho_{M^{\bullet}}}\arrow[swap]{rd}{\varrho_{M^{\bullet},q}^{\dag}} & 
        \nRec\left(\nSol\left(M^{\bullet}\right)\right) \arrow{d}{\nRec\left(\omega_{M^{\bullet}}^{\dag,q}\right)} \\
        \empty & \nRec\left(\nSol^{\dag,q}\left(M^{\bullet}\right)\right).
      \end{tikzcd}
    \end{equation}  
    \end{itemize}
\end{prop}

We present the proof of Proposition~\ref{prop:relateSolq-andSol-recpaper}
on page~\pageref{proof:prop:relateSolq-andSol-recpaper} below.

\begin{lem}\label{lem:sigmapullback-as-colimit-recpaper}
  Let $\tau_{q}\colon X_{\proet}\times\NN_{\gg0}^{\op}\to X_{\proet}$
  denote the projection $(U,q)\mapsfrom U$. Then, for every complex $\cal{F}^{\bullet}$ of
  $\I\left(\BB_{\dR}^{*,+}\right)$-module objects, we have the functorial
  isomorphism
  \begin{equation*}
    \rL\sigma^{*}\cal{F}^{\bullet}\cong\varinjlim_{q\gg0}\tau_{q*}\cal{F}^{\bullet}.
  \end{equation*}
\end{lem}

\begin{proof}
  We observed in the proof of Corollary~\ref{cor:padic-ishimuraprosmansschneiders}
  that $\rL\sigma^{*}\cal{F}^{\bullet}\cong\sigma^{-1}\cal{F}^{\bullet}$, functorially.
  Since $\sigma^{-1}$ is exact, we may assume that $\cal{F}^{\bullet}=\cal{F}$
  is concentrated in degree zero. Now the computations
  \begin{equation*}
    \varinjlim_{\sigma^{-1}(V)\leftarrow U}\cal{F}(V)
    =\varinjlim_{q\gg0}\cal{F}\left(\left(U,q\right)\right)
    =\varinjlim_{q\gg0}\tau_{q*}\cal{F}\left(U\right)
  \end{equation*}
  for all $U\in X_{\proet}$
  induces the isomorphism $\sigma^{-1}\cal{F}\cong\varinjlim_{q}\tau_{q*}\cal{F}$.
\end{proof}

We continue to refer to the projections $X_{\proet}\times\NN_{\gg0}^{\op}\to X_{\proet}$,
$(U,q)\mapsfrom U$ by $\tau_{q}$.

\begin{lem}\label{lem:tauqpush-preserves-flasqueres-recpaper}
  For a flasque resolution $\I\left(\OB_{\dR}^{*,+}\right)\to\cal{I}^{\bullet}$,
  $\I\left(\OB_{\dR}^{q,+}\right)\to\tau_{q*}\cal{I}^{\bullet}$ is a flasque resolution.
\end{lem}

\begin{proof}
  It suffices to check that $\tau_{q*}$ is exact. As it is left exact, we only have to prove
  that it is right exact. Since a cokernel in categories of sheaves is computed via
  the cokernel in the category of presheaves and sheafification, we only have to check that
  $\tau_{q*}$ commutes with sheafification. Sheafification is given by
  $L\circ L$ for a certain operation $L$ defined in~\cite[Definition 2.2.3]{Sch99},
  thus it remains to check that $\tau_{q*}$ commutes with $L$. This follows
  from the definition of $X_{\proet}\times\NN_{\gg0}^{\op}$.
\end{proof}

\begin{lem}\label{lem:--lem:shHomIDcapXP--preservesflasqueres-recpaper}
  For any K-projective resolution $P^{\bullet}\to M^{\bullet}\in\D\left(\I\left(\Dcap(X)\right)\right)$,
  the functor
  \begin{equation*}
    \shHom_{\I\left(\Dcap(X)\right)}\left(P^{\bullet},-\right)
  \end{equation*}
  preserves quasi-isomorphisms.
\end{lem}

\begin{proof}
  As $\IndBan_{\I(k)}$ is elementary,
  we can fix a small strictly generating set $\cal{G}\subseteq\IndBan_{\I(k)}$,
  such that all $G\in\cal{G}$ are projective.
  By~\cite[Proposition 2.1.8]{Sch99}, we have to show~that
  \begin{equation*}
    \Hom_{\I(k)}\left( G , \shHom_{\I\left(\Dcap(X)\right)}\left(P^{\bullet},-\right) \right)
  \end{equation*}
  preserves quasi-isomorphisms for all $G\in\cal{G}$.
  But using~\cite[Lemma 3.11(ii)]{Bo21}, we find
  \begin{equation*}
  \begin{split}
    \Hom_{\I(k)}\left( G , \shHom_{\I\left(\Dcap(X)\right)}\left(P^{\bullet},-\right) \right)
    &\cong\Hom_{\I\left(\Dcap(X)\right)}\left( G \widehat{\otimes}_{\I(k)} P^{\bullet} , - \right) \\
    &\cong\Hom_{\I\left(\Dcap(X)\right)}\left(
      P^{\bullet} , 
      \intHom_{\I(k)}\left( G , - \right)
      \right).
  \end{split}
  \end{equation*}
  This functor preserves quasi-isomorphisms
  because $P^{\bullet}$ is $K$-projective and $G$ is projective.
\end{proof}

\begin{lem}\label{lem:shHomIDcapXP--preservesflasqueres-recpaper}
  For any K-projective resolution $P^{\bullet}\to M^{\bullet}\in\D\left(\I\left(\Dcap(X)\right)\right)$,
  the functor
  \begin{equation*}
    \shHom_{\I\left(\Dcap(X)\right)}\left(P_{X_{\proet}}^{\bullet},-\right)
  \end{equation*}
  preserves quasi-isomorphisms between complexes of flasque
  sheaves.
\end{lem}

\begin{proof}
  Given an acyclic complex $\cal{F}^{\bullet}$ of flasque sheaves
  of $\I\left(\Dcap(X)\right)$-module objects, we have to check that
  $C^{\bullet}:=\shHom_{\I\left(\Dcap(X)\right)}\left(P_{X_{\proet}}^{\bullet},\cal{F}^{\bullet}\right)$
  is acyclic. But for every $U\in X_{\proet}$, we find
  \begin{equation*}
    C^{\bullet}(U)=\intHom_{\I\left(\Dcap(X)\right)}\left(P^{\bullet},\cal{F}^{\bullet}(U)\right).
  \end{equation*}
  This complex is exact
  by Lemma~\ref{lem:--lem:shHomIDcapXP--preservesflasqueres-recpaper}, and
  because $\cal{F}^{\bullet}(U)$ is acyclic.
\end{proof}

\begin{lem}\label{lem:resolutionforSolq-recpaper}
  In the category of
  $\I\left(\Dcap(X)\right)_{X_{\proet}\times\NN_{\gg0}^{\op}}$-$\I\left(\BB_{\dR}^{*,+}\right)$-bimodule objects,
  we consider an injective resolution $\I\left(\OB_{\dR}^{*,+}\right)\to\cal{I}^{\bullet}$.
  Let $P^{\bullet}\to M^{\bullet}\in\D\left(\I\left(\Dcap(X)\right)\right)$ denote a K-projective resolution.
  Then, for all $q\in\NN_{\geq2}$,
  we have the following isomorphism in $\D\left(\BB_{\dR}^{q,+}\right)$:
  \begin{equation*}
    \nSol^{q}\left(M^{\bullet}\right)
    \cong
    \shHom_{\I\left(\Dcap(X)\right)}\left(P_{X_{\proet}}^{\bullet},\tau_{q*}\cal{I}^{\bullet}\right).
  \end{equation*}
  These isomorphisms are functorial in $q$.
\end{lem}

\begin{proof}
  This follows from Lemma~\ref{lem:tauqpush-preserves-flasqueres-recpaper}
  and~\ref{lem:shHomIDcapXP--preservesflasqueres-recpaper}.
\end{proof}

We introduce more solution functors.

\begin{defn}\label{defn:solfunctorq-globalsections-recpaper}
  For every $q\gg0$,
  denote the following functor by $\nSol^{q}$:
  \begin{equation*}
    \D\left(\I\left(\Dcap(X)\right)\right)^{\op}
    \to
    \D\left(\I\left(\BB_{\dR}^{\dag,+}\right)\right),
    M^{\bullet}
    \mapsto
    \R\shHom_{\I\left(\Dcap(X)\right)}\left(
      M^{\bullet}_{X_{\proet}},\I\left(\OB_{\dR}^{q,+}\right)
    \right).
  \end{equation*}
\end{defn}

\begin{lem}\label{lem:compareSolq-and-Soldagq-recpaper}
  Given a complex $M^{\bullet}$ of $\I\left(\Dcap(X)\right)$-module objects,
  the canonical morphism
  \begin{equation}\label{eq:themap--lem:compareSolq-and-Soldagq-recpaper}
  \begin{split}
    \nSol^{q}\left(M^{\bullet}\right)
      \widehat{\otimes}_{\I\left(\BB_{\dR}^{q,+}\right)}^{\rL}\I\left(\BB_{\dR}^{\dag,+}\right)
    \to\nSol^{\dag,q}\left(M^{\bullet}\right).
  \end{split}
  \end{equation}
  is an isomorphism when $M^{\bullet}$ is a $\cal{C}$-complex.
\end{lem}

\begin{proof}
  Set $M_{q}^{\bullet}:=\I\left(\cal{D}_{q}(X)\right)\widehat{\otimes}_{\I\left(\Dcap(X)\right)}^{\rL}M^{\bullet}$.
  Arguing as in~(\ref{eq:Soldagq-dependson-Mq--prop:construction-thm-finite-level}),
  we compute
  \begin{equation}\label{eq:Solq-dependson-Mq--lem:compareSolq-and-Soldagq-recpaper}
  \begin{split}
    \nSol^{q}\left(M^{\bullet}\right)
    \cong\R\shHom_{\I\left(\cal{D}_{q}(X)\right)}\left(M_{q,X_{\proet}}^{\bullet},\I\left(\OB_{\dR}^{q,+}\right)\right).
  \end{split}
  \end{equation}
  Now use~(\ref{eq:Soldagq-dependson-Mq--prop:construction-thm-finite-level})
  and~(\ref{eq:Solq-dependson-Mq--lem:compareSolq-and-Soldagq-recpaper})
  to find that~(\ref{eq:themap--lem:compareSolq-and-Soldagq-recpaper}) coincides with
  the canonical morphism
  \begin{equation*}
  \begin{split}
    \R\shHom_{\I\left(\cal{D}_{q}(X)\right)}\left(M_{q,X_{\proet}}^{\bullet},\I\left(\OB_{\dR}^{q,+}\right)\right)
      \widehat{\otimes}_{\I\left(\BB_{\dR}^{q,+}\right)}^{\rL}\I\left(\BB_{\dR}^{\dag,+}\right)
    \to
    \R\shHom_{\I\left(\cal{D}_{q}(X))\right)}\left(M_{q,X_{\proet}}^{\bullet},\I\left(\OB_{\dR}^{\dag,q,+}\right)\right).
  \end{split}
  \end{equation*}
  Since $M_{q}^{\bullet}$ is a bounded perfect complex,
  cf. Lemma~\ref{lem:sections-C-complexdefn-reconstructionpaper},
  it suffices to check that the above morphism
  is an isomorphism for $M_{q}^{\bullet}=\I\left(\cal{D}_{q}(X)\right)$. This follows from
  Lemma~\ref{lem:isoBdRdagBdRqOBdRdagq-LH--lem:prop:construction-thm-finite-level}.
\end{proof}

\begin{lem}\label{lem:hocolimqIBdRqplus-cong-IBdRdagplus}
  The canonical morphisms $\I\left(\BB_{\dR}^{q,+}\right)\to\I\left(\BB_{\dR}^{\dag,+}\right)$
  exhibit $\I\left(\BB_{\dR}^{\dag,+}\right)$ as the homotopy colimit of the
  $\I\left(\BB_{\dR}^{q,+}\right)$ along the canonical maps
  $\I\left(\BB_{\dR}^{q,+}\right)\to\I\left(\BB_{\dR}^{q+1,+}\right)$, in symbols
  \begin{equation*}
    \hocolim_{q}\I\left(\BB_{\dR}^{q,+}\right)\cong\I\left(\BB_{\dR}^{\dag,+}\right).
  \end{equation*}
\end{lem}

\begin{proof}
  $\I$ commutes with colimits, cf.~\cite[Proposition 2.1.16]{Sch99}, thus
  $\BB_{\dR}^{\dag,+}=\varinjlim_{q}\BB_{\dR}^{q,+}$ implies
  \begin{equation*}
    \varinjlim_{q}\I\left(\BB_{\dR}^{q,+}\right)\isomap\I\left(\BB_{\dR}^{\dag,+}\right).
  \end{equation*}
  Now apply~\cite[\href{https://stacks.math.columbia.edu/tag/0949}{Tag 0949}]{stacks-project}.
\end{proof}

\begin{lem}\label{lem:functorpreserveshocolim-recpaper}
  A triangulated functor preserves homotopy colimits if it preserves coproducts.
\end{lem}

\begin{proof}
  $E$ is a homotopy colimit of a diagram
  $E_{0}\to E_{1}\to\dots$ if it fits into a distinguished triangle
  \begin{equation*}
    E\longrightarrow\coprod_{n}E_{n}\stackrel{\id-\text{shift}}{\longrightarrow}\coprod_{n}E_{n}\longrightarrow E[1].
  \end{equation*}
  Any triangulated functor that preserves coproducts preserves such triangles.
\end{proof}

\begin{proof}[Proof of Proposition~\ref{prop:relateSolq-andSol-recpaper}]
\label{proof:prop:relateSolq-andSol-recpaper}
    Let $\I\left(\OB_{\dR}^{*,+}\right)\to\cal{I}^{\bullet}$ denote an injective resolution
  of $\I\left(\OB_{\dR}^{*,+}\right)$ as an
  $\I\left(\Dcap(X)\right)_{X_{\proet}\times\NN_{\gg0}^{\op}}$-$\I\left(\BB_{\dR}^{*,+}\right)$-bimodule object.
  By~\cite[Theorem 14.4.3]{KashiwaraSchapira2006}, we can fix a
  K-projective resolution $P^{\bullet}\to M^{\bullet}$ of $M^{\bullet}$  
  as an $\I\left(\Dcap(X)\right)$-module object. Now compute
  \begin{equation}\label{eq:1--prop:relateSolq-andSol-recpaper}
  \begin{split}
    \nSol\left(M^{\bullet}\right)
    &\stackrel{}{\cong}
    \rL\sigma^{*}\shHom_{\I\left(\Dcap(X)\right)}\left(P_{X_{\proet}\times\NN_{\gg0}^{\op}}^{\bullet},\cal{I}^{\bullet}\right) \\
    &\stackrel{\text{\ref{lem:sigmapullback-as-colimit-recpaper}}}{\cong}
    \varinjlim_{q}\tau_{q*}\shHom_{\I\left(\Dcap(X)\right)}\left(P_{X_{\proet}\times\NN_{\gg0}^{\op}}^{\bullet},\cal{I}^{\bullet}\right) \\
    &\stackrel{}{\cong}
    \varinjlim_{q}\shHom_{\I\left(\Dcap(X)\right)}\left(P_{X_{\proet}}^{\bullet},\tau_{q*}\cal{I}^{\bullet}\right).
  \end{split}
  \end{equation}
  The canonical composition
  $\I\left(\OB_{\dR}^{q,+}\right)
    \to\tau_{q*}\I\left(\OB_{\dR}^{*,+}\right)\to\tau_{q*}\cal{I}^{\bullet}$
  induces
  \begin{equation*}
    \chi_{M^{\bullet}}^{q}\colon
    \nSol^{q}\left(M^{\bullet}\right)
    \to\nSol\left(M^{\bullet}\right).
  \end{equation*}
  These $\chi_{M^{\bullet}}^{q}$ are by construction functorial in $q$.
  Since the $-\widehat{\otimes}_{\I\left(\BB_{\dR}^{q,+}\right)}^{\rL}\I\left(\BB_{\dR}^{\dag,+}\right)$
  are left adjoint to the forgetful functors, the $\chi_{M^{\bullet}}^{q}$ lift functorially
  to morphisms of complexes of $\I\left(\BB_{\dR}^{\dag,+}\right)$-module objects
  \begin{equation*}
    \omega_{M^{\bullet}}^{q}\colon
    \nSol^{q}\left(M^{\bullet}\right)
      \widehat{\otimes}_{\I\left(\BB_{\dR}^{q,+}\right)}^{\rL}\I\left(\BB_{\dR}^{\dag,+}\right)
    \to\nSol\left(M^{\bullet}\right).
  \end{equation*}  
  The following compositions are the desired morphisms $\omega_{M^{\bullet}}^{\dag,q}$:
  \begin{equation*}
    \nSol^{\dag,q}\left(M^{\bullet}\right)
    \stackrel{\text{\ref{lem:compareSolq-and-Soldagq-recpaper}}}{\cong}
    \nSol^{q}\left(M^{\bullet}\right)
      \widehat{\otimes}_{\I\left(\BB_{\dR}^{q,+}\right)}^{\rL}\I\left(\BB_{\dR}^{\dag,+}\right)
    \stackrel{\omega_{M^{\bullet}}^{q}}{\longrightarrow}
    \nSol\left(M^{\bullet}\right).
  \end{equation*}
  The commutativity of the diagrams~(\ref{cd1:prop:relateSolq-andSol-recpaper})
  and~(\ref{cd1:prop:relateSolq-andSol-recpaper})
  follows directly from the functoriality of the construction.
  It remains to show the second part of Proposition~\ref{prop:relateSolq-andSol-recpaper}(i).
  That is, we have to check that the maps $\omega_{M^{\bullet}}^{\dag,q}$ exhibit $\nSol\left(M^{\bullet}\right)$
  as a homotopy colimit of the $\nSol^{\dag,q}\left(M^{\bullet}\right)$ along the canonical maps
  $\iota_{M^{\bullet}}^{\dag,q}$. But with~\cite[\href{https://stacks.math.columbia.edu/tag/0949}{Tag 0949}]{stacks-project},
  we find
  \begin{equation*}
    \nSol\left(M^{\bullet}\right)
    \stackrel{\text{(\ref{eq:1--prop:relateSolq-andSol-recpaper})}}{\cong}
    \varinjlim_{q}\shHom_{\I\left(\Dcap(X)\right)}\left(P_{X_{\proet}}^{\bullet},\tau_{q*}\cal{I}^{\bullet}\right)
    \stackrel{\text{\ref{lem:resolutionforSolq-recpaper}}}{\cong}
    \hocolim_{q}\nSol^{q}\left(M^{\bullet}\right).   
  \end{equation*}
  Therefore, everything follows from Lemma~\ref{lem:dagq-relateto-q-recpaper} below.
\end{proof}

\begin{lem}\label{lem:dagq-relateto-q-recpaper}
  Fix a $\cal{C}$-complex $M^{\bullet}\in\D\left(\I\left(\Dcap(X)\right)\right)$
  as well as homotopy colimits of the diagrams
  $q\mapsto\nSol^{q}\left(M^{\bullet}\right)$
  and
  $q\mapsto\nSol^{\dag,q}\left(M^{\bullet}\right)$,
  where the structural maps are the canonical ones.
  These homotopy colimits are isomorphic,
  and an isomorphism $\varphi$ can be chosen such that for all
  $l\in\NN$ large enough, the following diagram commutes:
  \begin{equation*}
  \begin{tikzcd}
    \nSol^{\dag,l}\left(M^{\bullet}\right)
      \arrow{r} &
    \hocolim_{q}\nSol^{\dag,q}\left(M^{\bullet}\right) \\
    \nSol^{l}\left(M^{\bullet}\right)
      \arrow{r}\arrow{u} &
    \hocolim_{q}\nSol^{q}\left(M^{\bullet}\right).
      \arrow{u}{\varphi}
      \arrow[swap]{u}{\cong}
  \end{tikzcd}
  \end{equation*}
  Here, all maps except $\varphi$ are the canonical ones.
\end{lem}

\begin{proof}
  We find the desired map $\varphi$ with~\cite[\href{https://stacks.math.columbia.edu/tag/0CRI}{Tag 0CRI}]{stacks-project}.
  Now compute
  \begin{align*}
    \hocolim_{q}\nSol^{\dag,q}\left(M^{\bullet}\right)
    &\stackrel{\text{\ref{lem:compareSolq-and-Soldagq-recpaper}}}{\cong}
      \hocolim_{q}
      \nSol^{q}\left(M^{\bullet}\right)
      \widehat{\otimes}_{\I\left(\BB_{\dR}^{q,+}\right)}^{\rL}\I\left(\BB_{\dR}^{\dag,+}\right) \\
    &\stackrel{\text{\ref{lem:hocolimqIBdRqplus-cong-IBdRdagplus}}}{\cong}
      \hocolim_{q}\hocolim_{l\geq q}
        \nSol^{q}\left(M^{\bullet}\right)
        \widehat{\otimes}_{\I\left(\BB_{\dR}^{q,+}\right)}^{\rL}\I\left(\BB_{\dR}^{l,+}\right).
  \end{align*}
  Here, the derived completed tensor product commutes with
  homotopy colimits by Lemma~\ref{lem:functorpreserveshocolim-recpaper}.
  \begin{equation*}
    \hocolim_{q}\nSol^{q}\left(M^{\bullet}\right)
    \cong\hocolim_{q}
        \nSol^{q}\left(M^{\bullet}\right)
        \widehat{\otimes}_{\I\left(\BB_{\dR}^{q,+}\right)}^{\rL}\I\left(\BB_{\dR}^{q,+}\right).
  \end{equation*}
  is clear. It follows that $\varphi$ may be written as
  \begin{equation*}
    \hocolim_{q}
        \nSol^{q}\left(M^{\bullet}\right)
        \widehat{\otimes}_{\I\left(\BB_{\dR}^{q,+}\right)}^{\rL}\I\left(\BB_{\dR}^{q,+}\right)
    \to\hocolim_{q}\hocolim_{l\geq q}
        \nSol^{q}\left(M^{\bullet}\right)
        \widehat{\otimes}_{\I\left(\BB_{\dR}^{q,+}\right)}^{\rL}\I\left(\BB_{\dR}^{l,+}\right).
  \end{equation*}
  To check that it is an isomorphism, we consider its cohomology: for all $i\in\ZZ$,
  \begin{align*}
    &\Ho^{i}\left(
      \hocolim_{q}
        \nSol^{q}\left(M^{\bullet}\right)
        \widehat{\otimes}_{\I\left(\BB_{\dR}^{q,+}\right)}^{\rL}\I\left(\BB_{\dR}^{q,+}\right)
    \right) \\
    &\cong
    \varinjlim_{q}\Ho^{i}\left(
        \nSol^{q}\left(M^{\bullet}\right)
        \widehat{\otimes}_{\I\left(\BB_{\dR}^{q,+}\right)}^{\rL}\I\left(\BB_{\dR}^{q,+}\right)
    \right) \\
    &\cong\varinjlim_{q}\varinjlim_{l\geq q}\Ho^{i}\left(
        \nSol^{q}\left(M^{\bullet}\right)
        \widehat{\otimes}_{\I\left(\BB_{\dR}^{q,+}\right)}^{\rL}\I\left(\BB_{\dR}^{l,+}\right)
    \right) \\    
    &\cong\Ho^{i}\left(\hocolim_{q}\hocolim_{l\geq q}
        \nSol^{q}\left(M^{\bullet}\right)
        \widehat{\otimes}_{\I\left(\BB_{\dR}^{q,+}\right)}^{\rL}\I\left(\BB_{\dR}^{l,+}\right)
      \right).
  \end{align*}
  as desired. Here, we used~\cite[\href{https://stacks.math.columbia.edu/tag/0CRK}{Tag 0CRK}]{stacks-project}.
\end{proof}

\begin{lem}\label{lem:nRec-hocolimtoholim-recpaper}
  $\nRec$ sends homotopy colimits to homotopy limits.
\end{lem}

\begin{proof}
  Arguing as in the proof of Lemma~\ref{lem:functorpreserveshocolim-recpaper},
  this follows because $\nRec$ is triangulated and sends coproducts to products.
\end{proof}

\begin{cor}\label{cor:HoiRecSoldagqMbullet-HoiSolM-recpaper}
  For every $\cal{C}$-complex $M^{\bullet}\in\D\left(\I\left(\Dcap(X)\right)\right)$, fix morphisms
  $\omega_{M^{\bullet}}^{\dag,q}\colon\nSol^{\dag,q}\left(M^{\bullet}\right)\to\nSol\left(M^{\bullet}\right)$
  as in Proposition~\ref{prop:relateSolq-andSol-recpaper}.
  For every $i\in\ZZ$, these induce the isomorphism
  \begin{equation*}
    \Ho^{i}\left(\nRec\left(\nSol\left(M^{\bullet}\right)\right)\right)
    \isomap
    \varprojlim_{q}\Ho^{i}\left(\nRec\left(\nSol^{\dag,q}\left(M^{\bullet}\right)\right)\right).
  \end{equation*}
\end{cor}

\begin{proof}
  As explained in the proof of Corollary~\ref{cor:construction-thm-finite-level-themap-fromDcapmoduleis-homlim-Hoi},
  we have the following isomorphisms for all $i\in\ZZ$:
  \begin{equation*}
    \R^{1}\varprojlim_{q}\Ho^{i}\left(\nRec\left(\nSol^{\dag,q}\left(M^{\bullet}\right)\right)\right)
    \cong0.
  \end{equation*}
  Together with~\cite[Lemma 3.3]{Bo21}, this gives the first isomorphism in the computation
  \begin{align*}
    \varprojlim_{q}\Ho^{i}\left(\nRec\left(\nSol^{\dag,q}\left(M^{\bullet}\right)\right)\right)
    &\cong
      \Ho^{i}\left(\holim_{q}\nRec\left(\nSol^{\dag,q}\left(M^{\bullet}\right)\right)\right) \\
    &\stackrel{\text{\ref{lem:nRec-hocolimtoholim-recpaper}}}{\cong}
      \Ho^{i}\left(\nRec\left(\hocolim_{q}\nSol^{\dag,q}\left(M^{\bullet}\right)\right)\right) \\
    &\stackrel{\text{\ref{prop:relateSolq-andSol-recpaper}(i)}}{\cong}
      \Ho^{i}\left(\nRec\left(\nSol\left(M^{\bullet}\right)\right)\right).
  \end{align*}
  This proves Corollary~\ref{cor:HoiRecSoldagqMbullet-HoiSolM-recpaper}.
\end{proof}

\begin{proof}[Proof of Theorem~\ref{thm:recthm-globalsections-recpaper}]
\label{proof:thm:recthm-globalsections-recpaper}
  Let $i\in\ZZ$ be arbitrary. By Proposition~\ref{prop:relateSolq-andSol-recpaper},
  we have commutative diagrams
    \begin{equation*}
      \begin{tikzcd}
        \Ho^{i}\left(M^{\bullet}\right)\arrow{r}{\varrho_{M^{\bullet}}}\arrow[swap]{rd}{\Ho^{i}\left(\varrho_{M^{\bullet},q}^{\dag}\right)} & 
        \Ho^{i}\left(\nRec\left(\nSol\left(M^{\bullet}\right)\right)\right) \arrow{d}{\Ho^{i}\left(\nRec\left(\omega_{M^{\bullet}}^{\dag,q}\right)\right)} \\
        \empty & \Ho^{i}\left(\nRec\left(\nSol^{\dag,q}\left(M^{\bullet}\right)\right)\right)
      \end{tikzcd}
    \end{equation*}
  for large $q$. They are suitably functorial by \emph{loc. cit.} such that
  they induce the commutative diagrams
    \begin{equation*}
      \begin{tikzcd}
        \Ho^{i}\left(M^{\bullet}\right)\arrow{r}{\Ho^{i}\left(\varrho_{M^{\bullet}}\right)}\arrow[swap]{rd}{\varprojlim_{q}\Ho^{i}\left(\varrho_{M^{\bullet},q}^{\dag}\right)} & 
        \Ho^{i}\left(\nRec\left(\nSol\left(M^{\bullet}\right)\right)\right) \arrow{d}{\varprojlim_{q}\Ho^{i}\left(\nRec\left(\omega_{M^{\bullet}}^{\dag,q}\right)\right)} \\
        \empty & \varprojlim_{q}\Ho^{i}\left(\nRec\left(\nSol^{\dag,q}\left(M^{\bullet}\right)\right)\right).
      \end{tikzcd}
    \end{equation*}
    These diagrams, together with Corollaries~\ref{cor:construction-thm-finite-level-themap-fromDcapmoduleis-homlim-Hoi}
    and~\ref{cor:HoiRecSoldagqMbullet-HoiSolM-recpaper},
    imply that $\Ho^{i}\left(\varrho\right)$ is an isomorphism.
    Since $i$ was arbitrary, $\varrho$ is an isomorphism.
\end{proof}

\begin{remark}
  Since homotopy limits are not functorial, as least not on the level
  of triangulated categories, Proposition~\ref{prop:relateSolq-andSol-recpaper}
  does not imply the existence of a commutative diagram
  \begin{equation*}
    \begin{tikzcd}
      M^{\bullet}\arrow{r}{\varrho_{M^{\bullet}}}\arrow[swap]{rd} & 
      \nRec\left(\nSol\left(M^{\bullet}\right)\right) \arrow{d} \\
      \empty & \holim_{q}\nRec\left(\nSol^{\dag,q}\left(M^{\bullet}\right)\right).
    \end{tikzcd}
  \end{equation*}
  This is why in the proof of Theorem~\ref{thm:recthm-globalsections-recpaper},
  we argue on the level of cohomology.
\end{remark}


\subsubsection{Proof of Theorem~\ref{thm:reconstruction-thm-solfunctor-reconstructionpaper}}\label{subsubsec:subsec:proofsofmainthms-reconstructionpaper-endgame}

Lemma~\ref{lem:dealingwith-F-toRnustarG-reconstructionpaper} below
allows us to reduce Theorem~\ref{thm:reconstruction-thm-solfunctor-reconstructionpaper}
to Theorem~\ref{thm:recthm-globalsections-recpaper}.

\begin{lem}\label{lem:dealingwith-F-toRnustarG-reconstructionpaper}
  Fix complexes $\cal{F}^{\bullet}$ and $\cal{G}^{\bullet}$ of sheaves of $\I\left(k\right)$-ind-Banach spaces
  on $X$ and $X_{\proet}$, respectively, together with a morphism
  $\varrho\colon\cal{F}^{\bullet}\to\R\nu_{*}\cal{G}^{\bullet}$.
  This $\varrho$ canonically induces morphisms
  \begin{equation*}
    \varrho_{U}\colon\R\Gamma\left(U,\cal{F}^{\bullet}\right)
      \to\R\Gamma\left(U,\cal{G}^{\bullet}\right).
  \end{equation*}  
  for all admissible open subsets $U\subseteq X$.
  $\varrho$ is an isomorphism if $\varrho_{U}$
  is an isomorphism for all $U$.
\end{lem}

\begin{proof}
  Consider the adjunct $\nu^{-1}\cal{F}^{\bullet}\to\cal{G}$ of $\varrho$
  and take global sections. $\varrho_{X}$ is then the composition
  \begin{equation*}
    \R\Gamma\left(X,\cal{F}^{\bullet}\right)
    \to\R\Gamma\left(X,\nu^{-1}\cal{F}^{\bullet}\right)
      \to\R\Gamma\left(X,\cal{G}^{\bullet}\right).
  \end{equation*}
  Similarly, one defines $\varrho_{U}$ for all $U$.
  Now suppose all $\varrho_{U}$ are isomorphisms.
  Since $\Sh\left(X,\IndBan_{\I\left(k\right)}\right)$ is Grothendieck abelian,
  cf.~\cite[Theorem 3.24]{Bo21},
  \cite[\href{https://stacks.math.columbia.edu/tag/079P}{Tag 079P}]{stacks-project}
  allows to pick K-injective resolutions
  $\cal{F}^{\bullet}\to\cal{I}^{\bullet}$
  and $\cal{G}^{\bullet}\to\cal{J}^{\bullet}$.
  We denote their underlying complexes of
  presheaves by $\cal{I}_{\psh}^{\bullet}$ and
  $\cal{J}_{\psh}^{\bullet}$. Now
  we use the exactness of sheafification~\cite[Lemma 3.15(iii)]{Bo21}
  to compute, for every $i\in\ZZ$,
  \begin{equation*}
    \Ho^{i}\left(\cal{F}^{\bullet}\right)
    =\Ho^{i}\left(\cal{I}^{\bullet}\right)
    =\Ho^{i}\left(\cal{I}_{\psh}^{\bullet}\right)^{\sh}
    \cong\left( U\mapsto \Ho^{i}\left( U , \cal{F}^{\bullet} \right) \right)^{\sh}.
  \end{equation*}
  As every $\varrho_{U}$ is an isomorphism, this sheaf is isomorphic to
  \begin{equation*}
    \left( U\mapsto \Ho^{i}\left( U , \cal{G}^{\bullet} \right) \right)^{\sh}
    \cong\Ho^{i}\left(\nu_{*}\cal{J}_{\psh}^{\bullet}\right)^{\sh}
    \cong\Ho^{i}\left(\nu_{*}\cal{J}^{\bullet}\right)
    =\Ho^{i}\left(\R\nu_{*}\cal{G}^{\bullet}\right).
  \end{equation*}
  That is, $\Ho^{i}\left(\cal{F}^{\bullet}\right)\cong\Ho^{i}\left(\R\nu_{*}\cal{G}^{\bullet}\right)$
  for every $i\in\ZZ$, as desired.
\end{proof}

Recall the Definition~\ref{defn:solfrefunctor-globalsections-recpaper}
of the functor $\nSol$.

\begin{lem}\label{lem:nSol-vs-Sol-recpaper}
  Assume that $X$ is affinoid and equipped with an étale morphism $X\to\TT^{d}$.
  Fix a $\cal{C}$-complex $\cal{M}^{\bullet}\in\D\left(\I\left(\Dcap\right)\right)$
  and set $M^{\bullet}:=\R\Gamma\left(X,\cal{M}^{\bullet}\right)$.
  Then we have a canonical isomorphism
  \begin{equation*}
     \Sol\left(\cal{M}^{\bullet}\right)
     \cong
     \nSol\left(M^{\bullet}\right).
  \end{equation*}
\end{lem}

\begin{proof}
  This follows from the computation
  \begin{align*}
    &\R\shHom_{\lambda^{-1}\I\left(\Dcap\right)}
    \left(
      \lambda^{-1}\cal{M}^{\bullet},
      \I\left(\OB_{\dR}^{*,+}\right)
    \right) \\
    &\cong\R\shHom_{\lambda^{-1}\I\left(\Dcap\right)}
    \left(
      \lambda^{-1}\left(\I\left(\Dcap\right)\widehat{\otimes}_{\I\left(\Dcap(X)\right)}^{\rL}M^{\bullet}\right),
      \I\left(\OB_{\dR}^{*,+}\right)
    \right) \\
    &\stackrel{\text{\ref{lem:derivedpullbackcommuteswithspecifictensorproduct-recpaper}}}{\cong}
    \R\shHom_{\lambda^{-1}\I\left(\Dcap\right)}
    \left(
      \lambda^{-1}\I\left(\Dcap\right)\widehat{\otimes}_{\I\left(\Dcap(X)\right)}^{\rL}M^{\bullet}_{X_{\proet}\times\NN_{\gg0}^{\op}},
      \I\left(\OB_{\dR}^{*,+}\right)
    \right) \\
    &\cong
    \R\shHom_{\I\left(\Dcap(X)\right)}
    \left(
      M^{\bullet}_{X_{\proet}\times\NN_{\gg0}^{\op}},
      \I\left(\OB_{\dR}^{*,+}\right)
    \right) \\
    &\stackrel{\text{\ref{lem:constantsheaf-commuteswith-inverseimage-recpaper}}}{\cong}
    \R\shHom_{\I\left(\Dcap(X)\right)}
    \left(
      \lambda^{-1}M^{\bullet}_{X},
      \I\left(\OB_{\dR}^{*,+}\right)
    \right);
  \end{align*}
  see~\cite[Theorem 8.12]{Bo21} for the first isomorphism.
\end{proof}

\begin{proof}[Proof of Theorem~\ref{thm:reconstruction-thm-solfunctor-reconstructionpaper}]
  We assume that $X$ is affinoid and equipped with an étale morphism $X\to\TT^{d}$.
  Recall the Definition~\ref{defn:rec-reconstructionpaper}
  of the reconstruction functor. We have to check that the morphism
  \begin{equation*}
    \rho_{\cal{M}^{\bullet}}\colon
    \cal{M}^{\bullet}
    \to
    \R\nu_{*}\R\shHom_{\I\left(\BB_{\dR}^{\dag,+}\right)}
      \left(\Sol\left(\cal{M}^{\bullet}\right),\I\left(\OB_{\pdR}^{\dag}\right)\right)
  \end{equation*}
  is an isomorphism. By Lemma~\ref{lem:dealingwith-F-toRnustarG-reconstructionpaper},
  we may check that the induced morphisms
  \begin{equation*}
    \R\Gamma\left(U,\cal{M}^{\bullet}\right)
    \to
    \R\Gamma\left(U,\R\shHom_{\I\left(\BB_{\dR}^{\dag,+}\right)}
      \left(\Sol\left(\cal{M}^{\bullet}\right),\I\left(\OB_{\pdR}^{\dag}\right)\right)\right)
  \end{equation*}
  are isomorphisms for all affinoid subdomains $U\subseteq X$.
  Since $X$ is suitably general, it suffices to consider the case $U=X$
  only. Writing $M^{\bullet}:=\R\Gamma\left(X,\cal{M}^{\bullet}\right)$
  and using the notation as in
  \S\ref{subsubsec:subsec:proofsofmainthms-reconstructionpaper-recgloablsection},
  we have to check that
  \begin{equation}\label{eq:thisISvarrhoMbullet-thm:reconstruction-thm-solfunctor-reconstructionpaper}
    M^{\bullet}\to
    \nRec\left(\Sol\left(\cal{M}^{\bullet}\right)\right)
    \stackrel{\text{\ref{lem:nSol-vs-Sol-recpaper}}}{\cong}
    \nRec\left(\nSol\left(M^{\bullet}\right)\right)  
  \end{equation}  
  is an isomorphism. Comparing the Construction~\ref{construct:rhocalMbullet-reconstructionpaper}
  of $\rho_{\cal{M}^{\bullet}}$
  with the definition of $\varrho_{M^{\bullet}}$
  in \S\ref{subsubsec:subsec:proofsofmainthms-reconstructionpaper-recgloablsection},
  we find that~(\ref{eq:thisISvarrhoMbullet-thm:reconstruction-thm-solfunctor-reconstructionpaper})
  is $\varrho_{M^{\bullet}}$.
  It is thus an isomorphism by Theorem~\ref{thm:recthm-globalsections-recpaper},
  which applies because $M^{\bullet}$ is a $\cal{C}$-complex,
  cf.~\cite[Lemma 8.13]{Bo21}.
\end{proof}


\subsection{Proof of the fully faithfulness}
\label{subsec:proofofsolfunctorfullyfaithfull-reconstructionpaper}

Now we prove Theorem~\ref{thm:pdRdagsolfunctor-fullyfaithful-reconstructionpaper} (Theorem~\ref{introthm:easymain}).

We have to check that the canonical map
\begin{equation}\label{eq:whattocheck--thm:pdRdagsolfunctor-fullyfaithful-reconstructionpaper}
  \Hom_{\I\left(\Dcap\right)}\left( \cal{N}^{\bullet} , \cal{M}^{\bullet} \right) \\
  \to\Hom_{\I\left(\BB_{\pdR}^{\dag}\right)}\left(
    \SolB_{\pdR}^{\dag}\left(\cal{M}^{\bullet}\right),
    \SolB_{\dR}^{\dag}\left(\cal{N}^{\bullet}\right)
  \right)
\end{equation}
is a bijection for every two $\cal{C}$-complexes $\cal{M}^{\bullet}$ and $\cal{N}^{\bullet}$.
To argue locally, we first study
a version of~(\ref{eq:whattocheck--thm:pdRdagsolfunctor-fullyfaithful-reconstructionpaper})
for the internal homomorphisms: the morphism
\begin{equation}\label{eq:themapLH--subsubsec:proofofsolfunctorfullyfaithfull-reconstructionpaper-constructingmorphisms}
  \R\shHom_{\I\left(\Dcap\right)}\left( \cal{N}^{\bullet} , \cal{M}^{\bullet} \right)
  \to\R\shHom_{\I\left(\BB_{\pdR}^{\dag}\right)}\left(
    \SolB_{\pdR}^{\dag}\left(\cal{M}^{\bullet}\right),
    \SolB_{\dR}^{\dag}\left(\cal{N}^{\bullet}\right)
  \right)
\end{equation}
in $\D\left( \Sh\left( X , \IndBan_{\I\left(k_{0}\right)} \right)\right)$.
In \S\ref{subsubsec:proofofsolfunctorfullyfaithfull-reconstructionpaper-local},
see in particular Theorem~\ref{thm:solfunctorfullyfaithfull-reconstructionpaper-localversion},
we then check that the map is bijective, at least when $X$ is affinoid and equipped with
suitable étale coordinates.
Then we deduce a global result in \S\ref{subsubsec:proofofsolfunctorfullyfaithfull-reconstructionpaper-global}, namely
Theorem~\ref{thm:solfunctorfullyfaithfull-reconstructionpaper-globalversioninternalhom}.
Finally, we deduce Theorem~\ref{thm:pdRdagsolfunctor-fullyfaithful-reconstructionpaper}
on page~\pageref{proof:--thm:pdRdagsolfunctor-fullyfaithful-reconstructionpaper}.



\subsubsection{Local computations}\label{subsubsec:proofofsolfunctorfullyfaithfull-reconstructionpaper-local}

Throughout \S\ref{subsubsec:proofofsolfunctorfullyfaithfull-reconstructionpaper-local},
we assume $X$ to be affinoid and equipped with an étale morphism $g\colon X\to\TT^{d}$.
We check that a slight modification
of~(\ref{eq:themapLH--subsubsec:proofofsolfunctorfullyfaithfull-reconstructionpaper-constructingmorphisms})
is an isomorphism. To do this, we use a variant of the solution functor $\nSol$ in
\S\ref{subsubsec:subsec:proofsofmainthms-reconstructionpaper-recgloablsection}:
\begin{equation*}
  \nSolB_{\pdR}^{\dag}\colon\D\left(\I\left(\Dcap(X)\right)\right)^{\op}
  \to
  \D\left(\I\left(\BB_{\pdR}^{\dag}\right)\right),
  M^{\bullet}
  \mapsto
  \rL\sigma^{*}\R\shHom_{\I\left(\Dcap(X)\right)}\left(
    M^{\bullet}_{X_{\proet}\times\NN_{\gg0}^{\op}},\I\left(\OB_{\pdR}^{*}\right)
  \right).
\end{equation*}

From now on, we fix two $\cal{C}$-complexes
$M^{\bullet},N^{\bullet}\in\D_{\mathcal{C}}\left( \I\left(\Dcap(X)\right) \right)$.
By~\cite[Theorem 8.12]{Bo21}, we also have the $\cal{C}$-complexes
$\cal{M}^{\bullet}:=\I\left(\Dcap\right)\widehat{\otimes}_{\I\left(\Dcap(X)\right)}^{\rL}M^{\bullet}$ and
$\cal{N}^{\bullet}:=\I\left(\Dcap\right)\widehat{\otimes}_{\I\left(\Dcap(X)\right)}^{\rL}N^{\bullet}$
on $X$.

\begin{lem}\label{lem:constructmap1-thm:solfunctorfullyfaithfull-reconstructionpaper-localversion}
  In $\D\left(\IndBan_{\I\left(k\right)}\right)$, we have the canonical isomorphism
  \begin{equation*}
    \R\Gamma\left(X,\R\shHom_{\I\left(\Dcap\right)}\left( \cal{N}^{\bullet} , \cal{M}^{\bullet} \right)\right)
    \cong
    \R\shHom_{\I\left(\Dcap(X)\right)}\left( N^{\bullet} , M^{\bullet} \right)
  \end{equation*}
\end{lem}

\begin{proof}
  By~\cite[Theorem 8.12]{Bo21},
  $\R\Gamma\left(X,\cal{M}^{\bullet}\right)\cong M^{\bullet}$ and 
  $\R\Gamma\left(X,\cal{N}^{\bullet}\right)\cong N^{\bullet}$. 
  Therefore, we may proceed similarly as in
  the~\cite[proof of Proposition 18.6.6]{KashiwaraSchapira2006}.
\end{proof}

\begin{lem}\label{lem:constructmap2-thm:solfunctorfullyfaithfull-reconstructionpaper-localversion}
  In $\D\left(\IndBan_{\I\left(k_{0}\right)}\right)$, we have the canonical isomorphism
  \begin{equation*}
  \begin{split}
    &\R\Gamma\left(X,
      \R\nu_{*}\shHom_{\I\left(\BB_{\pdR}^{\dag}\right)}\left(
    \SolB_{\pdR}^{\dag}\left(\cal{M}^{\bullet}\right),
    \SolB_{\dR}^{\dag}\left(\cal{N}^{\bullet}\right)
  \right)
    \right) \\
    &\cong
    \R\intHom_{\I\left(\BB_{\pdR}^{\dag}\right)}\Hom_{\I\left(\BB_{\pdR}^{\dag}\right)}\left(
    \SolB_{\pdR}^{\dag}\left(\cal{M}^{\bullet}\right),
    \SolB_{\dR}^{\dag}\left(\cal{N}^{\bullet}\right)
  \right).
  \end{split}
  \end{equation*}
\end{lem}

\begin{proof}
  Use Lemma~\ref{lem:nSol-vs-Sol-recpaper}
  and then proceed as in
  the~\cite[proof of Proposition 18.6.6]{KashiwaraSchapira2006}.
\end{proof}

Now we consider the
morphism~(\ref{eq:themapLH--subsubsec:proofofsolfunctorfullyfaithfull-reconstructionpaper-constructingmorphisms}).
Using Lemma~\ref{lem:dealingwith-F-toRnustarG-reconstructionpaper},
~\ref{lem:constructmap1-thm:solfunctorfullyfaithfull-reconstructionpaper-localversion},
and~\ref{lem:constructmap2-thm:solfunctorfullyfaithfull-reconstructionpaper-localversion},
we find
\begin{equation}\label{eq:themap--thm:solfunctorfullyfaithfull-reconstructionpaper-localversion}
\begin{split}
  &\R\intHom_{\I\left(\Dcap(X)\right)}\left( N^{\bullet} , M^{\bullet} \right) \\
  &\to
    \R\intHom_{\I\left(\BB_{\pdR}^{\dag}\right)}\Hom_{\I\left(\BB_{\pdR}^{\dag}\right)}\left(
    \SolB_{\pdR}^{\dag}\left(\cal{M}^{\bullet}\right),
    \SolB_{\dR}^{\dag}\left(\cal{N}^{\bullet}\right)
  \right).
\end{split}
\end{equation}
This morphism in $\D\left(\Sh\left(X,\IndBan_{\I\left(k_{0}\right)}\right)\right)$
is canonical and functorial in both $M^{\bullet}$ and $N^{\bullet}$.

\begin{thm}\label{thm:solfunctorfullyfaithfull-reconstructionpaper-localversion}
  (\ref{eq:themap--thm:solfunctorfullyfaithfull-reconstructionpaper-localversion})
  is an isomorphism
  for all $\cal{C}$-complexes $M^{\bullet},N^{\bullet}\in\D_{\mathcal{C}}\left( \I\left(\Dcap(X)\right) \right)$.
\end{thm}

The proof of Theorem~\ref{thm:solfunctorfullyfaithfull-reconstructionpaper-localversion}
is on page~\pageref{proof--thm:solfunctorfullyfaithfull-reconstructionpaper-localversion}
below. We continue with preliminary computations.

\begin{defn}
  An $\I\left(\BB_{\pdR}^{\dag}\right)$-module object $\cal{G}$
  is \emph{acyclic} if $\Ho^{i}\left(X,\cal{G}\right)=0$ for all $i>0$.
\end{defn}

\begin{lem}\label{lem:directsumlocacyclic-implies-locallyacyclic-recpaper}
  Arbitrary coproducts of acyclic
  $\I\left(\BB_{\pdR}^{\dag}\right)$-module objects are acyclic.
\end{lem}

\begin{proof}
  Let $X_{\proet,\qc}\subseteq X_{\proet}$ denote the full subcategory
  of quasi-compact objects. For example,
  $X$ itself and any affinoid perfectoid $U\in X_{\proet}$ are quasi-compact
  by~\cite[Proposition 3.12]{Sch13pAdicHodge}.
  By Proposition~\ref{prop:siteqc-reducesheafcondtofincov},
  we may view all $\I\left(\BB_{\pdR}^{\dag}\right)$-module objects
  as sheaves on $X_{\proet,\qc}$, which we do.
  Since arbitrary coproducts are filtered colimits of
  finite coproducts, it suffices to show that
  filtered colimits of acyclic
  $\I\left(\BB_{\pdR}^{\dag}\right)$-module objects are acyclic.
  Since we are working on $X_{\proet,\qc}$,
  a straightforward generalisation of
  \cite[\href{https://stacks.math.columbia.edu/tag/0737}{Tag 0737}]{stacks-project}
  to for $\IndBan_{\I\left(k_{0}\right)}$ applies.
\end{proof}

\begin{defn}
  For all $q\gg0$, we introduce the following sheaf of $\BB_{\pdR}^{\dag}$-ind-Banach algebras
  \begin{equation*}
    \OB_{\pdR}^{\dag,q}:=\OB_{\dR}^{q,+} \widehat{\otimes}_{\BB_{\dR}^{q,+}} \BB_{\pdR}^{\dag}
  \end{equation*}
  Using Corollary~\ref{cor:OBdRqplus-bimodulestructureover-calDrPDXXBdRqplus},
  we may as well view $\OB_{\pdR}^{\dag,q}$
  as a $\nu^{-1}\cal{D}_{q}^{\PD}(X)_{X}$-$\BB_{\pdR}^{\dag}$-bimodule
  object.  
\end{defn}

\begin{lem}\label{lem:describehom--lem:largederivedhom-preserves-hocolim-recpaper}
  Given an acyclic $\I\left(\BB_{\pdR}^{\dag}\right)$-module object $\cal{G}$,
  \begin{equation*}
    \R\intHom_{\I\left(\BB_{\pdR}^{\dag}\right)}
    \left(
      \I\left(\OB_{\pdR}^{\dag,q}\right),
      \cal{G}
    \right)
    \cong
    \intHom_{\I\left(k_{0}\right)}\left(\I\left(k_{0}\left\<\frac{Z_{1},\dots,Z_{d}}{p^{q}}\right\>\right),\cal{G}(X)\right).
  \end{equation*}
\end{lem}

\begin{proof}
  Arguing as in the proof of Lemma~\ref{lem1:lem:prop:derivedRq-iso-degree0-reconstructionpaper},
  we find
  \begin{equation*}
    \intHom_{\I\left(\BB_{\pdR}^{\dag}\right)}
    \left(
      \I\left(\OB_{\pdR}^{\dag,q}\right),
      \cal{G}
    \right)
    \cong
    \intHom_{\I\left(k_{0}\right)}\left(\I\left(k_{0}\left\<\frac{Z_{1},\dots,Z_{d}}{p^{q}}\right\>\right),\cal{G}(X)\right).
  \end{equation*}
  It remains to check that the functor preserves quasi-isomorphisms of
  complexes of acyclic $\I\left(\BB_{\pdR}^{\dag}\right)$-module objects.
  By looking at mapping cones, we find that it suffices to check the following:
  Given any acyclic complex $\cal{G}^{\bullet}$ of acyclic $\I\left(\BB_{\pdR}^{\dag}\right)$-module
  objects,
  \begin{equation*}
    \intHom_{\I\left(\BB_{\pdR}^{\dag}\right)}
    \left(
      \I\left(\OB_{\pdR}^{\dag,q}\right),
      \cal{G}^{\bullet}
    \right)
    \cong
    \intHom_{\I\left(k_{0}\right)}\left(\I\left(k_{0}\left\<\frac{Z_{1},\dots,Z_{d}}{p^{q}}\right\>\right),\cal{G}^{\bullet}(X)\right).
  \end{equation*}
  is again acyclic. This follows from
  Lemma~\ref{lem:Tate-algebra-internally-projective-LHIndBanF-reconstructionpaper}
  and because the complex $\cal{G}^{\bullet}(X)$ is exact.
\end{proof}

\begin{lem}\label{lem:largederivedhom-preserves-hocolim-recpaper}
  The following functor preserves homotopy colimits:
  \begin{equation*}
    \R\intHom_{\I\left(\BB_{\pdR}^{\dag}\right)}
    \left(
      \I\left(\OB_{\pdR}^{\dag,q}\right),
      -
    \right)
    \colon \D\left(\I\left(\BB_{\pdR}^{\dag}\right)\right)
    \to \D\left(\Sh\left( X_{\proet} , \IndBan_{\I(k)} \right)\right)
  \end{equation*}
\end{lem}

\begin{proof}
  By Lemma~\ref{lem:functorpreserveshocolim-recpaper}, it suffices to check
  that the functor commutes with coproducts. It follows from~\cite[Lemma 3.15]{Bo21}
  and the exactness of direct sums in $\IndBan_{\I\left(k_{0}\right)}$
  that coproducts are exact in the category of
  $\I\left(\BB_{\pdR}^{\dag}\right)$-module objects. By
  Lemma~\ref{lem:describehom--lem:largederivedhom-preserves-hocolim-recpaper},
  it remains to check that
  \begin{equation}\label{eq:functor--lem:largederivedhom-preserves-hocolim-recpaper}
    \intHom_{\I\left(k_{0}\right)}\left(\I\left(k_{0}\left\<\frac{Z_{1},\dots,Z_{d}}{p^{q}}\right\>\right),-\right)
  \end{equation}
  preserves coproducts of objects in $\IndBan_{\I\left(k_{0}\right)}$.
  But by~\cite[Proposition 2.1.17]{Sch99},
  \begin{equation*}
    \IndBan_{\I\left(k_{0}\right)}\simeq\Ind\left(\LH\left(\Ban_{k_{0}}\right)\right).
  \end{equation*}
  Since $\I\left(k_{0}\left\<\frac{Z_{1},\dots,Z_{d}}{p^{q}}\right\>\right)\in\LH\left(\Ban_{k_{0}}\right)$,
  (\ref{eq:functor--lem:largederivedhom-preserves-hocolim-recpaper})
  commutes with filtered colimits.
  Because coproducts are filtered colimits of finite coproducts,
  (\ref{eq:functor--lem:largederivedhom-preserves-hocolim-recpaper})
  comutes with coproducts.
  Lemma~\ref{lem:largederivedhom-preserves-hocolim-recpaper}
  follows.
\end{proof}

\begin{lem}\label{lem:hocolimOBpdRl-cong-OBdRdagl}
  The canonical morphisms $\I\left(\OB_{\pdR}^{\dag,q}\right)\to\I\left(\OB_{\pdR}^{\dag}\right)$
  exhibit $\I\left(\OB_{\pdR}^{\dag}\right)$ as the homotopy colimit of the
  $\I\left(\OB_{\pdR}^{\dag,q}\right)$ along the canonical maps
  $\I\left(\OB_{\pdR}^{\dag,q}\right)\to\I\left(\OB_{\pdR}^{\dag,q+1}\right)$, in symbols
  \begin{equation*}
    \hocolim_{q}\I\left(\OB_{\pdR}^{\dag,q}\right)\cong\I\left(\OB_{\pdR}^{\dag}\right).
  \end{equation*}
\end{lem}

\begin{proof}
  $\I$ commutes with colimits, cf.~\cite[Proposition 2.1.16]{Sch99}, thus
  \begin{equation*}
    \varinjlim_{q}\OB_{\pdR}^{\dag,q}
    =\varinjlim_{q} \OB_{\dR}^{q,+} \widehat{\otimes}_{\BB_{\dR}^{q,+}} \BB_{\pdR}^{\dag}
    \cong \OB_{\dR}^{\dag,+} \widehat{\otimes}_{\BB_{\dR}^{\dag,+}} \BB_{\pdR}^{\dag}
    \cong \OB_{\pdR}^{\dag}
  \end{equation*}
  implies
  \begin{equation*}
    \varinjlim_{q}\I\left(\OB_{\pdR}^{\dag,q}\right)\isomap\I\left(\OB_{\pdR}^{\dag}\right).
  \end{equation*}
  Now apply~\cite[\href{https://stacks.math.columbia.edu/tag/0949}{Tag 0949}]{stacks-project}.
\end{proof}

\begin{construction}\label{construction:map-for--lem:fullyfaithful-localcomputation1-recpaper}
   Consider the functor
   \begin{equation*}
     \R\intHom_{\I\left(\cal{D}_{q}^{\PD}(X)\right)}\left(
       - , \I\left( \OB_{\pdR}^{\dag,q} \right)
     \right)
     \colon
     \D\left( \I\left(\cal{D}_{q}^{\PD}(X)\right) \right)^{\op}
     \to\D\left( \I\left( \BB_{\pdR}^{\dag} \right) \right).
   \end{equation*}
   For every $M^{\bullet},N^{\bullet}\in\D\left(\I\left(\Dcap(X)\right)\right)$ with
   \begin{equation*}
   \begin{split}
     M_{q}^{\PD,\bullet}
     &:=\I\left(\cal{D}_{l}^{\PD}(X)\right)
       \widehat{\otimes}_{\I\left(\Dcap(X)\right)}^{\rL}
       M^{\bullet} \text{ and} \\
     N_{q}^{\PD,\bullet}
     &:=\I\left(\cal{D}_{l}^{\PD}(X)\right)
        \widehat{\otimes}_{\I\left(\Dcap(X)\right)}^{\rL}
        N^{\bullet},
    \end{split}
    \end{equation*}
    this functor induces a functorial map
   \begin{equation}\label{eq:map1--construction:map-for--lem:fullyfaithful-localcomputation1-recpaper}
   \begin{split}
     &\R\intHom_{\I\left(\cal{D}_{q}^{\PD}(X)\right)}\left(
       N_{q}^{\PD,\bullet},
       M_{q}^{\PD,\bullet}
     \right) \\
     &\to
     \R\shHom_{\I\left( \BB_{\pdR}^{\dag}\right)}
     \left(
       \R\intHom_{\I\left(\cal{D}_{q}^{\PD}(X)\right)}\left(
         M_{q}^{\PD,\bullet} , \I\left( \OB_{\pdR}^{\dag,q} \right)
         \right),
       \R\intHom_{\I\left(\cal{D}_{q}^{\PD}(X)\right)}\left(
         N_{q}^{\PD,\bullet} , \I\left( \OB_{\pdR}^{\dag,q} \right)
         \right)
     \right).
   \end{split}
   \end{equation}
   Now let $l\geq q$ be an integer. We have the canonical map
   \begin{equation*}
   \begin{split}
     &\R\intHom_{\I\left(\cal{D}_{q}^{\PD}(X)\right)}\left(
       N_{q}^{\PD,\bullet} , \I\left( \OB_{\pdR}^{\dag,q} \right)
       \right) \\
     &\cong
     \R\intHom_{\I\left(\Dcap(X)\right)}\left(
       N^{\bullet} , \I\left( \OB_{\pdR}^{\dag,q} \right)
       \right)
     \to
     \R\intHom_{\I\left(\Dcap(X)\right)}\left(
       N^{\bullet} , \I\left( \OB_{\pdR}^{\dag,l} \right)
       \right).
   \end{split}
   \end{equation*}
   Apply
   $\R\shHom_{\I\left( \BB_{\pdR}^{\dag}\right)}
     \left(
       \R\intHom_{\I\left(\cal{D}_{q}^{\PD}(X)\right)}\left(
         M_{q}^{\PD,\bullet} , \I\left( \OB_{\pdR}^{\dag,q} \right)
         \right),
         -     
     \right)$
   to this morphism to get
   \begin{equation}\label{eq:map2--construction:map-for--lem:fullyfaithful-localcomputation1-recpaper}
   \begin{split}
     &\R\shHom_{\I\left( \BB_{\pdR}^{\dag}\right)}
     \left(
       \R\intHom_{\I\left(\cal{D}_{q}^{\PD}(X)\right)}\left(
         M_{q}^{\PD,\bullet} , \I\left( \OB_{\pdR}^{\dag,q} \right)
         \right),
     \R\intHom_{\I\left(\cal{D}_{q}^{\PD}(X)\right)}\left(
       N_{q}^{\PD,\bullet} , \I\left( \OB_{\pdR}^{\dag,q} \right)
       \right)
     \right) \\
     &\to
     \R\shHom_{\I\left( \BB_{\pdR}^{\dag}\right)}
     \left(
       \R\intHom_{\I\left(\cal{D}_{q}^{\PD}(X)\right)}\left(
         M_{q}^{\PD,\bullet} , \I\left( \OB_{\pdR}^{\dag,q} \right)
         \right),
     \R\intHom_{\I\left(\Dcap(X)\right)}\left(
       N^{\bullet} , \I\left( \OB_{\pdR}^{\dag,l} \right)
       \right)
       \right).
     \end{split}
   \end{equation}   
   Now compose~(\ref{eq:map1--construction:map-for--lem:fullyfaithful-localcomputation1-recpaper})
   and~(\ref{eq:map2--construction:map-for--lem:fullyfaithful-localcomputation1-recpaper})
   to obtain the following morphism in $\D\left(\IndBan_{\I\left(k_{0}\right)}\right)$:
    \begin{equation}\label{eq:map3--construction:map-for--lem:fullyfaithful-localcomputation1-recpaper}
    \begin{split}
     &\R\intHom_{\I\left(\Dcap(X)\right)}\left(
       N^{\bullet},
       M_{q}^{\PD,\bullet}
     \right) \\
     &\cong\R\intHom_{\I\left(\cal{D}_{q}^{\PD}(X)\right)}\left(
       N_{q}^{\PD,\bullet},
       M_{q}^{\PD,\bullet}
     \right) \\
     &\to
     \R\shHom_{\I\left( \BB_{\pdR}^{\dag}\right)}
     \left(
       \R\intHom_{\I\left(\cal{D}_{q}^{\PD}(X)\right)}\left(
         M_{q}^{\PD,\bullet} , \I\left( \OB_{\pdR}^{\dag,q} \right)
         \right),
       \R\intHom_{\I\left(\Dcap(X)\right)}\left(
         N^{\bullet} , \I\left( \OB_{\pdR}^{\dag,l} \right)
         \right),     
     \right).
   \end{split}
   \end{equation}
   It is by construction canonical and
   functorial in $N^{\bullet}$, $M_{q}^{\PD,\bullet}$, and $l$.
   Therefore, we obtain
    \begin{equation}\label{eq:themap--construction:map-for--lem:fullyfaithful-localcomputation1-recpaper}
    \begin{split}
     &\R\intHom_{\I\left(\Dcap(X)\right)}\left(
       N^{\bullet},
       M_{q}^{\PD,\bullet}
     \right) \\
     &\to\hocolim_{l\geq q}
     \R\shHom_{\I\left( \BB_{\pdR}^{\dag}\right)}
     \left(
       \R\intHom_{\I\left(\cal{D}_{q}^{\PD}(X)\right)}\left(
         M_{q}^{\PD,\bullet} , \I\left( \OB_{\pdR}^{\dag,q} \right)
         \right),
       \R\intHom_{\I\left(\Dcap(X)\right)}\left(
         N^{\bullet} , \I\left( \OB_{\pdR}^{\dag,l} \right)
         \right)
     \right)
   \end{split}
   \end{equation}   
   by picking a homotopy colimit.
\end{construction}

\begin{lem}\label{lem:fullyfaithful-localcomputation1-recpaper}
  (\ref{eq:themap--construction:map-for--lem:fullyfaithful-localcomputation1-recpaper})
  is an isomorphism if
  $N^{\bullet}$ is a $\cal{C}$-complex
  and $M_{q}^{\PD,\bullet}=\I\left(\cal{D}_{q}^{\PD}(X)\right)$, that is
  \begin{equation*}
    \begin{split}
    &\R\intHom_{\I\left(\Dcap(X)\right)}\left( N^{\bullet} , \I\left(\cal{D}_{q}^{\PD}(X)\right) \right) \\
    &\isomap
      \hocolim_{l\geq q}
      \R\intHom_{\I\left(\BB_{\pdR}^{\dag}\right)}
        \left(
          \I\left(\OB_{\pdR}^{\dag,q}\right) , 
          \R\shHom_{\I\left(\Dcap(X)\right)}\left(
            N^{\bullet},
            \I\left(\OB_{\pdR}^{\dag,l}\right)
          \right)
        \right).
    \end{split}
  \end{equation*}
\end{lem}

\begin{proof}
  Write
  \begin{equation*}
    N_{q}^{\PD,\bullet}
    :=\I\left(\cal{D}_{l}^{\PD}(X)\right)
      \widehat{\otimes}_{\I\left(\Dcap(X)\right)}^{\rL}
      N^{\bullet}
  \end{equation*}
  as in Construction~\ref{construction:map-for--lem:fullyfaithful-localcomputation1-recpaper}.
  Now we find
  \begin{equation*}
    \begin{split}
      \R\intHom_{\I\left(\Dcap(X)\right)}\left( N^{\bullet} , \I\left(\cal{D}_{q}^{\PD}(X)\right) \right)
      =\R\intHom_{\I\left(\cal{D}_{l}^{\PD}(X)\right)}\left( N_{q}^{\PD,\bullet} , \I\left(\cal{D}_{q}^{\PD}(X)\right) \right)
    \end{split}
  \end{equation*}  
  and
  \begin{equation*}
    \begin{split}
      &\R\intHom_{\I\left(\BB_{\pdR}^{\dag}\right)}
        \left(
          \I\left(\OB_{\pdR}^{\dag,q}\right) , 
          \R\shHom_{\I\left(\Dcap(X)\right)}\left(
            N^{\bullet},
            \I\left(\OB_{\pdR}^{\dag,l}\right)
          \right)
        \right) \\
      &\cong\R\intHom_{\I\left(\BB_{\pdR}^{\dag}\right)}
        \left(
          \I\left(\OB_{\pdR}^{\dag,q}\right) , 
          \R\shHom_{\I\left(\cal{D}_{q}^{\PD}(X)\right)}\left(
            N_{q}^{\bullet},
            \I\left(\OB_{\pdR}^{\dag,l}\right)
          \right)
        \right)
    \end{split}
  \end{equation*}  
  By Lemma~\ref{lem:sections-C-complexdefn-reconstructionpaper}, $N_{q}^{\bullet}$ is a bounded perfect complex
  of $\I\left(\cal{D}_{q}^{\PD}(X)\right)$-module objects.
  Therefore, we may assume $N_{q}^{\bullet}=\I\left(\cal{D}_{q}^{\PD}(X)\right)$.
  In this case the morphism of interest becomes
  \begin{equation*}
    \I\left(\cal{D}_{q}^{\PD}(X) \right) \\
    \to
      \hocolim_{l}
      \R\intHom_{\I\left(\BB_{\pdR}^{\dag}\right)}
        \left(
          \I\left(\OB_{\pdR}^{\dag,q}\right),
            \I\left(\OB_{\pdR}^{\dag,l}\right)
        \right)
  \end{equation*}  
  It remains to check that this morphism is an isomorphism.
  This follows from
  \begin{equation*}
  \begin{split}
      &\hocolim_{l}
      \R\intHom_{\I\left(\BB_{\pdR}^{\dag}\right)}
        \left(
          \I\left(\OB_{\pdR}^{\dag,q}\right),
            \I\left(\OB_{\pdR}^{\dag,l}\right)
        \right) \\
      &\stackrel{\text{\ref{lem:largederivedhom-preserves-hocolim-recpaper}}}{\cong}
      \R\intHom_{\I\left(\BB_{\pdR}^{\dag}\right)}
        \left(
          \I\left(\OB_{\pdR}^{\dag,q}\right),
            \hocolim_{l}\I\left(\OB_{\pdR}^{\dag,l}\right)
        \right) \\
      &\stackrel{\text{\ref{lem:hocolimOBpdRl-cong-OBdRdagl}}}{\cong}
      \R\intHom_{\I\left(\BB_{\pdR}^{\dag}\right)}
        \left(
          \I\left(\OB_{\pdR}^{\dag,q}\right),
            \I\left(\OB_{\pdR}^{\dag}\right)
        \right) \\        
      &\stackrel{\text{\ref{prop:derivedRq-iso-reconstructionpaper}}}{\cong}
      \I\left(\cal{D}_{q}^{\PD}(X) \right),  
  \end{split}
  \end{equation*}
  as desired.
\end{proof}

\begin{lem}\label{lem:fullyfaithful-localcomputation2-recpaper}
  (\ref{eq:themap--construction:map-for--lem:fullyfaithful-localcomputation1-recpaper})
  is an isomorphism if
  $N^{\bullet}$ is a $\cal{C}$-complex
  and $M_{q}^{\PD,\bullet}$ is an arbitrary
  bounded perfect complex of $\I\left(\cal{D}_{q}^{\PD}(X)\right)$-module objects,
  that is
  \begin{equation*}
    \begin{split}
    &\R\intHom_{\I\left(\Dcap(X)\right)}\left( N^{\bullet} , M_{q}^{\bullet,\PD} \right) \\
    &\isomap
      \hocolim_{l\geq q}
     \R\shHom_{\I\left( \BB_{\pdR}^{\dag}\right)}
     \left(
       \R\intHom_{\I\left(\cal{D}_{q}^{\PD}(X)\right)}\left(
         M_{q}^{\PD,\bullet} , \I\left( \OB_{\pdR}^{\dag,q} \right)
         \right),
       \R\intHom_{\I\left(\Dcap(X)\right)}\left(
         N^{\bullet} , \I\left( \OB_{\pdR}^{\dag,l} \right)
         \right) 
     \right).
    \end{split}
  \end{equation*}
\end{lem}

\begin{proof}
  Since $M_{q}^{\bullet,\PD}$ is a bounded perfect complex, we may assume
  $M_{q}^{\bullet,\PD}=\I\left(\cal{D}_{q}^{\PD}(X)\right)$
  and apply Lemma~\ref{lem:fullyfaithful-localcomputation1-recpaper}.
\end{proof}

\begin{proof}[Proof of Theorem~\ref{thm:solfunctorfullyfaithfull-reconstructionpaper-localversion}]
\label{proof--thm:solfunctorfullyfaithfull-reconstructionpaper-localversion}
  Arguing as in the proof of Corollary~\ref{cor:construction-thm-finite-level-themap-fromDcapmoduleis-homlim},
  we find $M^{\bullet}\cong\holim_{q}M_{q}^{\bullet,\PD}$. This becomes useful in the following
  computation:
  \begin{equation*}
  \begin{split}
    &\R\intHom_{\I\left(\Dcap(X)\right)}\left(N^{\bullet},M^{\bullet}\right) \\
    &\cong\R\intHom_{\I\left(\Dcap(X)\right)}\left(N^{\bullet},\holim_{q}M_{q}^{\bullet,\PD}\right) \\
    &\cong\holim_{q}\R\intHom_{\I\left(\Dcap(X)\right)}\left(N^{\bullet},M_{q}^{\bullet,\PD}\right) \\
    &\stackrel{\text{\ref{lem:fullyfaithful-localcomputation2-recpaper}}}{\cong}\holim_{q}\hocolim_{l\geq q}
     \R\shHom_{\I\left( \BB_{\pdR}^{\dag}\right)}
     \left(
       \R\intHom_{\I\left(\cal{D}_{q}^{\PD}(X)\right)}\left(
         M_{q}^{\PD,\bullet} , \I\left( \OB_{\pdR}^{\dag,q} \right)
         \right),
       \R\intHom_{\I\left(\Dcap(X)\right)}\left(
         N^{\bullet} , \I\left( \OB_{\pdR}^{\dag,l} \right)
         \right) 
     \right) \\
    &\cong
     \R\shHom_{\I\left( \BB_{\pdR}^{\dag}\right)}
     \left(
       \hocolim_{q}\R\intHom_{\I\left(\cal{D}_{q}^{\PD}(X)\right)}\left(
         M_{q}^{\PD,\bullet} , \I\left( \OB_{\pdR}^{\dag,q} \right)
         \right),
       \hocolim_{l\geq q}\R\intHom_{\I\left(\Dcap(X)\right)}\left(
         N^{\bullet} , \I\left( \OB_{\pdR}^{\dag,l} \right)
         \right) 
     \right) \\
    &\cong
     \R\shHom_{\I\left( \BB_{\pdR}^{\dag}\right)}
     \left(
       \hocolim_{q}\R\intHom_{\I\left(\Dcap(X)\right)}\left(
         M^{\bullet} , \I\left( \OB_{\pdR}^{\dag,q} \right)
         \right),
       \hocolim_{l\geq q}\R\intHom_{\I\left(\Dcap(X)\right)}\left(
         N^{\bullet} , \I\left( \OB_{\pdR}^{\dag,l} \right)
         \right) 
     \right) \\
    &\cong
    \R\intHom_{\I\left(\BB_{\pdR}^{\dag}\right)}\left(
    \nSolB_{\pdR}^{\dag}\left(M^{\bullet}\right),
    \nSolB_{\pdR}^{\dag}\left(N^{\bullet}\right)\right).
  \end{split}
  \end{equation*}
  The last isomorphism follows from an appropriate version of Proposition~\ref{prop:relateSolq-andSol-recpaper}.
\end{proof}


\subsubsection{Global computations}\label{subsubsec:proofofsolfunctorfullyfaithfull-reconstructionpaper-global}

Here in \S\ref{subsubsec:proofofsolfunctorfullyfaithfull-reconstructionpaper-global},
we finally complete the proof of Theorem~\ref{subsec:proofofsolfunctorfullyfaithfull-reconstructionpaper}.
We therefore work in full generality, that is $X$ is an arbitrary smooth rigid-analytic $k$-variety.

\begin{thm}\label{thm:solfunctorfullyfaithfull-reconstructionpaper-globalversioninternalhom}
  For any two $\cal{C}$-complexes
  $\cal{M}^{\bullet},\cal{N}^{\bullet}\in\D_{\mathcal{C}}\left( \I\left(\Dcap\right) \right)$,
  the canonical morphism
  (\ref{eq:themapLH--subsubsec:proofofsolfunctorfullyfaithfull-reconstructionpaper-constructingmorphisms})
  on
  page~\pageref{eq:themapLH--subsubsec:proofofsolfunctorfullyfaithfull-reconstructionpaper-constructingmorphisms}
  is an isomorphism in $\D\left(\Sh\left(X,\IndBan_{\I\left(k_{0}\right)}\right)\right)$.
\end{thm}

\begin{proof}
  The statement is local in $X$, which is why we can safely assume that $X$ is affinoid and equipped with
  an étale morphism. By Lemma~\ref{lem:dealingwith-F-toRnustarG-reconstructionpaper},
  we have to check that the induced morphisms
  \begin{equation*}
  \begin{split}
  &\R\Gamma\left(U,\R\shHom_{\I\left(\Dcap\right)}\left( \cal{N}^{\bullet} , \cal{M}^{\bullet} \right)\right) \\
  &\to
    \R\Gamma\left(U,\R\shHom_{\I\left(\BB_{\pdR}^{\dag}\right)}\left(
      \SolB_{\pdR}^{\dag}\left(\cal{M}^{\bullet}\right),
      \SolB_{\pdR}^{\dag}\left(\cal{N}^{\bullet}\right),
    \right)
  \right)
  \end{split}
  \end{equation*}
  are isomorphisms for all admissible open subsets $U\subseteq X$. Since $X$
  is suitably general, it suffices to treat the case $U=X$ only.
  This is Theorem~\ref{thm:solfunctorfullyfaithfull-reconstructionpaper-localversion},
  by Lemma~\ref{lem:constructmap1-thm:solfunctorfullyfaithfull-reconstructionpaper-localversion}
  and~\ref{lem:constructmap2-thm:solfunctorfullyfaithfull-reconstructionpaper-localversion}.
\end{proof}

Lemma~\ref{lem:shHom-extHom--thm:solfunctorfullyfaithfull-reconstructionpaper}
and~\ref{lem:shHom-extHom2--thm:solfunctorfullyfaithfull-reconstructionpaper}
below allow us to deduce Theorem~\ref{thm:pdRdagsolfunctor-fullyfaithful-reconstructionpaper}
from Theorem~\ref{thm:solfunctorfullyfaithfull-reconstructionpaper-globalversioninternalhom}.

\begin{lem}\label{lem:shHom-extHom--thm:solfunctorfullyfaithfull-reconstructionpaper}
  Given $\cal{M}^{\bullet},\cal{N}^{\bullet}\in\D\left(\I\left(\Dcap\right)\right)$
  we have a functorial isomorphism
  \begin{equation*}
    \Hom_{\I\left(k_{0}\right)}\left(
      \I\left(k_{0}\right)_{X},
      \R\shHom_{\I\left(\Dcap\right)}\left( \cal{N}^{\bullet} , \cal{M}^{\bullet} \right)
    \right)
    \cong
    \Hom_{\I\left(\Dcap\right)}\left( \cal{N}^{\bullet} , \cal{M}^{\bullet} \right)
  \end{equation*}
  of abstract $k_{0}$-vector spaces. Here, $\I\left(k_{0}\right)_{X}$ denotes the constant sheaf.
\end{lem}

\begin{proof}
  We find
  \begin{equation*}
  \begin{split}
    \Hom_{\I\left(k_{0}\right)}\left( \I\left(k_{0}\right)_{X}  , \R\shHom_{\I\left(\Dcap\right)}\left( \cal{N}^{\bullet} , \cal{M}^{\bullet} \right)\right)
    &\cong\Hom_{\I\left(\Dcap\right)}\left( \I\left(k_{0}\right)_{X}\widehat{\otimes}_{\I(k)}^{\rL}\cal{N}^{\bullet}, \cal{M}^{\bullet} \right) \\
    &\cong\Hom_{\I\left(\Dcap\right)}\left( \cal{M}^{\bullet} , \cal{N}^{\bullet} \right)
  \end{split}
  \end{equation*}
  using~\cite[Lemma 3.11(ii)]{Bo21}.
\end{proof}

\begin{lem}\label{lem:shHom-extHom2--thm:solfunctorfullyfaithfull-reconstructionpaper}
  Given $\cal{F}^{\bullet},\cal{G}^{\bullet}\in\D\left(\I\left(\BB_{\pdR}^{\dag}\right)\right)$,
  we have a functorial isomorphism
  \begin{equation*}
    \Hom_{\I\left(k_{0}\right)}\left(
      \I\left(k_{0}\right)_{X},
      \R\nu_{*}\R\shHom_{\I\left(\BB_{\pdR}^{\dag}\right)}\left( \cal{F}^{\bullet} , \cal{G}^{\bullet} \right)
    \right)
    \cong
    \Hom_{\I\left(\BB_{\pdR}^{\dag}\right)}\left( \cal{F}^{\bullet} , \cal{G}^{\bullet} \right)
  \end{equation*}
  of abstract $k_{0}$-vector spaces. Here, $\I\left(k\right)_{X}$ denotes again the constant sheaf.
\end{lem}

\begin{proof}
  If $j\colon X\to *$ is the projection to a point, then $j^{-1}\I\left(k_{0}\right)=\I\left(k_{0}\right)_{X}$.
  Thus
  \begin{equation*}
    \nu^{-1}\I\left(k_{0}\right)_{X}
    =\nu^{-1}j^{-1}\I\left(k_{0}\right)
    \cong\left(j\circ \nu\right)^{-1}\I\left(k_{0}\right)
    =\I\left(k_{0}\right)_{X_{\proet}},
  \end{equation*}
  by~\cite[Lemma 3.45]{Bo21},
  where $\I\left(k_{0}\right)_{X_{\proet}}$ is the constant sheaf on $X_{\proet}$.
  Now compute
  \begin{equation*}
  \begin{split}
    \Hom_{\I\left(k_{0}\right)}\left(
      \I\left(k_{0}\right)_{X},
      \R\nu_{*}\R\shHom_{\I\left(\BB_{\pdR}^{\dag}\right)}\left( \cal{F}^{\bullet} , \cal{G}^{\bullet} \right)
    \right)
    &\cong
    \Hom_{\I\left(k_{0}\right)}\left(
      \nu^{-1}\I\left(k_{0}\right)_{X},
      \R\shHom_{\I\left(\BB_{\pdR}^{\dag}\right)}\left( \cal{F}^{\bullet} , \cal{G}^{\bullet} \right)
    \right) \\
    &\cong
    \Hom_{\I\left(k_{0}\right)}\left(
      \I\left(k_{0}\right)_{X_{\proet}},
      \R\shHom_{\I\left(\BB_{\pdR}^{\dag}\right)}\left( \cal{F}^{\bullet} , \cal{G}^{\bullet} \right)
    \right) \\
    &\stackrel{\clubsuit}{\cong}
    \Hom_{\I\left(\BB_{\pdR}^{\dag}\right)}\left(
      \I\left(k_{0}\right)_{X_{\proet}} \widehat{\otimes}_{\I\left(k_{0}\right)}^{\rL} \cal{F}^{\bullet},
      \cal{G}^{\bullet}
    \right) \\
    &\cong
    \Hom_{\I\left(\BB_{\pdR}^{\dag}\right)}\left( \cal{F}^{\bullet} , \cal{G}^{\bullet} \right)
  \end{split}
  \end{equation*}
  We also used~\cite[Lemma 3.11(ii)]{Bo21} for the isomorphism $\clubsuit$.
\end{proof}

\begin{proof}[Proof of Theorem~\ref{thm:pdRdagsolfunctor-fullyfaithful-reconstructionpaper}]
\label{proof:--thm:pdRdagsolfunctor-fullyfaithful-reconstructionpaper}
  For any two $\cal{C}$-complexes $\cal{M}^{\bullet}$ and $\cal{N}^{\bullet}$,
  we know from Theorem~\ref{thm:solfunctorfullyfaithfull-reconstructionpaper-globalversioninternalhom}
  that~(\ref{eq:themapLH--subsubsec:proofofsolfunctorfullyfaithfull-reconstructionpaper-constructingmorphisms})
  is an isomorphism.
  We apply $\Hom_{\I\left(k_{0}\right)}\left( \I\left(k_{0}\right)_{X} , - \right)$ to this isomorphism.
  Invoking
  Lemma~\ref{lem:shHom-extHom--thm:solfunctorfullyfaithfull-reconstructionpaper}
  and~\ref{lem:shHom-extHom2--thm:solfunctorfullyfaithfull-reconstructionpaper},
  we get
  \begin{equation*}
    \Hom_{\I\left(\Dcap\right)}\left( \cal{N}^{\bullet} , \cal{M}^{\bullet} \right) \\
    \isomap\Hom_{\I\left(\BB_{\pdR}^{\dag}\right)}\left(
      \SolB_{\pdR}^{\dag}\left(\cal{M}^{\bullet}\right),
      \SolB_{\pdR}^{\dag}\left(\cal{N}^{\bullet}\right)
    \right).
  \end{equation*}
  Since both $\cal{M}^{\bullet}$ and $\cal{N}^{\bullet}$ are arbitrary $\cal{C}$-complexes,
  this completes the proof of Theorem~\ref{thm:pdRdagsolfunctor-fullyfaithful-reconstructionpaper}.
\end{proof}

\appendix


\section{Miscellaneous}
\label{appendix:reconstructionpaper}


\subsection{Closed symmetric monoidal categories}

Fix a closed symmetric monoidal category $\left( \C,1,\otimes\right)$.

\begin{lem}\label{lem:epi-uniquemoduleaction}
  Let $R\in\C$ denote a monoid and $M$ an $R$-module object.
  Fix an epimorphism $\phi\colon M\to N$. If there exists an
  $R$-module structure on $N$ making
  $\phi$ $R$-linear, then this structure is unique.  
\end{lem}

\begin{proof}
  Consider two $R$-module structures on $N$ with
  actions $\act_{N,1}\colon R\otimes N\to N$, respectively $\act_{N,2}$,
  and units $1_{N,1}\colon R\to N$, respectively $1_{N,2}$. We assume
  that $\phi$ is $R$-linear with respect to both $R$-module structures on $N$.

  Denote the action of $R$ on $M$ by $\act_{M}\colon R\otimes M \to M$.
  The $\phi$-linearities imply
  \begin{equation*}
    \act_{N,1} \circ \left( \id_{R}\otimes\phi\right)
    =\phi\circ\act_{M}, \text{ and }
    \act_{N,2} \circ \left( \id_{R}\otimes\phi\right)
    =\phi\circ\act_{M}.
  \end{equation*}
  But $\id_{R}\otimes\phi$ is an epimorphism
  because $\C$ is closed; $\act_{N,1}=\act_{N,2}$ follows.
  
  Furthermore, the $\phi$-linearities imply
  $1_{N,1}=\phi\circ 1_{M}=1_{N,2}$.
\end{proof}

The following lemma is well-known, thus we omit the proof.

\begin{lem}\label{lem:monoid-on-tensorproduct}
  Let $R,S$ denote two monoids in $\C$
  Then $R\otimes S$ becomes a monoid as follows.
  The multiplication is the composition
  \begin{equation*}
    \left(R \otimes S\right)
    \otimes
    \left(R \otimes S\right)
    \cong
    \left(R \otimes R\right)
    \otimes
    \left(S \otimes S\right)
    \stackrel{\mu_{R}\otimes\mu_{S}}{\longrightarrow}
    R\otimes S
  \end{equation*}
  and the unit is the composition
  \begin{equation*}
    1 \cong 1\otimes 1 \stackrel{1_{R}\otimes1_{S}}{\longrightarrow}
    R\otimes S.
  \end{equation*}
  Here $\mu_{*}$ is the multiplication and
  $1_{*}$ is the unit of $*\in\left\{R,S\right\}$.
\end{lem}

\begin{defn}\label{defn:bimoduleobject}
  Let $R$ and $S$ denote two monoid objects in $\C$.
  An $R$-$S$-bimodule object $M$ is an
  $R\otimes S^{\op}$-module object. Here, the underlying
  $\C$-objects of $S$ and $S^{\op}$ coincide, but the
  multiplication is performed in the reverse order.
  
  A morphism between two $R$-$S$-bimodule objects
  is $R$-$S$-linear if it is a morphism of
  $R\otimes S^{\op}$-module objects.
\end{defn}

The following corollary is a special case of Lemma~\ref{lem:epi-uniquemoduleaction}.

\begin{cor}\label{cor:epi-uniquebimoduleaction}
  Let $R,S\in\C$ denote monoids. $M$ is an $R$-$S$-bimodule object.
  Fix an epimorphism $\phi\colon M\to N$. If there exists an
  $R$-$S$-bimodule structure on $N$ making
  $\phi$ $R$-$S$-linear, then this structure is unique.  
\end{cor}

The following two lemma are easy to verify; we leave the details to the reader.

\begin{lem}\label{lem:monoidobject-lim}
  Suppose $\C$ admits all limits. Consider the tower of monoids
  $\dots\to R_{2}\to R_{1} \to R_{0}$,
  where the maps are multiplicative. Then $R:=\varprojlim_{r\in\NN}R$
  becomes a monoid object:
  \begin{equation*}
    R \otimes R
    \to\varprojlim_{r\in\NN}\left( R_{r} \otimes R_{r} \right)
    \to\varprojlim_{r\in\NN}R_{r}
    =R
  \end{equation*}
  is the multiplication and the unit is
  \begin{equation*}
    1=\varprojlim_{r\in\NN}1
    \to\varprojlim_{r\in\NN}R_{r}
    =R.
  \end{equation*}
\end{lem}

\begin{lem}\label{lem:monoidobject-colim}
  Suppose $\C$ admits all colimits. Consider the tower of monoids
  $S^{0} \to S^{1} \to S^{2} \to\dots$,
  where the maps are multiplicative. Then $S:=\varinjlim_{q\in\NN}S^{q}$
  becomes a monoid object:
  \begin{equation*}
    S \otimes S
    =\varinjlim_{q\in\NN}\left( S^{q}\otimes S^{q}\right)
    \to\varinjlim_{q\in\NN}S^{q}
    =S
  \end{equation*}
  is the multiplication and the unit is
  \begin{equation*}
    1=\varinjlim_{q\in\NN}1
    \to\varinjlim_{q\in\NN}S^{q}
    =S.
  \end{equation*}
\end{lem}

We utilise Lemma~\ref{lem:bimodulestructure-limcolim}
in the proof of Theorem~\ref{thm:OBlaplus-bimodule}.

\begin{lem}\label{lem:bimodulestructure-limcolim}
  Suppose $\C$ admits all limits, and consider the tower
  $\dots\to R_{2}\to R_{1} \to R_{0}$
  of monoids,
  where the maps are multiplicative.
  Fix another monoid $S$ and an object $M$,
  together with $R_{r}$-$S$-bimodule structures on $M$ such that
  \begin{itemize}
    \item[(i)] the following diagrams commute for all $r^{\prime}\geq r$:
    \begin{equation*}
      \begin{tikzcd}
        \left( R_{r} \otimes S^{\op}\right)\otimes M
        \arrow{r} &
        M \\
        \left( R_{r^{\prime}} \otimes S^{\op}\right)\otimes M
        \arrow{r}\arrow{u} &
        M.\arrow[equal]{u}  
      \end{tikzcd}
    \end{equation*}
    Here, the horizontal maps
    denote the bimodule actions. Furthermore,
    \item[(ii)] the units $1\to M$ of all given
     $R_{r}$-$S$-bimodule structures on $M$ coincide.
  \end{itemize}
  Apply Lemma~\ref{lem:monoidobject-lim}
  to turn $R:=\varprojlim_{r\in\NN}R_{r}$ into a monoid.
  Then the composition
  \begin{equation*}
    \left( R\otimes S^{\op}\right)\otimes M
    \to\varprojlim_{r\geq 0}\left(
    R_{r}\otimes S^{\op}
    \right)\otimes M
    \stackrel{\text{(i)}}{\to} M
  \end{equation*}
  defines an $R$-$S$-bimodule action on $M$.
  The unit $1\to M$ is the unit of any of the $R_{r}$-$S$-bimodule structure, which
  is allowed by (ii).
\end{lem}

\begin{proof}
  We denote multiplication maps by $\mu$ and module actions
  by $\act$. Consider
\begin{equation}\label{cd:bimodulestructure-limcolim-associative}
\begin{tikzcd}[scale cd=0.85]
\varprojlim_{r\geq 0}\left(
  \begin{multlined}\left(R_{r}\otimes S^{\op}\right)\otimes \\
  \left(R_{r}\otimes S^{\op}\right)\otimes M
\end{multlined}\right) \arrow[ddd] \arrow[rrr] &&&
\varprojlim_{r\geq 0}
  \left(R_{r}\otimes S^{\op}\right)
  \otimes M
\arrow[ddd] \\ 
&
\begin{multlined}
  \left(R\otimes S^{\op}\right) \\
  \otimes\left(R\otimes S^{\op}\right)\otimes M
\end{multlined}
\arrow[lu] \arrow{d}{\mu\otimes\id} \arrow{r}{\id\otimes a} &
{\left(R\otimes S^{\op}\right)\otimes M} \arrow{d}{\act} \arrow[ru] & \\ 
&
{\left(R\otimes S^{\op}\right)\otimes M} \arrow[ld] \arrow{r}{\act} &
M \arrow[rd, equal] & \\ 
\varprojlim_{r\geq 0}
\left(R_{r}\otimes S^{\op}\right)
\otimes M
\arrow[rrr] &&&
M
\end{tikzcd}
\end{equation}
The diagram is commutative: this is obvious for all rectangles except the
middle one. However, the commutativity of the middle rectangle follows
from the commutativity of the others. We have thus the associativity of the
$R$-$S$-action on $M$.
Next, we consider the diagram
\begin{equation}\label{cd:bimodulestructure-limcolim-preserveunits}
\begin{tikzcd}
&&&
\varprojlim_{r\geq 0}\left( R_{r}\otimes S^{\op} \right)\otimes M
\arrow[lldd, bend left] \\
1\otimes M \arrow[rr] \arrow[rd] \arrow[rrru] &&
\left( R\otimes S^{\op}\right)\otimes M \arrow[ld] \arrow[ru] &
\\
&
M. &&                          
\end{tikzcd}
\end{equation}
The diagram is commutative: this is obvious for all triangles except the
one in the bottom left corner. However, the commutativity of the triangle follows
from the commutativity of the others. We have thus shown that the
$R$-$S$-action on $M$ preserves the unit.
\end{proof}

\begin{lem}\label{lem:bimodulestructure-colim}
  Suppose $\C$ admits all colimits, and consider the tower
  $S^{0}\to S^{1} \to S^{2} \to\dots$
  of monoids, where the maps are multiplicative.
  Fix another monoid $R$ and a tower of object $M^{0}\to M^{1}\to M^{2}\dots$,
  together with $R$-$S^{q}$-bimodule structures on $M^{q}$ such that
  \begin{itemize}
    \item[(i)] the following diagrams commute for all $q^{\prime}\geq q$:
    \begin{equation*}
      \begin{tikzcd}
        \left( R \otimes S^{q^{\prime},\op}\right)\otimes M^{q^{\prime}}
        \arrow{r} &
        M^{q^{\prime}} \\
        \left( R \otimes S^{q,\op}\right)\otimes M^{q}
        \arrow{r}\arrow{u} &
        M^{q}.\arrow{u}  
      \end{tikzcd}
    \end{equation*}
    Here, the horizontal maps
    denote the bimodule actions. Furthermore,
    \item[(ii)] the following diagram commutes for all $q^{\prime}\geq q$:
    \begin{equation*}
      \begin{tikzcd}
        1
        \arrow{r} &
        M^{q^{\prime}} \\
        1
        \arrow{r}\arrow{u} &
        M^{q}.\arrow{u}  
      \end{tikzcd}
    \end{equation*}    
    Here, the horizontal maps are the units maps.
  \end{itemize}
  Apply Lemma~\ref{lem:monoidobject-colim}
  to turn $S:=\varinjlim_{q\geq0}S^{q}$ into a monoid.
  Then the composition
  \begin{equation*}
    \left( R\otimes S^{\op}\right)\otimes M
    =\varinjlim_{q\geq0}\left(
    R\otimes S^{q,\op}
    \right)\otimes M^{q}
    \stackrel{\text{(i)}}{\to} \varinjlim_{q\geq0} M^{q}=0
  \end{equation*}
  defines an $R$-$S$-bimodule action on $M$.
  Thanks to (ii), we can define the
  unit $1\to M$ is the colimit of the units of the $R$-$S^{q}$-bimodule structures on $M^{q}$
\end{lem}

\begin{proof}
  This is easy.
\end{proof}


\subsection{Projective objets in $\IndBan_{\I(F)}$}

Fix a field $F=R$, complete with respect to nontrivial non-Archimedean valuation.
Furthermore, $\pi\in F$ denotes a pseudo-uniformiser, that is $0<|\pi|<1$.
We refer the reader to~\cite[\S\ref*{ch:functional-analysis}]{WiersigPeriods}
for the definitions of the categories $\Ban_{F}$, $\IndBan_{F}$, and $\IndBan_{\I(F)}$.

\begin{lem}\label{lem:Tate-algebra-internally-projective-BanF-reconstructionpaper}
  The following functor is exact for any $d,q\in\NN$:
  \begin{equation*}
    \intHom_{F}\left( F\left\< \frac{Z_{1},\dots,Z_{d}}{\pi^{q}}\right\>,-\right)
    \colon\Ban_{F}\to\Ban_{F}.
  \end{equation*}
\end{lem}

\begin{proof}
  Thanks to the open mapping theorem, it suffices to check that
  $\Hom_{F}\left( F\left\< \frac{Z_{1},\dots,Z_{d}}{\pi^{q}}\right\>,-\right)$
  is exact, that is $F\left\< \frac{Z_{1},\dots,Z_{d}}{\pi^{q}}\right\>$ is projective.
  But Lemma~\ref{lem:infiniteTatealgebrasascontractingcoproductandc0}
  gives an isomorphism
  \begin{equation*}
    F\left\< \frac{Z_{1},\dots,Z_{d}}{\pi^{q}}\right\>
    \cong{\coprod_{\alpha\in\NN^{d}}}^{\leq1}F,
  \end{equation*}
  where the right-hand side denotes the \emph{contracting coproduct};
  cf.~\cite[\S A.4]{BBK} for details.
  Following the proof of \emph{loc. cit.} Lemma A.38,
  one checks that ${\coprod_{\alpha\in\NN^{d}}}^{\leq1}F$ is projective, as desired.
\end{proof}

\begin{lem}\label{lem:Tate-algebra-internally-projective-IndBanF-reconstructionpaper}
  The following functor is exact for any $d,q\in\NN$:
  \begin{equation*}
    \intHom_{F}\left( F\left\< \frac{Z_{1},\dots,Z_{d}}{\pi^{q}}\right\>,-\right)
    \colon\IndBan_{F}\to\IndBan_{F}.
  \end{equation*}
\end{lem}

\begin{proof}
  The functor is right exact since $\IndBan_{F}$ is closed,
  cf. Lemma~\ref{lem:indBanFcirc-structureresult}.
  It remains to check that
  $\intHom_{F}\left( F\left\< \frac{Z_{1},\dots,Z_{d}}{\pi^{q}}\right\>,-\right)$
  preserves strict epimorphisms.
  By Lemma~\ref{lem:messmse-indcat},
  it suffices to check that it
  preserves strict epimorphisms of $F$-Banach spaces.
  Now apply Lemma~\ref{lem:Tate-algebra-internally-projective-BanF-reconstructionpaper}.
\end{proof}

Finally, we pass to the left heart.

\begin{lem}\label{lem:IHom-Tatealgebra-is-HomITatealgebra-reconstructionpaper}
  For any $F$-Banach space $V$, the following canonical morphism
  is an isomorphism:
  \begin{equation*}
    \I\left(\intHom_{F}\left(F\left\< \frac{Z_{1},\dots,Z_{d}}{\pi^{q}}\right\> , V \right)\right)
    \isomap
    \intHom_{\I(F)}\left( \I\left( F\left\< \frac{Z_{1},\dots,Z_{d}}{\pi^{q}}\right\>\right) , \I\left(V\right) \right).
  \end{equation*}
\end{lem}

\begin{proof}
  This follows from Lemma~\ref{lem:Tate-algebra-internally-projective-IndBanF-reconstructionpaper}
  and the definition of $\intHom_{\I(F)}$.
\end{proof}

\begin{lem}\label{lem:Tate-algebra-internally-projective-LHIndBanF-reconstructionpaper}
  The following functor is exact for any $d,q\in\NN$:
  \begin{equation*}
    \intHom_{\I(F)}\left( \I\left(F\left\< \frac{Z_{1},\dots,Z_{d}}{\pi^{q}}\right\>\right),-\right)
    \colon\IndBan_{\I(F)}\to\IndBan_{\I(F)}.
  \end{equation*}
\end{lem}

\begin{proof}
  Since $\IndBan_{F}$ is elementary, cf. Lemma~\ref{lem:indBanFcirc-structureresult},
  we have a small strictly generating set $\cal{G}$ of \emph{tiny} projective
  objects, cf.~\cite[Definition 2.1.1]{Sch99}. The proof of \emph{loc. cit.} Proposition 2.1.12
  shows that $\left\{ \I(P)\colon P\in\cal{G}\right\}$ is a strictly generating set of tiny
  projective objects in $\IndBan_{\I\left( F\right)}$. By~\cite[Proposition 2.1.14]{Sch99},
  we may now check that the functor
  \begin{equation*}
    \Hom_{\I(F)}\left( \I(P) ,
    \intHom_{\I(F)}\left( \I\left( F\left\< \frac{Z_{1},\dots,Z_{d}}{\pi^{q}}\right\>\right),-\right)
    \right)
  \end{equation*}
  is exact for all $P\in\cal{G}$. By the tensor-hom adjunction, it coincides with
  \begin{equation*}
    \Hom_{\I(F)}\left( \I(P) \widehat{\otimes}_{\I(F)}\I\left( F\left\< \frac{Z_{1},\dots,Z_{d}}{\pi^{q}}\right\>\right),-\right) \\
    \stackrel{\text{\ref{lem:Ifunctor-stronglymonoidal-reconstructionpaper}}}{\cong}
    \Hom_{\I(F)}\left( \I\left(P \widehat{\otimes}_{F} F\left\< \frac{Z_{1},\dots,Z_{d}}{\pi^{q}}\right\>\right),-\right).
  \end{equation*}
  It suffices to show that $\I\left(P \widehat{\otimes}_{F} F\left\< \frac{Z_{1},\dots,Z_{d}}{\pi^{q}}\right\>\right)$
  is projective. By~\cite[Proposition 1.3.24(a)]{Sch99}, this is the case precisely when
  $P \widehat{\otimes}_{F} F\left\< \frac{Z_{1},\dots,Z_{d}}{\pi^{q}}\right\>$ is projective. This follows because
  \begin{equation*}
    \Hom_{F}\left( P \widehat{\otimes}_{F} F\left\< \frac{Z_{1},\dots,Z_{d}}{\pi^{q}}\right\>,-\right)
    =\Hom_{F}\left( P ,\intHom_{F}\left(F\left\< \frac{Z_{1},\dots,Z_{d}}{\pi^{q}}\right\>,-\right)\right),
  \end{equation*}
  together with the projectivity of $P$ and Lemma~\ref{lem:Tate-algebra-internally-projective-IndBanF-reconstructionpaper}.
\end{proof}

Finally, we study a particular example of the monoidal operation on $\IndBan_{\I\left(F\right)}$.

\begin{lem}\label{lem:constantsheaf-Tatealgebra-LH-exact}
  The following functor is exact for any $d,q\in\NN$:
  \begin{equation*}
    \I\left(F\left\< \frac{Z_{1},\dots,Z_{d}}{\pi^{q}}\right\>\right)
      \widehat{\otimes}_{\I(F)}-
    \colon\IndBan_{\I(F)}\to\IndBan_{\I(F)}.
  \end{equation*}  
\end{lem}

\begin{proof}
  
  $F\left\< \frac{Z_{1},\dots,Z_{d}}{\pi^{q}}\right\>\otimes_{F}-\colon\IndBan_{F}\to\IndBan_{F}$
  is strongly exact, cf. Lemma~\ref{lem:restrictedpowerseries-exact-strongerconvergence}.
  Therefore, we can follow the arguments as in the~\cite[proof of Corollary 5.36]{Bo21}.
\end{proof}


\subsection{Some constant sheaves}

We fix again a field $F=R$, complete with respect to nontrivial non-Archimedean valuation.
Let $X$ be a site and recall the category $\Sh\left(X,\IndBan_{\I(F)}\right)$ as
in~\cite[\S\ref{subsec:sheavesindbanspaces-reconstructionpaper}]{WiersigPeriods}.
We study some properties of constant sheaves in this setting.

\begin{defn}\label{defn:constantsheaf}
  Given any elementary quasi-abelian category $\E$ and an object
  $E\in\E$, its \emph{constant sheaf $E_{X}$ on $X$} is the sheafification
  of the constant presheaf $X\ni U\mapsto E$.
\end{defn}

\begin{lem}\label{lem:I-commutes-constantsheaf-recpaper}
  Given any $k$-ind-Banach space $V$, we have a functorial isomorphism
  $\I(V)_{X}\cong\I\left(V_{X}\right)$.
\end{lem}

\begin{proof}
  This is because sheafification commutes with $\I$,
  cf.~\cite[Lemma 3.18]{Bo21}.
\end{proof}

See~\cite[\S 3.7]{Bo21}
for the definition of inverse images sheaves in our setting.

\begin{lem}\label{lem:constantsheaf-commuteswith-inverseimage-recpaper}
  Given a morphism of sites $f\colon X\to Y$
  and an object $E\in\IndBan_{\I(F)}$,
  we have a functorial isomorphism
  $E_{X}\cong f^{-1}E_{Y}$.
\end{lem}

\begin{proof}
  This follows from the canonical isomorphism
  $f^{-1}\left(\cal{F}^{a}\right)\cong\left(f_{\text{pre}}^{-1}\cal{F}\right)^{a}$  
  in the~\cite[proof of Lemma 3.40]{Bo21} where we set $\cal{F}$ to be the
  constant presheaf on $Y$ with value $E$.
  We refer the reader to \emph{loc. cit.} for details regarding the notation.
\end{proof}

\begin{lem}\label{lem:constantsheaf-commutes-compeltedtensorproduct-arbitrarybase-recpaper}
  Fix a site $X$, a monoid $R\in\IndBan_{\I(F)}$, a right $R$-module object $M$,
  and a left $R$-module object $N$. Then we have the canonical isomorphism
  \begin{equation*}
    M_{X}\widehat{\otimes}_{R}N_{X}\cong\left(M\widehat{\otimes}_{R}N\right)_{X}
  \end{equation*}
\end{lem}

\begin{proof}
  From the definitions of both $\widehat{\otimes}_{R}$, cf.~\cite[Definition 2.2]{BBK},
  and the exactness of $E\mapsto E_{X}$,
  which follows from~\cite[Lemma 3.15(iii)]{Bo21},
  one sees that it suffices to check that $E\mapsto E_{X}$
  is strongly monoidal with respect to $\widehat{\otimes}_{\I(F)}$.
  This is a special of the more general fact that
  the inverse image functor is strongly monoidal, cf.~\cite[Lemma 3.40]{Bo21}.
\end{proof}

In the following lemma, we use that the functor sending an object $E\in\IndBan_{\I(F)}$
to the constant sheaf $E_{X}$ is exact, which follows from~\cite[Lemma 3.15(iii)]{Bo21}.
This implies that $E\mapsto E_{X}$ extends to a functor between the derived categories.

\begin{lem}\label{lem:derivedpullbackcommuteswithspecifictensorproduct-recpaper}
  Given a morphism of sites $f\colon X\to Y$, a monoid $R\in\IndBan_{\I(F)}$,
  a complex of $R$-module objects $E^{\bullet}$, and a sheaf $\cal{R}$ of $R$-algebras, then we have
  a functorial and canonical isomorphism
  \begin{equation*}
    f^{-1}\left(\cal{R}\widehat{\otimes}_{R}^{\rL}E_{Y}^{\bullet}\right)
    \cong \left(f^{-1}\cal{R}\right)\widehat{\otimes}_{R}^{\rL}E_{X}
  \end{equation*}
  of sheaves of $f^{-1}\cal{R}$-module objects.
\end{lem}

\begin{proof}
  Lemma~\ref{lem:derivedpullbackcommuteswithspecifictensorproduct-recpaper}
  is implied by the following computation and Yoneda's lemma:
  \begin{align*}
    \Hom_{f^{-1}\cal{R}}\left( f^{-1}\left(\cal{R}\widehat{\otimes}_{R}^{\rL}E_{Y}^{\bullet}\right) , - \right)
    &=\Hom_{\cal{R}}\left( \cal{R}\widehat{\otimes}_{R}^{\rL}E_{Y}^{\bullet} , \R f_{*}- \right) \\
    &=\Hom_{R}\left( E_{Y}^{\bullet} , \R f_{*}- \right) \\
    &=\Hom_{R}\left( E^{\bullet} , \R\Gamma\left(Y,\R f_{*}-\right) \right) \\
    &\stackrel{\text{($*$)}}{=}\Hom_{R}\left( E^{\bullet} , \R\Gamma\left(X,-\right) \right) \\
    &=\Hom_{R}\left( E_{X}^{\bullet} , - \right) \\
    &=\Hom_{f^{-1}\cal{R}}\left( \left(f^{-1}\cal{R}\right)\widehat{\otimes}_{R}^{\rL}E_{X}^{\bullet} , - \right).
  \end{align*}
  See~\cite[Lemma 3.45]{Bo21}
  for ($*$).
  We also used the adjunctions as in \emph{loc. cit.} Proposition 3.35 and Lemma 3.41,
  which applied because $E\mapsto E_{X}$ is a special case of an inverse image sheaf functor.
\end{proof}


\subsection{Internal homomorphisms}

We fix again the notation as in~\cite[\S\ref{subsec:sheavesindbanspaces-reconstructionpaper}]{WiersigPeriods}.

\begin{lem}\label{lem:sheafinternalhom-over-monoid}
  Fix a monoid $\cal{R}$ in $\Sh\left(X,\IndBan_{\I(F)}\right)$ and two $\cal{R}$-module objects
  $\cal{F}$ and $\cal{G}$. Then $\intHom_{\cal{R}}\left(\cal{F},\cal{G}\right)$ is the
  kernel of the map
  \begin{equation*}
    h\colon
    \prod_{U\in X}\intHom_{\cal{R}(U)}\left(\cal{F}(U),\cal{G}(U)\right)
    \to
    \prod_{\substack{V,W\in X \\ W\to V}}\intHom_{\cal{R}(V)}\left(\cal{F}(V),\cal{G}(W)\right)
  \end{equation*}
  defined by setting
  \begin{equation*}
    \sigma_{VW}\circ h
    =\intHom_{\cal{R}(V)}\left(\id_{\cal{F}(V)},r_{WV}\right)\circ\sigma_{V}
      -\intHom_{\cal{R}(V)}\left(r_{WV},\id_{\cal{G}(W)}\right)\circ\sigma_{W}
  \end{equation*}
  for all $W\to V$ in $X$. Here, $\sigma_{U}$
  denotes the projection onto $\intHom_{\cal{R}(U)}\left(\cal{F}(U),\cal{G}(U)\right)$
  and $\sigma_{VW}$ denotes the projection onto
  $\intHom_{\cal{R}(V)}\left(\cal{F}(V),\cal{G}(W)\right)$.
  Also, $r_{WV}$ refers to both restriction maps
  $\cal{F}(W)\to\cal{F}(V)$ and $\cal{G}(W)\to\cal{G}(V)$.
\end{lem}

\begin{proof}
  This follows directly from~\cite[Definition 2.2.13]{Sch99}
  and~\cite[Definition 2.2]{BBK}.
\end{proof}

\bibliography{RecThm}
\bibliographystyle{amsplain}  

\end{document}